\documentclass[reqno,10pt]{amsart}
\usepackage{amsmath,mathrsfs}
\usepackage{amssymb, amsmath}
\usepackage{graphicx}
\usepackage{pdfsync}
\usepackage{color}
\usepackage[colorlinks,
            linkcolor=red,
            anchorcolor=blue,
            citecolor=green
            ]{hyperref}
\usepackage{colonequals}
\def\cs#1#2{\big(#1\big)^{\frac{1}{2}}\big(#2\big)^{\frac{1}{2}}}

\def\inte#1{
\displaystyle\mathop{#1\kern0pt}^\circ }

\let\pa=\partial

\let\f=\frac

\def\pa{\partial}

\def\virgp{\raise 2pt\hbox{,}}
\def\cdotpv{\raise 2pt\hbox{;}}

\def\C{\mathop{\mathbb C\kern 0pt}\nolimits}
\def\DD{\mathop{\mathbb D\kern 0pt}\nolimits}
\def\EE{\mathop{{\mathbb E \kern 0pt}}\nolimits}
\def\K{\mathop{\mathbb K\kern 0pt}\nolimits}
\def\N{\mathop{\mathbb N\kern 0pt}\nolimits}
\def\Q{\mathop{\mathbb Q\kern 0pt}\nolimits}
\def\R{\mathop{\mathbb R\kern 0pt}\nolimits}
\def\SS{\mathop{\mathbb S\kern 0pt}\nolimits}
\def\ZZ{\mathop{\mathbb Z\kern 0pt}\nolimits}
\def\TT{\mathop{\mathbb T\kern 0pt}\nolimits}
\def\P{\mathop{\mathbb P\kern 0pt}\nolimits}

\DeclareMathOperator*{\esssup}{ess\,sup}

\newcommand{\beq}{\begin{equation}}
\newcommand{\eeq}{\end{equation}}
\newcommand{\ben}{\begin{eqnarray}}
\newcommand{\een}{\end{eqnarray}}
\newcommand{\beno}{\begin{eqnarray*}}
\newcommand{\eeno}{\end{eqnarray*}}
\newtheorem{thm}{Theorem}[section]
\newtheorem{lem}{Lemma}[section]

\newtheorem{prop}{Proposition}[section]

\theoremstyle{definition}

\newtheorem{rmk}{Remark}[section]

\numberwithin{equation}{section}
\begin{document}
\title[BBE limit]
{Quantitative global well-posedness of Boltzmann-Bose-Einstein equation and incompressible Navier-Stokes-Fourier limit}
% stability and hydrodynamic limit of

\author[L.-B. He, N. Jiang and Y.-L. Zhou]{Ling-Bing He, Ning Jiang and Yu-long Zhou}
\address[L.-B. He]{Department of Mathematical Sciences, Tsinghua University\\
Beijing, 100084,  P. R.  China.} \email{hlb@tsinghua.edu.cn}
\address[N. Jiang]{School of Mathematics and Statistics, Wuhan University, Wuhan, 430072, P. R. China} \email{njiang@whu.edu.cn}
\address[Y.-L. Zhou]{School of Mathematics, Sun Yat-Sen University, Guangzhou, 510275, P. R.  China.} \email{zhouyulong@mail.sysu.edu.cn}

\begin{abstract}
In the diffusive scaling and in the whole space, we prove the global well-posedness of the scaled Boltzmann-Bose-Einstein (briefly, BBE) equation with high temperature in the low regularity space $H^2_xL^2$. In particular, we quantify the fluctuation around the Bose-Einstein equilibrium $\mathcal{M}_{\lambda,T}(v)$ with respect to the parameters $\lambda$ and temperature $T$. Furthermore, the estimate for the diffusively scaled BBE equation is uniform to the Knudsen number $\epsilon$. As a consequence, we rigorously justify the hydrodynamic limit to the incompressible Navier-Stokes-Fourier equations. This is the first rigorous fluid limit result for BBE.
\end{abstract}

\maketitle
\markright{Quantum Boltzmann equation}

\setcounter{tocdepth}{1}
\tableofcontents

%%%%%%%%%%%%%%

%\noindent {\sl Keywords:} {inhomogeneous  Boltzmann
%equation, long-range interactions, asymptotic analysis, %spectral gap, global error.}

%\vskip 0.2cm

\noindent {\sl AMS Subject Classification (2020):} {35Q20, 82C40.}

%%%%%%%%%%%%%%%%%%%%%%%%%%%%%%%%%%%%%%%%%%%%%%
%%%%%%%%%%%%%%%%%%%%%%%%%%%%%%%%%%%%%%%%%%

 %

\section{Introduction}

Quantum Boltzmann equations were proposed to describe the time evolution of a dilute system of weakly interacting bosons or fermions, which obeys the Bose-Einstein or Fermi-Dirac statistics respectively. Consequently these equations are named Boltzmann-Fermi-Dirac (briefly, BFD) or Boltzmann-Bose-Einstein (briefly, BBE) equations, respectively. The derivations of such equation date back to as early as 1920s by Nordheim \cite{nordhiem1928kinetic} and 1933 by Uehling-Uhlenbeck \cite{uehling1933transport}.
Consequently, the quantum Boltzmann equations are also called Boltzmann-Nordheim equations or Uehling-Uhlenbeck equations in the literature.
Later on, further developments were made by Erd{\H{o}}s-Salmhofer-Yau \cite{erdHos2004quantum}, Benedetto-Castella-Esposito-Pulvirenti \cite{benedetto2004some}, \cite{benedetto2005weak},  \cite{benedetto2006some} and a short review \cite{benedetto2007short},  Lukkarinen-Spohn \cite{lukkarinen2009not}.
One can refer the classical book \cite{chapman1990mathematical} for physical backgrounds.

When the quantum effects are not considered, the evolution of dilute gas particles are governed by classical Boltzmann equations:
\beno
\partial _t f +  v \cdot \nabla_{x} f=Q_{B}(f)\,,
\eeno
where the Boltzmann collision term
\beno
Q_{B}(f)(v) = \int_{{\mathbb S}^2 \times {\mathbb R}^3} B(v- v_{*},\sigma)
 \big\{f(v_{*}^{\prime}) f(v^{\prime}) - f(v_{*}) f(v) \big\}
\mathrm{d}\sigma \mathrm{d}v_{*}\,.
\eeno
In the above equation, the unknown non-negative function $f\equiv f(t,x,v)$ is the so-called number density of the particles. It describes the evolution of the particles at time $t\geq 0$, at position $x\in \Omega$, with velocity $v\in \mathbb{R}^3$. The domain $\Omega$ could be the whole space, torus and other domains with boundaries. In this paper, we assume $\Omega=\mathbb{R}^3$. In this setting, we assume all the particles have the same mass, with velocities $v$ and $v_*$ before the collisions, and $v^\prime$ and $v_{*}^{\prime}$ after the collisions. Here we only consider the elastic collisions which conserve the momentum and kinetic energy, i.e.
\begin{equation}\label{conservation}
  v+v_*= v^\prime + v_{*}^{\prime}\,,\quad\! |v|^2+|v_*|^2= |v^\prime|^2 + |v_{*}^{\prime}|^2\,.
\end{equation}
The above conservation laws include four equations, while $v^\prime$ and $v_*^\prime$ are six unknowns. So we need parameters in two dimensional manifold $\mathbb{S}^2$ to represent $v^\prime$ and $v_*^\prime$ in terms of $v$ and $v_*$:
\beno
v^{\prime}=\frac{v+v_{*}}{2}+\frac{|v-v_{*}|}{2}\sigma, \quad  v_{*}^{\prime}=\frac{v+v_{*}}{2}-\frac{|v-v_{*}|}{2}\sigma\,,
\eeno
where $\sigma \in \mathbb{S}^2$. In fact, it is easy to see $\sigma=\frac{v^\prime-v_*^\prime}{|v^\prime-v_*^\prime|}\,.$ For the general introduction to Boltzmann equation, see the standard references \cite{CIP1994, Sone-2002, Sone-2007}.

\subsection{Quantum Boltzmann equations}

Now we introduce the so-called quantum Boltzmann equations:
\ben \label{quantum-Boltzmann-UU}
\partial_t F +  v \cdot \nabla_{x} F=Q_{\Phi, \hbar}(F,F), ~~t > 0, x \in \mathbb{R}^{3}, v \in \R^3 ; \quad
F|_{t=0}(x,v) = F_{0}(x,v).
\een
Here $F(t,x,v)\geq 0$ is the density function of particles with velocity
$v\in\R^3$ at time $t\geq 0$ in position $x \in \mathbb{R}^{3}$.
The  quantum Boltzmann collision operator $Q_{\Phi, \hbar}$ acting only on velocity variable
$v$ is defined by
\ben \label{U-U-operator}
Q_{\Phi, \hbar}(g,h)(v) \colonequals   \int_{{\mathbb S}^2 \times {\mathbb R}^3} B_{\Phi,\hbar}(v- v_{*},\sigma)
 \mathrm{D}\big(g_{*}^{\prime} h^{\prime}(1 +  \delta\hbar^{3}g_{*})(1 +  \delta\hbar^{3}h)\big)
\mathrm{d}\sigma \mathrm{d}v_{*},
\een
where according to \cite{erdHos2004quantum} and \cite{benedetto2005weak}, the quantum Boltzmann collision kernel $B_{\Phi,\hbar}(v- v_{*},\sigma)$ has the following form
\ben \label{scaling-Boltzmann-kernel}
B_{\Phi,\hbar}(v- v_{*},\sigma) \colonequals   \hbar^{-4} |v-v_{*}| \big(
\hat{\Phi} (\hbar^{-1} |v-v^{\prime}|)
+ \hat{\Phi} (\hbar^{-1} |v-v^{\prime}_{*}|)
\big)^{2}\,.
\een
Here the radial function $\hat{\Phi}(|\xi|) \colonequals   \hat{\Phi}(\xi) = \int_{\mathbb{R}^{3}} e^{-\mathrm{i}x\cdot\xi}\Phi(x)\mathrm{d}x$ is the Fourier transform of a radial potential function  $\Phi(x)$. Furthermore, in \eqref{U-U-operator} $\hbar$ is the Plank constant, and $\delta= 1$ or $-1$. Specifically, $\delta=1$ corresponds to Bose-Einstein statistics, while $\delta=-1$ corresponds to Fermi-Dirac statistics.

In \eqref{U-U-operator} and the rest of the article,
we use the convenient  shorthands $h=h(v)$, $g_*=g(v_*)$,
$h'=h(v')$, $g'_*=g(v'_*)$ where $v'$, $v_*'$ are given by
\ben\label{v-prime-v-prime-star}
v'=\frac{v+v_{*}}{2}+\frac{|v-v_{*}|}{2}\sigma, \quad  v'_{*}=\frac{v+v_{*}}{2}-\frac{|v-v_{*}|}{2}\sigma, \quad \sigma\in\SS^{2}.
\een
Now we explain the notation $\mathrm{D}(\cdot)$ in \eqref{U-U-operator}. For $n=1$ or $n=2$, we denote
\ben \label{shorthand-D}
\mathrm{D}^{n}(f(v,v_{*},v^{\prime},v^{\prime}_{*})) \colonequals \left(f(v,v_{*},v^{\prime},v^{\prime}_{*}) - f(v^{\prime},v^{\prime}_{*},v,v_{*})\right)^{n}.
\een
If $n=1$, we write
$\mathrm{D}(\cdot) = \mathrm{D}^{1} (\cdot)$.
The term $\mathrm{D}$ is interpreted as ``difference'' before and after collision. In particular, in \eqref{U-U-operator},
\begin{equation}
  \mathrm{D}\big(g_{*}^{\prime} h^{\prime}(1 +  \delta\hbar^{3}g_{*})(1 +  \delta\hbar^{3}h)\big) = g_{*}^{\prime} h^{\prime}(1 +  \delta \hbar^{3}g_{*})(1 +  \delta\hbar^{3}h) - g_{*} h(1 +  \delta\hbar^{3}g_{*}^{\prime})(1 + \delta \hbar^{3}h^{\prime})\,.
\end{equation}

By the following scaling
\ben \label{scaling-tranform}
\tilde{F}(t,x,v) = \hbar^{3}F(\hbar^{3}t,x,\hbar^{-3}v), \quad \phi(|x|) = \hbar^{4}\Phi(\hbar^{4}|x|),
\een
we can normalize the Plank constant $\hbar$. Indeed, it is easy to check $F$ is a solution to \eqref{quantum-Boltzmann-UU} if and only if $\tilde{F}$ is a solution of the following normalized equation
\ben \label{quantum-Boltzmann-UU-scaling}
\partial _t F +  v \cdot \nabla_{x} F=Q_{\phi}(F,F), ~~t > 0, x \in \mathbb{R}^{3}, v \in \R^3 ; \quad
F|_{t=0}(x,v) = \hbar^{3}F_{0}(x,\hbar^{-3}v),
\een
where the operator $Q_{\phi}$ is defined by
\ben \label{U-U-operator-scaling}
Q_{\phi}(g,h)(v) \colonequals   \int_{{\mathbb S}^2 \times {\mathbb R}^3} B_{\phi}(v- v_{*},\sigma)
\mathrm{D}\big(g_{*}^{\prime} h^{\prime}(1 + \delta g_{*})(1 + \delta h)\big)
\mathrm{d}\sigma \mathrm{d}v_{*}
\een
with the kernel $B_{\phi}(v- v_{*},\sigma)$ given by
\ben \label{scaling-Boltzmann-kernel-scaling}
B_{\phi}(v- v_{*},\sigma) \colonequals   |v-v_{*}| \big(
\hat{\phi} (|v-v^{\prime}|)
+ \hat{\phi} (|v-v^{\prime}_{*}|)
\big)^{2}.
\een

In the current paper, we will take the potential $\phi(x) = \f{1}{2} \delta(x)$. Then the kernel in \eqref{scaling-Boltzmann-kernel-scaling} reduces to that of {\em hard sphere} model
\ben \label{hard-sphere-kernel}
B(v- v_{*},\sigma) \colonequals   |v-v_{*}|
\een
and the collision operator reduces to
\ben \label{U-U-operator-scaling-hard-sphere}
Q(g,h)(v) \colonequals   \int_{{\mathbb S}^2 \times {\mathbb R}^3} B(v- v_{*},\sigma)
\mathrm{D}\big(g_{*}^{\prime} h^{\prime}(1 + \delta g_{*})(1 + \delta h)\big)
\mathrm{d}\sigma \mathrm{d}v_{*}\,.
\een
As introduced before, when $\delta=1$ in the collision kernel \eqref{U-U-operator-scaling-hard-sphere}, the corresponding quantum Boltzmann equation \eqref{quantum-Boltzmann-UU-scaling} is called Boltzmann-Bose-Einstein equation, briefly, BBE. When $\delta=-1$, the equation \eqref{quantum-Boltzmann-UU-scaling} is called Boltzmann-Fermi-Dirac equation, briefly, BFD. In this paper, we focus on BBE equation, i.e. $\delta=1$. More specifically, we study the fluid dynamics from BBE to incompressible Navier-Stokes-Fourier system. The same limit from BFD (i.e. for $\delta=-1$ case) was treated in \cite{jiang2021incompressible}.

\subsection{Well-posedness of quantum Boltzmann equations}

As mentioned above, quantum Boltzmann equations include BFD and BBE equations. The first mathematical question is the well-posedness of these equations, in the framework of corresponding functional spaces of regularity, such as weak solutions, smooth solutions, or between.  Compared to the extensive studies on the classical Boltzmann equations, much less has been done on quantum Boltzmann equations. Mathematically, BFD is relatively easier. We first review the studies on the BFD in the past three decades.

For mathematical theory of well-posedness of BFD equation, early results were obtained by Dolbeault \cite{Dolbeault} and Lions \cite{Lions}. They studied the global weak existence of solutions in mild or distributional sense for the whole space $\R^3$ under some assumptions on the collision kernel. Furthermore, Dolbeault \cite{Dolbeault} obtained that the solution of BFD equation converges to the solution of the Boltzmann equation as $\delta \to 0$ (the $\delta$ appears in \eqref{U-U-operator}) for very special bounded collision kernel case. Allemand \cite{Allemand} tried to extend the results of \cite{Dolbeault} to bounded domains with specular reflection boundary condition for integrable collision kernels. Alexandre \cite{Alexandre} obtained another kind of weak solutions satisfying the entropy inequality, the so-called $H$-solutions. Up to now, the best results on the global weak solutions to the BFD equation belong to Lu.  For general initial data, Lu \cite{Lu-2001, Lu-Wennberg-2003, Lu-2008} studied the global existence and stability of weak solutions on torus for very soft potential with a weak angular cutoff. More results on weak solutions are referred to \cite{EMV, EMV-2004}. We remark that in \cite{jiang2021incompressible}, the authors derived uniform in $\epsilon\in (0,1)$ energy estimate in the incompressible regime. This also gave a global in time smooth solution near global Fermi-Dirac distribution.

For the more difficult BBE equation which has more interesting physical phenomena such as the famous Bose-Einstein condensation, the mathematical research is even less. In a series papers, Lu \cite{lu2004isotropic, lu2005boltzmann, Lu-2011JSP, Lu-2013JSP, Lu-2014JSP, Lu-2016JSP, Lu-2018KRM, Cai-Lu-2019JSP} made major contributions on the systematic analytical studies of weak solutions of (mainly homogeneous) BBE, including existence of weak solutions, convergence to equilibrium (i.e. global Bose-Einstein distribution), long time behavior, etc. In particular, the condensation at low temperature was investigated. The focus of the current paper will be on the classical solution of BBE near the equilibrium (i.e. global Bose-Einstein distribution), but at high temperature. We will explain the relation of equilibrium of BBE with temperature later in this section.

\subsection{Hydrodynamic limits of quantum Boltzmann equations}
In the other direction, the hydrodynamic limits from kinetic equations to fluid equations have been very active in recent decades. One of the important features of kinetic equations is their connection to the fluid equations. The so-called hydrodynamic limits are the process that the Knudsen number $\epsilon$ goes to zero. Here $\epsilon>0$ is the Knudsen number which is the dimensionless quantity defined as the ratio of mean free path and macroscopic length scale. Depending on the physically scalings, different fluid equations (incompressible or compressible Navier-Stokes, Euler, etc.) can be derived from kinetic equations.

Bardos and Ukai \cite{Bardos-Ukai} proved the global existence of classical solution to diffusively scaled Boltzmann equation uniformly in $0<\epsilon<1$ for hard potential with cutoff collision kernel. Consequently they justified the limit to incompressible Navier-Stokes equations with small initial data. By employing semigroup approach, Briant \cite{Briant-2015} also proved the same limit on the torus for hard cutoff potential, in particular, with convergence rate. Jiang, Xu and Zhao \cite{Jiang-Xu-Zhao-2018} proved again the same limit for a more general class of collision kernel. Starting from the solutions to the limiting fluid equations, Caflisch \cite{Caflisch} and Nishida \cite{Nishida-1978} proved the compressible Euler limit from the Boltzmann equation in the context of classical solution by the Hilbert expansion, and analytic solutions, respectively. Caflisch's approach was applied to the acoustic limit by Guo, Jang and Jiang \cite{GJJ-2009, GJJ-2010, JJ-2009} by combining with nonlinear energy method. We also mention some more results using Hilbert expansions \cite{Guo-CPAM2006, Jiang-Xiong-2015}.

For the fluid limits of BFD equation, Zakrevskiy \cite{Zakrevskiy, Zakrevskiy-AA} formally derived the compressible Euler and Navier-Stokes limits and incompressible Navier-Stokes limits. We also mention that, Filbet, Hu and Jin \cite{FHJ} introduced a new scheme for quantum Boltzmann equation to capture the Euler limit by numerical computations. In \cite{jiang2021incompressible}, Jiang-Xiong-Zhou proved the global existence of classical solutions near equilibrium and in addition, they also obtained the uniform in $\epsilon$ energy estimates. As a consequence, the incompressible Navier-Stokes-Fourier limit from BFD equation was justified. Recently, Jiang-Zhou \cite{Jiang-Zhou-2021} studied the compressible Euler limit from BFD equation using Hilbert expansion method.

The main theme of this paper is on the incompressible Navier-Stokes-Fourier limit from BBE equation. The key feature is to obtain a global in time and uniform in Knudsen number $\epsilon$ estimate of the scaled BBE equation. However, this can be only achieved in the case of high temperature at this stage. We start from the scaled BBE in the diffusive scaling, from which the incompressible Navier-Stokes-Fourier system can be derived. For the high temperature case, the formal derivation is similar to that in \cite{Zakrevskiy} for BFD and \cite{BGL-1991JSP} for classical Boltzmann equations. We emphasize that our result can be considered as the analogue of the corresponding limit for classical Boltzmann equation \cite{Jiang-Xu-Zhao-2018} and BFD equation \cite{jiang2021incompressible}. All of these results belong to the so-called ``bottom-up" type fluid limits, as classified in \cite{Jiang-Xu-Zhao-2018}. More specifically, these limits do not rely on the existence of the limiting equations. In fact, these limits provide the solutions of the limiting equations from the solutions of the kinetic equations.

We will start from the following scaled BBE equation (in diffusive scaling):
\ben \label{quantum-Boltzmann-CP}
\partial _t F +  \frac{1}{\epsilon} v \cdot \nabla_{x} F= \frac{1}{\epsilon^{2}}Q(F,F), ~~t > 0, x \in \mathbb{R}^{3}, v \in \R^3,
\een
where the operator $Q$ is defined through \eqref{U-U-operator-scaling-hard-sphere} with $B$ given in \eqref{hard-sphere-kernel}, and $\delta=1$ since we study BBE equation in this paper. In the above expression, $\epsilon>0$ is the Knudsen number. The so-called {\em hydrodynamic limit} is the process that the Knudsen number $\epsilon\rightarrow 0$.

Our goal in this paper is to rigorously justify the limit as $\epsilon \to 0$, from solutions of the BBE \eqref{quantum-Boltzmann-CP} to solutions of the incompressible Navier-Stokes-Fourier (NSF) equations \eqref{NSF-problem}.
Precise definitions of solutions to these equations will be given soon.
%The derivation from BBE to NSF will be illustrated in Section 5.

\subsection{Temperature in Boltzmann-Bose-Einstein equation} Temperature plays an important role in the study of quantum Boltzmann equation. For example, for the particles obeying Bose-Einstein statistics, the so-called Bose-Einstein-Condensation (BEC) only happens in low temperature. For example, see \cite{EV-Inventions2016}. We now introduce some basic knowledge about temperature in the quantum context.
Let us consider a homogeneous density $f=f(v)$ with zero mean $\int f(v) \mathrm{d}v=0$. For $k \geq 0$, we recall the moment function
\beno
M_{k}(f) \colonequals   \int |v|^{k} f(v)  \mathrm{d}v.
\eeno
Let $M_{0} = M_{0}(f), M_{2}=M_{2}(f)$ for simplicity. Let $m$ be the mass of a particle, then
 $m M_{0}$ and $\f{1}{2}m M_{2}$ are the total mass and kinetic energy per unit space volume.
Referring \cite{lu2004isotropic}, for a given density function $f$,
the kinetic temperature $\bar{T}$ and the critical temperature $\bar{T}_{c}$ of the particle system are defined by
\ben \label{kinetic-temperature-and-critical}
\bar{T} = \frac{1}{3 k_{B}} \frac{m M_{2}}{M_{0}}, \quad \bar{T}_{c} = \frac{m\zeta(5/2)}{2 \pi k_{B} \zeta(3/2)}
 \big( \frac{M_{0} }{\zeta(3/2)} \big)^{\f23},
\een
where $k_{B}$ is the Boltzmann constant and $\zeta(s) = \sum_{n=1}^{\infty} \frac{1}{n^{s}}$ is the Riemann zeta function.
The ratio $\bar{T}/\bar{T}_{c}$ quantifies high and low temperature. More precisely, high temperature  $\bar{T}/\bar{T}_{c} > 1$; critical case  $\bar{T}/\bar{T}_{c} = 1$; low temperature  $\bar{T}/\bar{T}_{c} < 1$.

We now recall some known results about equilibrium distribution. The famous Bose-Einstein distribution has density
function
\ben
\label{equilibrium-mathcal-M}
 \mathcal{M}_{\lambda,T}(v) \colonequals  \frac{1}{\exp(\f{|v|^{2}}{2T} + \lambda) - 1};\quad \lambda \geq 0, T>0
\een
The ratio $\bar{T}/\bar{T}_{c}$ of $\mathcal{M}_{\lambda,T}$ depends only on $\lambda$. The critical value $\lambda=0$ corresponds to the critical temperature   $\bar{T}/\bar{T}_{c} = 1$. The critical value $\lambda>0$ corresponds to high temperature  case $\bar{T}/\bar{T}_{c} > 1$. In low temperature  $\bar{T}/\bar{T}_{c} < 1$,  the equilibrium of BBE equation is the Bose-Einstein distribution \eqref{equilibrium-mathcal-M} with
$\lambda=0$ plus some Dirac delta function. That is, the equilibrium contains a Dirac measure. One can refer to \cite{lu2005boltzmann} for the classification of equilibria.

In this article, we work with high temperature and consider the equilibrium $\mathcal{M}_{\lambda,T}$ with $\lambda>0$. Note that very high temperature assumption is imposed in \cite{li2019global} to prove global well-posedness of homogeneous BBE equation with (slightly general than) hard sphere collisions. We consider the inhomogeneous case and give a careful analysis on the dependence on
$\lambda,T>0$. In particular, we pay much attention to the effects of $\lambda, T \to 0$ or $\lambda, T \to \infty$. Our analysis illustrates that for BBE, the incompressible fluid limits (at least for smooth solutions) only happen at high temperature.

%Recall that the solution of \eqref{quantum-Boltzmann-UU} conserves mass, momentum and energy. That is, for any $t \geq 0$,
%\ben \label{conversation-mass-momentum-energy}
%\int (1, v, |v|^{2}) F(t,x,v) \mathrm{d}x \mathrm{d}v = \int (1, v, |v|^{2}) F_{0}(x,v) \mathrm{d}x \mathrm{d}v.
%\een
%Once $F_{0}$ is appropriately given, the constants $\rho,T>0, v_{0} \in \mathbb{R}^{3}$ in \eqref{equilibrium-mathcal-M} are uniquely determined through % (ref. \cite{lu2005boltzmann} for B-E particles and ref. Lu for F-D particles)
%\ben \label{F0-determines-equilibrium}
%\int (1, v, |v|^{2}) F_{0}(x,v) \mathrm{d}x \mathrm{d}v = \int (1, v, |v|^{2}) \mathcal{M}_{\rho, v_{0}, T}(v) \mathrm{d}x \mathrm{d}v.
%\een
%Without any loss of generality, we assume that $F_{0}$ has zero mean and thus gives $v_{0}=0$ from now on. As a result, we only keep $\rho, T$ as two variable parameter. That is, we only consider those initial conditions that give $v_{0}=0$  according to \eqref{F0-determines-equilibrium}. Taking $\mu_{\lambda,T} (v) \colonequals   \rho (2\pi T)^{-\frac{3}{2}} e^{-\f{1}{2T}|v|^{2}},$

\subsection{Perturbation around equilibrium and main results} Similarly to the classical Boltzmann equations and BFD, the incompressible fluid regimes are near global equilibrium and long time scale. For the detailed explanation, see \cite{BGL-1991JSP} for Boltzmann and \cite{Zakrevskiy} for BFD.

For perturbation around the equilibrium, we define
\ben \label{muliplier-mathcal-N}
\mathcal{N}_{\lambda,T}(v) \colonequals  \sqrt{\mathcal{M}_{\lambda,T}(v)(1+\mathcal{M}_{\lambda,T}(v))} = \frac{1}{\exp(\f{|v|^{2}}{4T} + \f{\lambda}{2}) - 1}.
\een
We remark that the function $\mathcal{N}_{\lambda,T}$ serves as the multiplier in the expansion $F = \mathcal{M}_{\lambda,T} + \mathcal{N}_{\lambda,T} f$.
For simplicity, let $\mathcal{M}\colonequals   \mathcal{M}_{\lambda,T}$ and $ \mathcal{N}\colonequals   \mathcal{N}_{\lambda,T}$. With the expansion $F = \mathcal{M} + \epsilon \mathcal{N} f$,
the perturbed quantum Boltzmann equation corresponding to \eqref{quantum-Boltzmann-CP} reads
\ben \label{linearized-quantum-Boltzmann-eq}
\partial _t f + \f{1}{\epsilon} v \cdot \nabla_{x} f + \f{1}{\epsilon^{2}} \mathcal{L}^{\lambda,T}f = \f{1}{\epsilon} \Gamma_{2}^{\lambda,T}(f,f) + \Gamma_{3}^{\lambda,T}(f,f,f), \quad
f|_{t=0} = f_{0}^{\epsilon}.
\een
Here the linearized quantum Boltzmann operator $\mathcal{L}^{\lambda,T}$ is define by
\ben \label{linearized-quantum-Boltzmann-operator-UU}
(\mathcal{L}^{\lambda,T}f)(v)  \colonequals    \int  B \mathcal{N}_{*} \mathcal{N}^{\prime} \mathcal{N}^{\prime}_{*} \mathrm{S}(\mathcal{N}^{-1}f)
\mathrm{d}\sigma \mathrm{d}v_{*},
\een
where $\mathrm{S}(\cdot)$ is defined by
\ben  \label{symmetry-operator}
\mathrm{S}(g) \colonequals      g + g_{*} - g^{\prime} - g^{\prime}_{*}.
\een
The bilinear term $\Gamma_{2}^{\lambda,T}(\cdot,\cdot)$ and the trilinear term $\Gamma_{3}^{\lambda,T}(\cdot,\cdot,\cdot)$ are defined by
\ben \label{definition-Gamma-2-epsilon}
\Gamma_{2}^{\lambda,T}(g,h) &\colonequals &   \mathcal{N}^{-1}\int
B \Pi_{2}(g,h)
\mathrm{d}\sigma \mathrm{d}v_{*}.
\\ \label{definition-Gamma-3-epsilon}
\Gamma_{3}^{\lambda,T}(g,h,\varrho)  &\colonequals &     \mathcal{N}^{-1}\int
B  \mathrm{D} \big( (\mathcal{N}g)_{*}^{\prime} (\mathcal{N}h)^{\prime} ((\mathcal{N}\varrho)_{*} + \mathcal{N}\varrho) \big) \mathrm{d}\sigma \mathrm{d}v_{*}. \quad
\quad
\een
The notation $\Pi_{2}$ in \eqref{definition-Gamma-2-epsilon} is defined by
 \ben   \label{definition-A-2}
\Pi_{2}(g,h) &\colonequals   & \mathrm{D} \big((\mathcal{N}g)^{\prime}_{*}(\mathcal{N}h)^{\prime}\big)
\\ \label{line-1} &&+
\mathrm{D}\big((\mathcal{N}g)^{\prime}_{*}(\mathcal{N}h)^{\prime}(\mathcal{M} +  \mathcal{M}_{*})\big)
\\ \label{line-2} &&+ \mathrm{D}\big( (\mathcal{N}g)_{*}(\mathcal{N}h)^{\prime}( \mathcal{M}^{\prime}_{*} -  \mathcal{M}) \big)
\\ \label{line-3} &&+ \big( (\mathcal{N}g)^{\prime}(\mathcal{N}h)\mathrm{D}(\mathcal{M}^{\prime}_{*}) + (\mathcal{N}g)^{\prime}_{*}(\mathcal{N}h)_{*}\mathrm{D}(\mathcal{M}^{\prime})\big).
\een
Remark that the three operators $\mathcal{L}^{\lambda,T}$, $\Gamma_{2}^{\lambda,T}(\cdot,\cdot)$ and $\Gamma_{3}^{\lambda,T}(\cdot,\cdot,\cdot)$ depends on $\lambda, T$ through $\mathcal{M} =  \mathcal{M}_{\lambda,T}$, and $\mathcal{N} =  \mathcal{N}_{\lambda,T}$.

The main result of this paper is to prove global well-posedness of \eqref{linearized-quantum-Boltzmann-eq} {\em uniformly} in $\epsilon$ in the Sobolev space $H^{N}_{x}L^{2}$, which is defined as
\ben \label{definition-energy}
 \|f\|_{H^{N}_{x}L^{2}}^{2} \colonequals    \sum_{|\alpha|\leq N} \|\partial^{\alpha}_{x} f \|_{L^{2}_{x}L^{2}}^{2}.
\een
Note that the functional only involves $x$-derivatives. Correspondingly, the dissipation functional reads
\ben \label{dissipation-ori}
\mathcal{D}_{N,T}(f)
\colonequals \mathcal{D}_{N,T,1}(f) + \mathcal{D}_{N,T,2}(f).
\een
where 
\ben \label{def-D1-D2}
\mathcal{D}_{N,T,1}(f)
\colonequals  |\nabla_{x} \mathbb{P}_{\lambda,T} f|_{H^{N-1}_{x}L^{2}}^{2}, \quad \mathcal{D}_{N,T,2}(f)
\colonequals \|f - \mathbb{P}_{\lambda,T} f\|_{H^{N}_{x}L^{2}}^{2}
+ T^{-\f{1}{2}}\||v|^{\f{1}{2}}(f - \mathbb{P}_{\lambda,T} f)\|_{H^{N}_{x}L^{2}}^{2}.
\een

Here $\mathbb{P}_{\lambda,T}$ is the projection on the kernel space of $\mathcal{L}^{\lambda,T}$, which will be defined in  \eqref{definition-projection-operator}.

\begin{thm}\label{global-well-posedness-original}

Let $0<\epsilon \leq 1$, $\lambda, T >0$, and $N \geq 2$. Let
\ben \label{C-star-final}
C_{*}(\lambda,T) \colonequals e^{-2\lambda} (1-e^{-\lambda})^{\f{33}{2}} \min\{T^{3/2},T^{-3/2}\}.
\een
Then there exists a universal constant $\delta_{*}>0$, independent of  $\epsilon, \lambda$ and $T$, such that if the initial datum $f_0$ satisfying
\ben \label{condition-on-initial-data-original} \mathcal{M}_{\lambda, T} + \mathcal{N}_{\lambda, T}f_{0} \geq 0, \quad \|f_{0}\|_{H^{2}_{x}L^{2}}^{2} \leq \delta_{*} C_{*}(\lambda,T), \quad \|f_{0}\|_{H^{N}_{x}L^{2}}<\infty\,, \een
then the Cauchy problem \eqref{linearized-quantum-Boltzmann-eq} with initial datum $f_0$ has a unique global solution $f = f^{\lambda,T}_{\epsilon} \in L^{\infty}([0,\infty); H^{N}_{x}L^{2})$ satisfying $\mathcal{M}_{\lambda, T} + \mathcal{N}_{\lambda, T} f(t) \geq 0$ and
\ben \label{uniform-estimate-global-initial}
\sup_{t}\|f(t)\|_{H^{N}_{x}L^{2}}^{2}+ \f{1}{K(\lambda,T)}  \int_{0}^{\infty} \mathcal{D}_{N,T}(f) \mathrm{d}\tau + \f{\mathrm{C}_{1}(\lambda,T)}{K(\lambda,T)}  \f{1}{\epsilon^{2}} \int_{0}^{\infty} \mathcal{D}_{N,T,2}(f) \mathrm{d}\tau \leq  O_{N}(f_{0})\|f_0\|_{H^{N}_{x}L^{2}}^{2}.
\een
where \ben \label{definition-of-Qn}
O_{2}(f_{0}) \equiv 12, \quad O_{N}(f_{0}) \colonequals  24 \exp\left(Q_{3}(\lambda,T,N,f_{0}) O_{N-1}(f_{0})T^{-\f{3}{2}}\|f_{0}\|_{H^{N-1}_{x}L^{2}}^{2} \right)   \text{ for } N \geq 3,
\\ \label{definition of-Q-3-f0}
Q_{3}(\lambda,T,N,f_{0}) \colonequals  2(Q_{1}(\lambda,T,N)  +
Q_{2}(\lambda,T,N) O_{N-1}(f_{0})T^{-\f{3}{2}}\|f_{0}\|_{H^{N-1}_{x}L^{2}}^{2}).
\een
where $Q_{1}(\lambda,T,N)$ and $Q_{2}(\lambda,T,N)$ are some constants defined in \eqref{defintion-Q-1-constant} and \eqref{defintion-Q-1-constant}.
 Here the constants $K(\lambda,T), \mathrm{C}_{1}(\lambda,T)$ are defined defined in \eqref{defintion-of-K-lambda-T}, \eqref{defintion-of-C1-lambda-T}.
\end{thm}

We make some remarks on the above theorem.

\begin{rmk} We emphasize that the constant $\delta_{*}$ in Theorem \ref{global-well-posedness-original} is universal and does not depend on anything.
 We give an explicit smallness assumption in terms of the two parameters $\lambda, T >0$. Note that the well-posedness region
 vanishes (i.e. the constant $C_{*}(\lambda,T)$ defined in \eqref{C-star-final} will tends to $0$) as $\lambda \to 0$ or $\lambda \to \infty$ or $T \to 0$ or $T \to \infty$.
\end{rmk}

\begin{rmk}
Recall that $H^{2}_{x}L^{2}$ might be the largest the Sobolev space (with integer index) in which global well-posedness can be established for the classical Boltzmann equation in the whole space. We manage to construct the global well-posedness theory in the space for the bosonic Nordheim Boltzmann equation. %Note that the work \cite{ouyang2021quantum} needs $H^{8}_{x,v}$ with order $8$ mixed $x,v$ derivatives. Here $ \|f\|_{H^{n}_{x,v}}^{2}  %\colonequals      \sum_{|\alpha| + |\beta| \leq n} \|\partial^{\alpha}_{x} \partial^{\beta}_{v}f\|_{L^{2}_{x}L^{2}}^{2}.$
\end{rmk}

\begin{rmk}
Notice that in the condition \eqref{condition-on-initial-data-original}, only smallness assumption on $\|f^\epsilon_{0}\|_{H^{2}_{x}L^{2}}$ is imposed.
In other words, $\|f^\epsilon_{0}\|_{H^{N}_{x}L^{2}}(N \geq 3)$ could be arbitrarily large (and bounded).
\end{rmk}

%In terms of physical relevance of Theorem \ref{Gamma-3-energy-estimate}, we give the following two remarks.
%\begin{rmk} (not true now, consider it later) Note that $\|f^{\lambda,T}(t)\|_{L^{\infty}_{x,v}}^{2} \lesssim \mathcal{E}_{N}(f^{\lambda,T}(t)) \lesssim \rho$ and thus
%$\|\frac{\rho \mu}{1 - \rho \mu} + \frac{\rho^{\f{1}{2}}\mu^{\f{1}{2}}}{1 - \rho \mu}f^{\lambda,T}(t)\|_{L^{\infty}_{x,v}}\lesssim \rho$ for any $t \geq 0$.
%Theorem \ref{global-well-posedness} shows that BEC will not happen if the initial datum is close enough to
%an equilibrium with high temperature, which is reasonable and consistent with physical observation.
%\end{rmk}
%\begin{rmk}
%By the scaling \eqref{scaling-tranform}, Theorem \ref{global-well-posedness} ensures the global well-posedness of \eqref{quantum-Boltzmann-UU} with $\Phi(x) = \hbar^{4p-4}|x|^{-p}$ and initial datum $F_{0}$ close enough to the equilibrium $\mathcal{M}_{\rho, \hbar}$ where $\mathcal{M}_{\rho, \hbar}(v) = \hbar^{-3}\mathcal{M}_{\lambda,T}(\hbar^{3} v)$. Simply looking at the equilibrium $\mathcal{M}_{\rho, \hbar}$, we have $\|\mathcal{M}_{\rho, \hbar}\|_{L^{\infty}_{x,v}} \sim \hbar^{-3} \rho, \|\mathcal{M}_{\rho, \hbar}\|_{L^{2}_{x}L^{2}} \sim \hbar^{-\f{15}{2}} \rho.$ Since $\hbar$ is a relatively small constant, the magnitude of the solution to the problem \eqref{quantum-Boltzmann-UU} can be relatively large.
%\end{rmk}

Based on Theorem \ref{global-well-posedness-original}, we can prove the hydrodynamical limit of \eqref{linearized-quantum-Boltzmann-eq} to incompressible Navier-Stokes-Fourier equations.

\begin{thm} \label{thm-hydrodynamical-limit-ori}
Let $\delta_{*}$ and $C_{*}(\lambda,T)$ be the constants in Theorem \ref{global-well-posedness-original}.
Let $0<\epsilon \leq 1$, $\lambda, T >0$, and $N \geq 2$.
Let $f^{\epsilon}_{0}$ be a family of initial datum satisfying $\mathcal{M}_{\lambda, T} +  \mathcal{N}_{\lambda, T} f^{\epsilon}_{0} \geq 0$ and
\ben \label{small-condition-ori}
\sup_{0<\epsilon<1} \|f^{\epsilon}_{0}\|_{H^{2}_{x}L^{2}}^{2} \leq \delta_{*} C_{*}(\lambda,T), \quad
M_{0} \colonequals  \sup_{0<\epsilon<1} \|f^{\epsilon}_{0}\|_{H^{N}_{x}L^{2}} < \infty,
\\
\mathbb{P}_{\lambda, T}f^{\epsilon}_{0} \to f_{0} = (\rho_{0} + u_{0} \cdot \f{v}{T^{1/2}} + \theta_{0} (\f{|v|^{2}}{2 T} - K_{\lambda})) \mathcal{N}_{\lambda, T} \text{ as }
\epsilon \to 0,
\text{ strongly in } H^{N}_{x}L^{2}
\een
for some $(\rho_{0}, u_{0}, \theta_{0}) \in H^{N}_{x}$ with $\rho_{0} + \theta_{0} = 0$. Here $K_{\lambda} = K_{A} -1$ is a constant depending only on $\lambda$. Let $f^{\epsilon}$ be the solution to the Cauchy problem \eqref{linearized-quantum-Boltzmann-eq}
 with initial datum $f^{\epsilon}_{0}$. Then there is a subsequence of $\{f^{\epsilon}\}$ still denoting it by $\{f^{\epsilon}\}$ such that,
 \ben \label{weakly-star-convergence-to-some-special-form-ori}
 f^{\epsilon} \to (\rho + u \cdot \f{v}{T^{1/2}} + \theta (\f{|v|^{2}}{2T} - K_{\lambda})) \mathcal{N}_{\lambda, T} \text{ as }
\epsilon \to 0, \text{ weakly-* in } L^{\infty}(\mathbb{R}_{+}; H^{N}_{x}L^{2}).
 \een
 for some $(\rho, u, \theta) \in L^{\infty}(\mathbb{R}_{+}; H^{N}_{x}) \cap C(\mathbb{R}_{+}; H^{N-1}_{x})$ satisfying \ben\label{NSF-problem} \left\{ \begin{aligned}
& \rho + \theta = 0, \quad \nabla_{x} \cdot u = 0, \\
&\partial_{t} u + T^{\f{1}{2}} \mathcal{P} \nabla_{x} \cdot ( u \otimes u) = \mu_{\lambda, T} \Delta_{x} u,
\\
&\partial_{t} \theta + T^{\f{1}{2}} u \cdot \nabla_{x}\theta = \kappa_{\lambda, T} \Delta_{x} \theta,
\\
& \rho|_{t=0} = \rho_{0}, \quad
u|_{t=0} = \mathcal{P}u_{0}, \quad \theta|_{t=0} = \f{K_{\lambda}\theta_{0} - \rho_{0}}{K_{\lambda}+1}.
\end{aligned} \right.
\een
 In addition,
 \ben \label{momoent-1-ori}
\f{3T^{-\f{3}{2}}}{m_{2}}\mathcal{P}\langle f^{\epsilon}, \f{v}{T^{1/2}} \mathcal{N}_{\lambda, T} \rangle \to u \text{ strongly in } C(\mathbb{R}_{+}; H^{N-1}_{x}) \text{ and weakly-* in } L^{\infty}(\mathbb{R}_{+}; H^{N}_{x}),
\\ \label{momoent-2-ori}
\f{T^{-\f{3}{2}}}{C_{A}} \langle f^{\epsilon}, (\f{|v|^{2}}{2T} - K_{A}) \mathcal{N}_{\lambda, T} \rangle \to \theta \text{ strongly in } C(\mathbb{R}_{+}; H^{N-1}_{x}) \text{ and weakly-* in } L^{\infty}(\mathbb{R}_{+}; H^{N}_{x}).
 \een
\end{thm}

\begin{rmk}
Only smallness assumption on $|(\rho_{0}, u_{0}, \theta_{0})|_{H^{2}_{x}}$ is imposed.
That is, $|(\rho_{0}, u_{0}, \theta_{0})|_{H^{N}_{x}}(N \geq 3)$ can be arbitrarily large. This point is new in the hydrodynamic limit.
\end{rmk}

\subsection{Novelties of the results} In terms of condition and conclusion, Theorem \ref{global-well-posedness} is closest to the main result (Theorem 1.4) of \cite{li2019global}. More precisely, both of these two results validate global existence of anisotropic solutions. The main differences of these two results are also obvious. To reiterate, \cite{li2019global} considers spatially homogeneous case under very high temperature condition. We work in the whole space $x \in \mathbb{R}^{3}$ under any temperature higher than the critical one.

In terms of mathematical methods, this article is closer to \cite{bae2021relativistic} and \cite{jiang2021incompressible} since all of these works fall into the close-to-equilibrium framework well established for the classical Boltzmann equation. There are many works that contribute to this mathematically satisfactory theory for global well-posedness of the classical Boltzmann equation. For readers reference, we mention
\cite{ukai1974existence,guo2003classical} for angular cutoff kernels and \cite{gressman2011global,alexandre2012boltzmann} for non-cutoff kernels.

%Note that there are also some results on general spatial domains with physical boundary conditions and solutions with polynomial decays. Since these topics are beyond the scope of this work, we do not review much on them.
%Interested readers can refer to A, B and the references therein.

Each of  \cite{bae2021relativistic} and \cite{jiang2021incompressible} has their own features and focuses. The work \cite{bae2021relativistic} is the first to investigate both relativistic and quantum effect.
The article \cite{jiang2021incompressible} studies the hydrodynamic limit from BFD (but not BBE) to incompressible Navier-Stokes-Fourier equation. These works contribute to the literature of quantum Boltzmann equation from different aspects.

Our main results Theorem \ref{global-well-posedness-original} and Theorem \ref{thm-hydrodynamical-limit-ori}  have some unique features that may better our understanding of quantum Boltzmann equation. Besides the low regularity requirement of the solution space $H^2_xL^2$,  in particular, this article may be the first to investigate
\begin{itemize}
\item well-posedness theory of BBE for any temperature higher than the critical one in terms of the two parameters $\lambda, T>0$.
\item hydrodynamic limit from BBE to incompressible Navier-Stokes-Fourier equation.
\end{itemize}
Specifically, Theorem \ref{global-well-posedness-original} precisely state the dependence of the existence regime on the parameters $\lambda$ and the temperature $T$. Furthermore, technically, we obtain the uniform in Knudsen number $\epsilon$ estimate. This can only be achieved under the incompressible Navier-Stokes scaling. In this sense, our limit is a ``bottom-up" type, i.e. we start from the solutions of microscopic kinetic equation, and take the limit $\epsilon \rightarrow 0$ to automatically prove the global existence of the limiting equations. In the whole process, we do not need any information of the limiting equations. This is quite different with the expansion method. The compressible Euler and acoustic limits will be treated in a separate paper which is under preparation.

%Different from the three works (\cite{bae2021relativistic}, \cite{ouyang2021quantum} and \cite{jiang2021incompressible}),
%we keep the parameter $\rho$ along our derivation throughout the article. This intentional choice enables us to relate the high temperature condition to the smallness of $\rho$ quantitatively in \eqref{high-temperature-condition}.
%As a result, BEC is ruled out globally.

\subsection{Notations} \label{notation} In this subsection, we give a list of notations.

\noindent $\bullet$ Given a set $A$,  $\mathrm{1}_A$ is the characteristic function of $A$.

\noindent $\bullet$ The notation $a\lesssim b$  means that  there is a universal constant $C$
such that $a\leq Cb$.

% The constant $C$ could depend on the kernel parameters $\gamma, s$ and the energy space index $N, l_{0,0}$.

\noindent $\bullet$ If both $a\lesssim b$ and $b \lesssim a$, we write $a\sim b$.

\noindent $\bullet$ We denote $C(\lambda_1,\lambda_2,\cdots, \lambda_n)$ or $C_{\lambda_1,\lambda_2,\cdots, \lambda_n}$  by a constant depending on $\lambda_1,\lambda_2,\cdots, \lambda_n$.

\noindent $\bullet$ The bracket $\langle \cdot\rangle$ is defined by $\langle v \rangle \colonequals    (1+|v|^2)^{\f{1}{2}}$. The weight function  $W_l(v)\colonequals     \langle v\rangle^l $.

\noindent $\bullet$ For $f,g \in L^{2}({\R^3})$,  $\langle f,g\rangle\colonequals     \int_{\R^3}f(v)g(v) \mathrm{d}v$ and $|f|_{L^{2}}^{2}\colonequals    \langle f,f\rangle$.

\noindent $\bullet$ For $f,g \in L^{2}({\R^3})$,   $\langle f,g\rangle_{x}\colonequals     \int_{\R^3}f(x)g(x) \mathrm{d}x$ and $|f|_{L^{2}_{x}}^{2}\colonequals    \langle f,f\rangle_{x}$.

\noindent $\bullet$ For $f,g \in L^{2}({\R^3 \times \R^3})$,   $(f,g)\colonequals     \int_{\R^3 \times \R^3} f(x,v)g(x,v) \mathrm{d}x\mathrm{d}v$ and $\|f\|_{L^{2}_{x}L^{2}}^{2}\colonequals    (f, f)$.

\noindent $\bullet$ For a multi-index
$\alpha =(\alpha_1,\alpha_2,\alpha_3) \in \mathbb{N}^{3}$, define
$|\alpha|\colonequals    \alpha_1+\alpha_2+\alpha_3$.

\noindent $\bullet$ For  $\alpha \in \mathbb{N}^{3}$
denote $\partial^{\alpha}\colonequals    \partial^{\alpha}_{x}$.

We now introduce some norm.

\noindent $\bullet$ For $l \geq 0$ and a function $f(v)$ on $\mathbb{R}^{3}$, define
\ben \label{l2-l-not-like-before}
|f|_{L^{2}_{l}}^{2} \colonequals  |f|_{L^{2}}^{2} +  \mathrm{1}_{l > 0}||\cdot|^{l}f|_{L^{2}}^{2}, \quad |f|_{L^{2}}=|f|_{L^{2}_{0}}.
\een

\noindent $\bullet$ For $n \in \mathbb{N}$ and a function $f(x)$ on $\mathbb{R}^{3}$, define
\beno
% \label{Sobolev-norm-x}
|f|_{H^{n}_{x}}^{2} \colonequals     \sum_{|\alpha| \leq n} |\partial^{\alpha}f|_{L^{2}_{x}}^{2},  \quad |f|_{L^{2}_{x}}:=|f|_{H^{0}_{x}}, \quad |f|_{L^{\infty}_{x}} \colonequals   \esssup_{x \in \mathbb{R}^{3}}  |f(x)|.
\eeno

\noindent $\bullet$
For $m \in \mathbb{N}, l \geq 0$ and a function $f(x,v)$ on $\mathbb{R}^{3}\times \mathbb{R}^{3}$, define
\ben \label{not-mix-x-v-norm-energy}
\|f\|_{H^{m}_{x}L^{2}_{l}}^{2} \colonequals      \sum_{|\alpha| \leq m} ||\partial^{\alpha} f|_{L^{2}_{l}}|_{L^{2}_{x}}^{2}, \quad \|f\|_{L^{2}_{x}L^{2}_{l}} \colonequals    \|f\|_{H^{0}_{x}L^{2}_{l}}, \quad \|f\|_{H^{m}_{x}L^{2}} \colonequals    \|f\|_{H^{m}_{x}L^{2}_{0}},
\een

\subsection{ Organization of this paper} \label{plan}
Section \ref{scaling} introduces a simple scaling. (We may also put this section in the introduction) Section \ref{operator-analysis} contains
 estimates of linear operators and non-linear operators, including coercivity and upper bound estimate.
% In Section \ref{local}, we derive local well-posedness.
 In Section \ref{global}, we first prove a priori estimate and then establish global well-posedness. In Section \ref{hydrodynamic-limit}, we give hydrodynamic limit.
 Section \ref{appendix} is an appendix in which we put some elementary proof for the sake of completeness.

\section{Scaled Boltzmann-Bose-Einstein equation} \label{scaling}

%Check single $|\cdot|$ or double $\|\cdot\|$.

In this section, we introduce a simple scaling to free us from the parameter $T>0$.
 Fix $a>0$. Let us define an operator $A_{a}$ by
 $A_{a}f(v) \colonequals   f(av)$ for all $v \in \mathbb{R}^{3}$.
Making $A_{T^{1/2}}$ to $\mathcal{M}_{\lambda,T}$ and $\mathcal{N}_{\lambda,T}$, we get
\ben \label{definition-of-M-lambda}
M_{\lambda}(v)\colonequals  (A_{T^{1/2}} \mathcal{M}_{\lambda,T})(v) =
 \f{1}{ \exp(\f{|v|^{2}}{2}+\lambda) -1}, \quad N_{\lambda}(v)\colonequals  (A_{T^{1/2}} \mathcal{N}_{\lambda,T})(v) =
 \f{\exp(\f{|v|^{2}}{4}+\f{\lambda}{2})}{ \exp(\f{|v|^{2}}{2}+\lambda) -1},
\een
Let us make the action $A_{T^{1/2}}$ to \eqref{linearized-quantum-Boltzmann-eq}. Then the equation \eqref{linearized-quantum-Boltzmann-eq} becomes
\ben \label{T-half-scaled-more-consistent}
 \partial _t \tilde{f} + \f{1}{\epsilon} T^{1/2} v \cdot \nabla_{x} \tilde{f} + \f{1}{\epsilon^{2}} \tilde{\mathcal{L}}^{\lambda,T}\tilde{f} = \f{1}{\epsilon} \tilde{\Gamma}_{2}^{\lambda,T}(\tilde{f},\tilde{f}) + \tilde{\Gamma}_{3}^{\lambda,T}(\tilde{f},\tilde{f},\tilde{f}), \quad \tilde{f}|_{t=0} = \tilde{f}_{0}.
\een
where $\tilde{f} = A_{T^{1/2}} f, \tilde{f}_{0}=A_{T^{1/2}} f_{0}^{\epsilon}$. Here the linear operator $\tilde{\mathcal{L}}^{\lambda,T}$ is defined by 
\ben \label{scaled-linearized-operator-more-consistent}
(\tilde{\mathcal{L}}^{\lambda,T} f)(v)  \colonequals   \int  B_{T}(v-v_{*}, \sigma) (N_{\lambda})_{*} (N_{\lambda})^{\prime} (N_{\lambda})^{\prime}_{*} \mathrm{S}(N_{\lambda}^{-1}f)
\mathrm{d}\sigma \mathrm{d}v_{*}.
\een
with
\beno B_{T}(v-v_{*}, \sigma) = T^{\f{3}{2}} B(T^{1/2}(v-v_{*}), \sigma) = T^{2}|v-v_{*}|.
\eeno
The bilinear operator
$\tilde{\Gamma}_{2}^{\lambda,T}$ is defined by
\ben \label{definition-Gamma-2-epsilon-rho-T}
\tilde{\Gamma}_{2}^{\lambda,T}(g,h) \colonequals    N_{\lambda}^{-1}\int
B_{T} \tilde{\Pi}_{2}(g,h)
\mathrm{d}\sigma \mathrm{d}v_{*}.
\een
with $\tilde{\Pi}_{2}$ given by
\ben   \label{definition-A-2-rho-T}
\tilde{\Pi}_{2}(g,h) &\colonequals   & \mathrm{D} \big((N_{\lambda}g)^{\prime}_{*}(N_{\lambda}h)^{\prime}\big)
\\ \label{line-1-rho-T} &&+
\mathrm{D}\big((N_{\lambda}g)^{\prime}_{*}(N_{\lambda}h)^{\prime}(M_{\lambda} +  (M_{\lambda})_{*})\big)
\\ \label{line-2-rho-T} &&+
\mathrm{D}\big( (N_{\lambda}g)_{*}(N_{\lambda}h)^{\prime}( (M_{\lambda})^{\prime}_{*} -  M_{\lambda}) \big)
\\ \label{line-3-rho-T} &&+
\big( (N_{\lambda}g)^{\prime}(N_{\lambda}h)\mathrm{D}((M_{\lambda})^{\prime}_{*}) + (N_{\lambda}g)^{\prime}_{*}(N_{\lambda}h)_{*}\mathrm{D}(M_{\lambda}^{\prime})\big).
\een
The trilinear operator
$\tilde{\Gamma}_{3}^{\lambda,T}$ is defined by
\ben \label{definition-Gamma-3-epsilon-rho-T}
\tilde{\Gamma}_{3}^{\lambda,T}(g,h,\varrho) (v)  \colonequals
N_{\lambda}^{-1}   \int B_{T}(v-v_{*}, \sigma)  \mathrm{D} \big( (N_{\lambda}g)_{*}^{\prime} (N_{\lambda}h)^{\prime} ((N_{\lambda}\varrho)_{*} + N_{\lambda}\varrho) \big) \mathrm{d}\sigma \mathrm{d}v_{*}.
\een
Indeed, by the change of variable $v_{*} \to T^{-1/2}v_{*}$, $\mathrm{d}v_{*}  = T^{3/2} \mathrm{d}(T^{-1/2}v_{*})$, it is easy to find
\beno % \label{relation-to-scaled-operator-more-consistent}
A_{T^{1/2}} \mathcal{L}^{\lambda,T}f  = \tilde{\mathcal{L}}^{\lambda,T} A_{T^{1/2}} f,
\\
A_{T^{1/2}} \Gamma_{2}^{\lambda,T}(g,h)  = \tilde{\Gamma}_{2}^{\lambda,T}  (A_{T^{1/2}}g, A_{T^{1/2}}h),
\\
A_{T^{1/2}} \Gamma_{3}^{\lambda,T}(g,h,\varrho)  = \tilde{\Gamma}_{3}^{\lambda,T}  (A_{T^{1/2}}g, A_{T^{1/2}}h, A_{T^{1/2}}\varrho),
\eeno

\begin{prop} \label{solution-relation}
$f$ is a solution to  \eqref{linearized-quantum-Boltzmann-eq} with initial datum $f_{0}$ if and only if $A_{T^{1/2}}f$ is a solution to \eqref{T-half-scaled-more-consistent} with initial datum $A_{T^{1/2}}f_{0}$.
\end{prop}

\begin{rmk} Another scaling choice is to consider the equation of $g(x,v) := f(T^{1/2}x, T^{1/2}v)$. Note that
\beno
(v \cdot \nabla_{x} g)(x,v) = (v \cdot \nabla_{x} f) (T^{1/2}x, T^{1/2}v).
\eeno
So the equation for $g$ is
\ben \label{another-scaling-more-consistent}
 \partial _t g + \f{1}{\epsilon} v \cdot \nabla_{x} g + \f{1}{\epsilon^{2}} \tilde{\mathcal{L}}^{\lambda,T} g = \f{1}{\epsilon} \tilde{\Gamma}_{2}^{\lambda,T}(g,g) + \tilde{\Gamma}_{3}^{\lambda,T}(g,g,g), \quad g|_{t=0} = g_{0}:=f_{0}(T^{1/2}x, T^{1/2}v).
\een
This choice gives a simple streaming term because the factor $T^{1/2}$ in \eqref{T-half-scaled-more-consistent} disappears.
However, this choice also has some drawbacks. First,
the energy functional contains derivatives of the $x$ variable and so some additional factor $T^{1/2}$ comes out.
Second, such choice may change the spatial domain if 
one wants to deal with bounded domain such as the
torus $[0, l_{x}]^{3}$.
\end{rmk}

Recalling \eqref{linearized-quantum-Boltzmann-operator-UU} and \eqref{scaled-linearized-operator-more-consistent}, the kernel spaces of $\mathcal{L}^{\lambda,T}$ and $\tilde{\mathcal{L}}^{\lambda,T}$ are
\ben
\label{null-based-on-mathcal-N}
\ker \mathcal{L}^{\lambda,T} =  \mathrm{span} \{ \mathcal{N}_{\lambda,T}, \mathcal{N}_{\lambda,T}v_{1}, \mathcal{N}_{\lambda,T}v_{2},
\mathcal{N}_{\lambda,T}v_{3}, \mathcal{N}_{\lambda,T}|v|^{2} \},
\\
\label{kernel-space}
\ker \tilde{\mathcal{L}}^{\lambda,T} = \mathrm{span} \{N_{\lambda}, N_{\lambda}v_{1}, N_{\lambda}v_{2}, N_{\lambda}v_{3}, N_{\lambda}|v|^{2}\}.
\een
Observe that $\ker \tilde{\mathcal{L}}^{\lambda,T}$ depends only on the parameter $\lambda$, while
$\ker \mathcal{L}^{\lambda,T}$ depends on the two parameters $\lambda,T$.

For notational simplicity, for $k \geq 0$,
we denote the $k$-th moments of the density $M_{\lambda} (1+M_{\lambda}) = N_{\lambda}^{2}$ by $m_{k}$. More precisely,
\ben \label{moment-k-of-N-square}
m_{k} := \int_{\mathbb{R}^{3}} |v|^{k} M_{\lambda}(v) (1+M_{\lambda}(v)) \mathrm{d}v = \int_{\mathbb{R}^{3}} |v|^{k} N_{\lambda}^{2}(v) \mathrm{d}v.
\een

\begin{rmk}
Note that some of the above moments is infinite if $\lambda =0$ and $k$ is small. Indeed, if $\lambda =0$, then near $|v|=0$, it holds that
\beno
M_{\lambda}(v) \sim |v|^{-2}, \quad M_{\lambda}(v)(1+M_{\lambda}(v)) \sim |v|^{-4}.
\eeno
We remark that $m_{0} = |N_{\lambda}|_{L^{2}}^{2} = \infty$ if $\lambda= 0$. This fact means the $L^{2}$ framework may not be suitable for the critical temperature $\lambda= 0$.
\end{rmk}

In the following Lemma, we give some estimates of $m_{k}$.
\begin{lem} \label{N-lambda-l2-norm}
We have for $k \geq 1$,
\ben \label{L2-of-N}
m_{0} = |N_{\lambda}|_{L^{2}}^{2} \sim e^{-\lambda} (1-e^{-\lambda})^{-1/2}.
\\
\label{weighted-L2-k-geq-1-of-N}
m_{2k} = ||\cdot|^{k}N_{\lambda}|_{L^{2}}^{2} \sim C_{k}e^{-\lambda}.
\een
Here $C_{k} \sim 1$ for $1 \leq k \leq 3$.
\end{lem}
\begin{proof} Recalling \eqref{definition-of-M-lambda}, we have
\ben \label{N-lambda-l2-evaluation}
m_{0} = |N_{\lambda}|_{L^{2}}^{2} = \int  \f{\exp(\f{|v|^{2}}{2}+\lambda)}{ (\exp(\f{|v|^{2}}{2}+\lambda) -1)^{2}} \mathrm{d}v
= \exp(-\lambda) \int  \f{\exp(-\f{|v|^{2}}{2})}{ (1 -\exp(-\f{|v|^{2}}{2}-\lambda))^{2}} \mathrm{d}v \colonequals  \exp(-\lambda)h(\lambda)
\een
For $\lambda \geq 1$, we have
\beno
\int   \exp(-\f{|v|^{2}}{2}) \mathrm{d}v \leq h(\lambda) \leq
(1 -\exp(-1))^{-2} \int   \exp(-\f{|v|^{2}}{2}) \mathrm{d}v .
\eeno
That is
\ben \label{h-lambda-geq-1}
h(\lambda) \sim 1
\een

Now consider $\lambda \leq 1$.
For $|v| \geq 1$, we have
\beno
\int \mathrm{1}_{|v| \geq 1}  \f{\exp(-\f{|v|^{2}}{2})}{ (1 -\exp(-\f{|v|^{2}}{2}-\lambda))^{2}} \mathrm{d}v \lesssim
\int \mathrm{1}_{|v| \geq 1}  \f{\exp(-\f{|v|^{2}}{2})}{ (1 -\exp(-\f{1}{2}))^{2}} \mathrm{d}v \lesssim 1 \leq \lambda^{-\f{1}{2}}.
\eeno
For $|v| \leq 1$, we have $1 -\exp(-\f{|v|^{2}}{2}-\lambda) \sim \f{|v|^{2}}{2}+\lambda$ and so
\beno
\int \mathrm{1}_{|v| \leq 1}  \f{\exp(-\f{|v|^{2}}{2})}{ (1 -\exp(-\f{|v|^{2}}{2}-\lambda))^{2}} \mathrm{d}v \sim
\int \mathrm{1}_{|v| \leq 1}  \f{1}{ (\f{|v|^{2}}{2}+\lambda)^{2}} \mathrm{d}v.
\eeno
Note that
\beno
\int \mathrm{1}_{|v| \leq 1}  \f{1}{ (\f{|v|^{2}}{2}+\lambda)^{2}} \mathrm{d}v = \lambda^{-\f{1}{2}}\int \mathrm{1}_{|v| \leq \lambda^{-\f{1}{2}}}  \f{1}{ (\f{|v|^{2}}{2}+1)^{2}} \mathrm{d}v \sim \lambda^{-\f{1}{2}}
\eeno
By these estimates, we find
\ben \label{h-lambda-leq-1}
h(\lambda) \sim \lambda^{-\f{1}{2}}.
\een
Patching together \eqref{h-lambda-geq-1} and \eqref{h-lambda-leq-1}, we conclude that
\beno
h(\lambda) \sim \max\{\lambda^{-\f{1}{2}},1\} \sim (1-e^{-\lambda})^{-1/2}.
\eeno
Recalling \eqref{N-lambda-l2-evaluation}, we arrive at \eqref{L2-of-N}.

Similarly to \eqref{N-lambda-l2-evaluation}
\ben \label{N-lambda-weighted-l2-evaluation}
||\cdot|^{k}N_{\lambda}|_{L^{2}}^{2} = \int  \f{\exp(\f{|v|^{2}}{2}+\lambda)}{ (\exp(\f{|v|^{2}}{2}+\lambda) -1)^{2}} \mathrm{d}v
= \exp(-\lambda) \int  \f{|v|^{2k}\exp(-\f{|v|^{2}}{2})}{ (1 -\exp(-\f{|v|^{2}}{2}-\lambda))^{2}} \mathrm{d}v \colonequals  \exp(-\lambda)h_{k}(\lambda)
\een
For $|v| \geq 1$, we have
\beno
\int  \mathrm{1}_{|v| \geq 1} \f{|v|^{2k}\exp(-\f{|v|^{2}}{2})}{ (1 -\exp(-\f{|v|^{2}}{2}-\lambda))^{2}} \mathrm{d}v
 \sim \int  \mathrm{1}_{|v| \geq 1} |v|^{2k}\exp(-\f{|v|^{2}}{2}) \mathrm{d}v \sim C_{k}.
\eeno
For $|v| \leq 1$ and $k \geq 1$,
we have
\beno
\int  \mathrm{1}_{|v| \leq 1} \f{|v|^{2k}\exp(-\f{|v|^{2}}{2})}{ (1 -\exp(-\f{|v|^{2}}{2}-\lambda))^{2}} \mathrm{d}v
 \leq \int  \mathrm{1}_{|v| \leq 1} \f{|v|^{2}\exp(-\f{|v|^{2}}{2})}{ (1 -\exp(-\f{|v|^{2}}{2}))^{2}} \mathrm{d}v \sim
 \int  \mathrm{1}_{|v| \leq 1} |v|^{-2} \sim 1.
\eeno
As a result
\beno
h_{k}(\lambda) \sim C_{k}.
\eeno
Recalling \eqref{N-lambda-weighted-l2-evaluation}, we arrive at \eqref{weighted-L2-k-geq-1-of-N}.
\end{proof}

We construct orthogonal basis for $\ker \mathcal{L}^{\lambda,T}$ and $\ker \tilde{\mathcal{L}}^{\lambda,T}$ for $\lambda, T > 0$ as follows
\ben \label{definition-of-di}
\{ d^{\lambda,T}_{i} \}_{1 \leq i \leq 5} \colonequals     \{ \mathcal{N}_{\lambda,T}, \mathcal{N}_{\lambda,T} v_{1}/T, \mathcal{N}_{\lambda,T} v_{2}/T, \mathcal{N}_{\lambda,T} v_{3}/T, \mathcal{N}_{\lambda,T}(|v|^{2}/T - \tilde{C}_{\lambda}) \}.
\\
\label{definition-of-di-tilde}
\{ \tilde{d}^{\lambda}_{i} \}_{1 \leq i \leq 5} \colonequals     \{ N_{\lambda}, N_{\lambda} v_{1}, N_{\lambda} v_{2}, N_{\lambda} v_{3}, N_{\lambda}(|v|^{2} - \tilde{C}_{\lambda}) \},
\een
where $\tilde{C}_{\lambda} \colonequals  m_{2}/m_{0}$ which depends only on $\lambda$.
Note that $\langle \mathcal{N}_{\lambda,T}|v|^{2} , \mathcal{N}_{\lambda,T}\rangle|\mathcal{N}_{\lambda,T}|^{-2}_{L^{2}}\mathcal{N}_{\lambda,T}$ is the projection of $\mathcal{N}_{\lambda,T}|v|^{2}$ on $\mathcal{N}_{\lambda,T}$.
By normalization, we get orthonormal basis of $\ker \mathcal{L}^{\lambda,T}$ and $\ker \tilde{\mathcal{L}}^{\lambda,T}$
\ben \label{definition-of-e-pm-rho}
\{ e^{\lambda,T}_{i} \}_{1 \leq i \leq 5} \colonequals     \{ \frac{d^{\lambda,T}_{i}}{|d^{\lambda,T}_{i}|_{L^{2}}} \}_{1 \leq i \leq 5}, \quad \{ \tilde{e}^{\lambda}_{i} \}_{1 \leq i \leq 5} \colonequals     \{ \frac{\tilde{d}^{\lambda}_{i}}{|\tilde{d}^{\lambda}_{i}|_{L^{2}}} \}_{1 \leq i \leq 5}.
\een
With these orthonormal basis, the projection $\mathbb{P}_{\lambda,T}$ on  $\ker \mathcal{L}^{\lambda,T}$ and $\tilde{\mathbb{P}}_{\lambda}$ on  $\ker \tilde{\mathcal{L}}^{\lambda,T}$
are defined by
\ben \label{definition-projection-operator}
\mathbb{P}_{\lambda,T}f  \colonequals     \sum_{i=1}^{5} \langle f, e^{\lambda,T}_{i}\rangle e^{\lambda,T}_{i},
\quad
\tilde{\mathbb{P}}_{\lambda}f \colonequals     \sum_{i=1}^{5} \langle f, \tilde{e}^{\lambda}_{i}\rangle \tilde{e}^{\lambda}_{i}.
\een

Recalling \eqref{definition-of-M-lambda}, it is elementary to check that
\ben \label{projection-operator-before-after-scaling}
A_{T^{1/2}} \mathbb{P}_{\lambda,T} f = \tilde{\mathbb{P}}_{\lambda} A_{T^{1/2}}f.
\een

Let us see $\tilde{\mathbb{P}}_{\lambda}$ more clearly.
By \eqref{definition-projection-operator} and rearrangement, we have
\ben \nonumber
\tilde{\mathbb{P}}_{\lambda}f &=& \langle f, \frac{N}{|N|_{L^{2}}} \rangle \frac{N}{|N|_{L^{2}}} + \sum_{i=1}^{3} \langle f, \frac{N v_{i}}{|N v_{i}|_{L^{2}}} \rangle \frac{N v_{i}}{|N v_{i}|_{L^{2}}} + \langle f, \frac{N(|v|^{2}-\tilde{C}_{\lambda})}{ |N(|v|^{2}-\tilde{C}_{\lambda})|_{L^{2}} } \rangle \frac{N(|v|^{2}-\tilde{C}_{\lambda})}{|N(|v|^{2}-\tilde{C}_{\lambda})|_{L^{2}}}
\\ \label{linear-combination-of-basis} &=& (a^{f}_{\lambda} + b^{f}_{\lambda} \cdot v + c^{f}_{\lambda}|v|^{2})N,
\een
where
\beno
a^{f}_{\lambda} \colonequals     \langle f, (\frac{1}{|N|_{L^{2}}^{2}} + \frac{\tilde{C}_{\lambda}^{2}}{|N(|v|^{2}-\tilde{C}_{\lambda})|_{L^{2}}^{2}})N -  \frac{\tilde{C}_{\lambda}}{|N(|v|^{2}-\tilde{C}_{\lambda})|_{L^{2}}^{2}}N|v|^{2} \rangle, \quad
b^{f}_{\lambda} \colonequals     \langle f, \frac{N v}{|N v_{i}|_{L^{2}}^{2}} \rangle,
\\
c^{f}_{\lambda} \colonequals     \langle f, \frac{1}{|N(|v|^{2}-\tilde{C}_{\lambda})|_{L^{2}}^{2}} N|v|^{2} - \frac{\tilde{C}_{\lambda}}{|N(|v|^{2}-\tilde{C}_{\lambda})|_{L^{2}}^{2}} N \rangle.
\eeno
Note that $b^{f}_{\lambda}$ is a vector of length 3.  Note that
\beno
\tilde{C}_{\lambda} = \f{m_{2}}{m_{0}}, \quad |N|_{L^{2}}^{2} = m_{0}, \quad |N v_{i}|_{L^{2}}^{2} = \f{1}{3}m_{2}, \quad |N(|v|^{2}-\tilde{C}_{\lambda})|_{L^{2}}^{2} = \f{m_{0}m_{4}-m_{2}^{2}}{m_{0}}
\eeno
Let us define
\ben \label{defintion-of-l-i}
l_{\lambda,1}\colonequals    \f{m_{4}}{m_{0}m_{4}-m_{2}^{2}}, \quad
l_{\lambda,2}\colonequals    \f{m_{2}}{m_{0}m_{4}-m_{2}^{2}}, \quad
l_{\lambda,3}\colonequals    \frac{3}{m_{2}}, \quad
l_{\lambda,4}\colonequals    \f{m_{0}}{m_{0}m_{4}-m_{2}^{2}}.
\een
For simplicity, let $l_{i} = l_{\lambda, i}$.
Then there holds
\ben \label{explicit-defintion-of-abc}
a^{f}_{\lambda} =  \langle f, l_{1}N -  l_{2}N|v|^{2} \rangle, \quad
b^{f}_{\lambda} = \langle f, l_{3}N v \rangle, \quad
c^{f}_{\lambda} = \langle f, l_{4} N|v|^{2} - l_{2} N \rangle.
\een

For simplicity, let $f_{1} = \tilde{\mathbb{P}}_{\lambda}f, f_{2} = f- \tilde{\mathbb{P}}_{\lambda}f$, and $\mathcal{A} = (a^{f}_{\lambda}, b^{f}_{\lambda}, c^{f}_{\lambda})$.
By Lemma \ref{N-lambda-l2-norm}, we have
\ben \label{f-1-is-bounded-by-abc}
|\partial^{\alpha}f_{1}|_{L^{2}_{1/2}} \lesssim e^{-\lambda/2} (1-e^{-\lambda})^{-1/4} |\partial^{\alpha} \mathcal{A}|.
\een
As a result, we can choose $C_{0}$ is large enough such that $C_{0}e^{-\lambda} (1-e^{-\lambda})^{-1/2}
|\nabla_{x} \mathcal{A}|_{H^{N-1}_{x}}^{2} \geq 
\|\nabla_{x} f_{1}\|_{H^{N-1}_{x}L^{2}_{1/2}}^{2}$. Then we define the dissipation functional $\mathcal{D}_{N}(\cdot)$ as
\ben \label{defintion-of-dissipation-after-scaling}
\mathcal{D}_{N}(f) \colonequals
C_{0}e^{-\lambda} (1-e^{-\lambda})^{-1/2}|\nabla_{x} \mathcal{A}|_{H^{N-1}_{x}}^{2} + \|f_{2}\|_{H^{N}_{x}L^{2}_{1/2}}^{2}.
\een
Note that
\ben \label{D-N-lower-bound}
 \mathcal{D}_{N}(f)\geq \|\nabla_{x} f_{1}\|_{H^{N-1}_{x}L^{2}_{1/2}}^{2} + \|f_{2}\|_{H^{N}_{x}L^{2}_{1/2}}^{2}.
\een

Global well-posedness of \eqref{T-half-scaled-more-consistent} is presented as the following Theorem which is sufficient to derive Theorem \ref{global-well-posedness-original}.
\begin{thm}\label{global-well-posedness} Let $0<\epsilon \leq 1$. Let $\lambda, T >0$. Let $N \geq 2$. Let
\ben \label{tilde-C-star}
\tilde{C}_{*}(\lambda,T) \colonequals  e^{-2\lambda} (1-e^{-\lambda})^{\f{33}{2}} \min\{T^{-3},1\}.
\een
There exist a universal constant  $\delta_{*}>0$ such that if
\ben \label{condition-on-initial-data} 
M_{\lambda} 
+ N_{\lambda} \tilde{f}_{0} \geq 0, \quad \|\tilde{f}_{0}\|_{H^{2}_{x}L^{2}}^{2} \leq \delta_{*} \tilde{C}_{*}(\lambda,T), \quad \|\tilde{f}_{0}\|_{H^{N}_{x}L^{2}}<\infty, \een
then the Cauchy problem \eqref{T-half-scaled-more-consistent} with initial datum $\tilde{f}_{0}$
has a unique global solution $\tilde{f} = \tilde{f}^{\lambda,T}_{\epsilon} \in L^{\infty}([0,\infty); H^{N}_{x}L^{2})$ satisfying $M_{\lambda} + N_{\lambda} \tilde{f}(t) \geq 0$ and
\ben \label{uniform-estimate-global}
\sup_{t}\|\tilde{f}(t)\|_{H^{N}_{x}L^{2}}^{2}+ \f{1}{K(\lambda,T)}  \int_{0}^{\infty} \mathcal{D}_{N}(\tilde{f}) \mathrm{d}\tau + \f{\mathrm{C}_{1}(\lambda,T)}{K(\lambda,T)}  \f{1}{\epsilon^{2}} \int_{0}^{\infty} \|\tilde{f}_{2}\|^{2}_{H^{N}_{x}L^{2}_{1/2}} \mathrm{d}\tau \leq  P_{N}(\tilde{f}_{0})\|f_0\|_{H^{N}_{x}L^{2}}^{2}.
\een
where the formula of $P_{N}$ is  explicitly given \eqref{definition-of-pn} in Theorem \ref{a-priori-estimate-LBE}. Here the constants $K(\lambda,T), \mathrm{C}_{1}(\lambda,T)$ are defined defined in \eqref{defintion-of-K-lambda-T}, \eqref{defintion-of-C1-lambda-T}.
\end{thm}

\begin{proof}[Proof of Theorem \ref{global-well-posedness-original}] Let $f_0$ be an initial datum satisfying \eqref{condition-on-initial-data-original}. Then $\tilde{f}_{0} = A_{T^{1/2}} f_0$ is an initial datum satisfying \eqref{condition-on-initial-data}. Then by Theorem \ref{global-well-posedness}, the Cauchy problem \eqref{T-half-scaled-more-consistent} with initial datum $\tilde{f}_{0}$ has a unique global solution $\tilde{f}$ satisfying \eqref{uniform-estimate-global}. By Proposition \ref{solution-relation}, $f = A_{T^{-1/2}} \tilde{f}$ is a solution to  \eqref{linearized-quantum-Boltzmann-eq} with the initial datum $f_{0}$.
	Recalling \eqref{projection-operator-before-after-scaling} and using the following identity
	\ben \label{inner-product-relation}
	a^{3} \langle A_{a}g, A_{a}h \rangle =  \langle g, h \rangle,
	\een
	we have
	\beno
		\|\tilde{f}\|_{H^{N}_{x}L^{2}}^{2} = T^{-\f{3}{2}} \|f\|_{H^{N}_{x}L^{2}}^{2}, \quad 
	\|\nabla_{x} \tilde{\mathbb{P}}_{\lambda} \tilde{f}\|_{H^{N-1}_{x}L^{2}}^{2} = T^{-\f{3}{2}} \|\nabla_{x} \mathbb{P}_{\lambda,T} f\|_{H^{N-1}_{x}L^{2}}^{2},
	\\
	\|\tilde{f} - \tilde{\mathbb{P}}_{\lambda} \tilde{f}\|_{H^{N}_{x}L^{2}}^{2} = T^{-\f{3}{2}} \|f - \mathbb{P}_{\lambda,T} f\|_{H^{N}_{x}L^{2}}^{2}, \quad
	\||v|^{\f{1}{2}}(\tilde{f} - \tilde{\mathbb{P}}_{\lambda} \tilde{f})\|_{H^{N}_{x}L^{2}}^{2} = T^{-\f{3}{2}} T^{-\f{1}{2}} \||v|^{\f{1}{2}}(f - \mathbb{P}_{\lambda,T} f)\|_{H^{N}_{x}L^{2}}^{2},
	\eeno
	which gives
	\ben \label{scaling-norm-relation}
\quad \mathcal{D}_{N}(\tilde{f})  \geq T^{-\f{3}{2}} \mathcal{D}_{N,T}(f), \quad \|\tilde{f}_{2}\|^{2}_{H^{N}_{x}L^{2}_{1/2}}  \geq T^{-\f{3}{2}} \mathcal{D}_{N,T,2}(f).
	\een
	From these relations, we obtain all the estimates on $f$ in Theorem \ref{global-well-posedness-original} and finish the proof.
\end{proof}

Hydrodynamics limit of \eqref{T-half-scaled-more-consistent} is presented as the following Theorem which is sufficient to derive Theorem \ref{thm-hydrodynamical-limit-ori}.
\begin{thm} \label{thm-hydrodynamical-limit}
Recall the constant $\delta_{*}$ in Theorem \ref{global-well-posedness}.
Let $0<\epsilon < 1$. Let $\lambda, T >0$. Let $N \geq 2$.
Let $\{f^{\epsilon}_{0}\}_{0<\epsilon < 1}$ be a family of initial datum satisfying $M_{\lambda} + N_{\lambda}f^{\epsilon}_{0} \geq 0$ and
\ben \label{small-condition}
\sup_{0<\epsilon<1} \|f^{\epsilon}_{0}\|_{H^{2}_{x}L^{2}}^{2} \leq \delta_{*} \tilde{C}_{*}(\lambda,T), \quad
M_{0} \colonequals  \sup_{0<\epsilon<1} \|f^{\epsilon}_{0}\|_{H^{N}_{x}L^{2}} < \infty,
\\ \label{initial-convergence}
\mathbb{P}_{\lambda}f^{\epsilon}_{0} \to f_{0} = (\rho_{0} + u_{0} \cdot v + \theta_{0} (\f{|v|^{2}}{2} - K_{\lambda})) N_{\lambda} \text{ as }
\epsilon \to 0,
\text{ strongly in } H^{N}_{x}L^{2},
\een
for some $(\rho_{0}, u_{0}, \theta_{0}) \in H^{N}_{x}$ with $\rho_{0} + \theta_{0} = 0$. Here $K_{\lambda} = K_{A} -1$ is a constant depending only on $\lambda$. Let $f^{\epsilon}$ be the solution to the Cauchy problem \eqref{T-half-scaled-more-consistent}
 with initial datum $f^{\epsilon}_{0}$. Then there is a subsequence of $\{f^{\epsilon}\}$ still denoting it by $\{f^{\epsilon}\}$ such that,
 \ben \label{weakly-star-convergence-to-some-special-form}
 f^{\epsilon} \to (\rho + u \cdot v + \theta (\f{|v|^{2}}{2} - K_{\lambda})) N_{\lambda} \text{ as }
\epsilon \to 0, \text{ weakly-* in } L^{\infty}(\mathbb{R}_{+}; H^{N}_{x}L^{2}),
 \een
 for some $(\rho, u, \theta) \in L^{\infty}(\mathbb{R}_{+}; H^{N}_{x}) \cap C(\mathbb{R}_{+}; H^{N-1}_{x})$. Moreover $(\rho, u, \theta)$ is a weak solution of \eqref{NSF-problem}. In addition,
 \ben \label{momoent-1}
\f{3}{m_{2}}\mathcal{P}\langle f^{\epsilon}, v N_{\lambda} \rangle \to u \text{ strongly in } C(\mathbb{R}_{+}; H^{N-1}_{x}) \text{ and weakly-* in } L^{\infty}(\mathbb{R}_{+}; H^{N}_{x}),
\\ \label{momoent-2}
\f{1}{C_{A}} \langle f^{\epsilon}, (\f{|v|^{2}}{2} - K_{A}) N_{\lambda} \rangle \to \theta \text{ strongly in } C(\mathbb{R}_{+}; H^{N-1}_{x}) \text{ and weakly-* in } L^{\infty}(\mathbb{R}_{+}; H^{N}_{x}).
 \een
\end{thm}

\begin{proof}[Proof of Theorem \ref{thm-hydrodynamical-limit-ori}]
By Proposition \ref{solution-relation} and Theorem \ref{thm-hydrodynamical-limit}, using the formula \eqref{inner-product-relation},
 we get Theorem \ref{thm-hydrodynamical-limit-ori}.
\end{proof}

Now in the rest of the article it remains to derive Theorem \ref{global-well-posedness} and Theorem \ref{thm-hydrodynamical-limit} on the Cauchy problem \eqref{T-half-scaled-more-consistent}.

\section{Linear and nonlinear collision operators analysis} \label{operator-analysis}
In the rest of the article, in the various functional estimates, the involved functions $g, h, \varrho, f$ are assumed to be functions on $\R^{3}$ or  $\mathbb{R}^{3} \times \R^{3}$ such that the corresponding norms of them
are well-defined. For simplicity, we use the notation $\mathrm{d}V \colonequals  \mathrm{d}\sigma \mathrm{d}v_{*} \mathrm{d}v$.

\subsection{Coercivity estimate}
 We first give some basic properties of $M_{\lambda}$ and $N_{\lambda}$ defined in \eqref{definition-of-M-lambda}.
\begin{lem} \label{M-N-mu} For simplicity, let $\mu := \exp(-\f{|v|^{2}}{2})$.
For $\lambda>0$, it holds that
\ben
\label{M-rho-mu-N}
e^{-\lambda }\mu \leq M_{\lambda} \leq \f{e^{-\lambda}\mu}{1-e^{-\lambda}}, \quad  e^{-\lambda/2} \mu^{\f{1}{2}} \leq  N_{\lambda} \leq \f{e^{-\lambda/2}\mu^{\f{1}{2}}}{1-e^{-\lambda}}.
\een
As a direct result, since $\mu\mu_{*}=\mu^{\prime}\mu^{\prime}_{*}$, it holds that
\ben \label{K-2-mu}
e^{-2\lambda} \mu \mu_{*} \leq  N_{\lambda} (N_{\lambda})_{*} (N_{\lambda})^{\prime} (N_{\lambda})^{\prime}_{*} \leq (1-e^{-\lambda})^{-4}e^{-2\lambda} \mu \mu_{*}.
\een
It holds that
\ben \label{comparing-NN-with-NN-prime}
1-e^{-\lambda} \leq \f{N_{\lambda} (N_{\lambda})_{*}}{(N_{\lambda})^{\prime} (N_{\lambda})^{\prime}_{*}} \leq (1-e^{-\lambda})^{-1}.
\een
\end{lem}
\begin{proof}
We only prove \eqref{comparing-NN-with-NN-prime} as
the other two results are obvious. Since $|v|^{2}+|v_{*}|^{2} = |v^{\prime}|^{2}+|v^{\prime}_{*}|^{2}$,
\beno
\f{N_{\lambda} (N_{\lambda})_{*}}{(N_{\lambda})^{\prime} (N_{\lambda})^{\prime}_{*}} = \f{(\exp(\f{|v^{\prime}|^{2}}{2}+\lambda) -1) (\exp(\f{|v^{\prime}_{*}|^{2}}{2}+\lambda) -1)}{(\exp(\f{|v|^{2}}{2}+\lambda) -1) (\exp(\f{|v_{*}|^{2}}{2}+\lambda) -1)}.
\eeno
Now it suffices to consider for some $s>1$ the following quantity
\beno
\f{(w_{1}s -1) (w_{2}s -1)}{(w_{3}s -1) (w_{4}s -1)},
\eeno
subject to $w_{1}w_{2}=w_{3}w_{4}, w_{1}, w_{2}, w_{3}, w_{4} \geq 1$. Let $k^{2} = w_{1}w_{2}=w_{3}w_{4}$ for some $k \geq 1$, then the numerator and denominator 
enjoy the same bounds $(k^{2}s -1) (s -1) \leq (w_{1}s -1) (w_{2}s -1), (w_{3}s -1) (w_{4}s -1)  \leq (ks -1)^{2}$. Therefore,
\beno
\f{(k^{2}s -1) (s -1)}{(ks -1)^{2}} \leq \f{(w_{1}s -1) (w_{2}s -1)}{(w_{3}s -1) (w_{4}s -1)} \leq \f{(ks -1)^{2}}{(k^{2}s -1) (s -1)}
\eeno
It is easy to see $\f{(ks -1)^{2}}{(k^{2}s -1) (s -1)}$ is increasing w.r.t. $k$ and so achieves its maximum $\f{s}{s -1}$ at $k \to \infty$. Therefore
\beno
\f{s -1}{s} \leq \f{(w_{1}s -1) (w_{2}s -1)}{(w_{3}s -1) (w_{4}s -1)} \leq \f{s}{s -1}.
\eeno
As a result,
by taking $s = e^{\lambda}$, we get \eqref{comparing-NN-with-NN-prime}.
\end{proof}

Recall that the classical linearized Boltzmann operator $\mathcal{L}$ with hard sphere kernel is  defined by
\ben \label{linearized-Boltzmann-operator-hard}
(\mathcal{L}f)(v)  \colonequals      \int  |v-v_{*}| \mu_{*} \mu^{\f{1}{2}}
\mathrm{S}(\mu^{-\f{1}{2}}f)
\mathrm{d}\sigma \mathrm{d}v_{*},
\een
Define the functional $\mathcal{H} (\cdot)$ as
\ben
\quad \mathcal{H} (f) = \int |v-v_{*}| \mu \mu_{*}
 \mathrm{S}^{2}( \mu^{-1/2} f ) \mathrm{d}V.
\een
The coercivity estimate of $\mathcal{L}$ is that
for $f \in (\ker \mathcal{L})^{\perp}$, 
\ben \label{coercivity-of-classical-operator}
C_{1} |f|_{L^{2}_{1/2}}^{2} \geq \mathcal{H} (f) =  4  \langle \mathcal{L}f, f \rangle \geq C_{0} |f|_{L^{2}_{1/2}}^{2},
\een
where $0<C_{0}<C_{1}$ are two universal constants.

We are ready to get the key coercivity estimate of $\tilde{\mathcal{L}}^{\lambda,T}$ by using \eqref{coercivity-of-classical-operator}. 
\begin{thm} \label{coercivity-estimate}  Let $\lambda, T >0$. Recall \eqref{kernel-space}.
For $f \in (\ker \tilde{\mathcal{L}}^{\lambda,T})^{\perp}$, it holds that
\ben \label{lower-and-upper-bound}
 \tilde{C}_{1,\lambda,T}|f|_{L^{2}_{1/2}}^{2} \geq \langle \tilde{\mathcal{L}}^{\lambda,T}f, f \rangle \geq  C_{1,\lambda,T} |f|_{L^{2}_{1/2}}^{2},
\een
where
\ben \label{constant-of-C-1}
 C_{1,\lambda,T} \colonequals   T^{2} \f{C_{0} e^{-\lambda} (1-e^{-\lambda})^{5/2}}{C_{2}}, \quad \tilde{C}_{1,\lambda,T} \colonequals  T^{2} C_{3} C_{1}   (1-e^{-\lambda})^{-4} e^{-\lambda},
\een
for some universal constant $C_{0}, C_{1}, C_{2}, C_{3}$. Safely speaking $\f{1}{100} \leq C_{0} \leq 1 \leq C_{1}, C_{2}, C_{3} \leq 100$. Here $C_{0}, C_{1}$ are the constants appearing in \eqref{coercivity-of-classical-operator}.
\end{thm}
\begin{proof} Note that
\ben \label{inner-product}
\langle \tilde{\mathcal{L}}^{\lambda,T}f, f \rangle = \frac{1}{4} \int B_{T} N_{\lambda} (N_{\lambda})_{*} (N_{\lambda})^{\prime} (N_{\lambda})^{\prime}_{*}
 \mathrm{S}^{2}( N_{\lambda}^{-1} f ) \mathrm{d}V.
\een
Thanks to \eqref{K-2-mu}, we have
\ben \label{inner-product-2-mu-mu}
\frac{1}{4} e^{-2\lambda} T^{2} \mathcal{J}_{\lambda} (f) \leq
\langle \tilde{\mathcal{L}}^{\lambda,T}f, f \rangle \leq  \frac{1}{4} (1-e^{-\lambda})^{-4} e^{-2\lambda} T^{2} \mathcal{J}_{\lambda} (f),
\een
where we define
\ben \label{defintion-J-and-H}
\mathcal{J}_{\lambda} (f) = \int |v-v_{*}| \mu \mu_{*}
 \mathrm{S}^{2}( N_{\lambda}^{-1} f ) \mathrm{d}V.
\een
We now study $\mathcal{J}_{\lambda} (f)$ using \eqref{coercivity-of-classical-operator}. We now relate $\mathcal{J}_{\lambda} (\cdot)$ to $\mathcal{H}(\cdot)$.
	For a function $f$, we define
%$w_{f}, \varPhi_{f}$ by
\ben \label{def-w-f-psi-f}
w_{f} \colonequals N_{\lambda}\mu^{-\f12} \mathbb{P}_{0} (N_{\lambda}^{-1}\mu^{\f12}f) 
, \quad \varPhi_{f}  \colonequals 
(f - w_{f}) N_{\lambda}^{-1} \mu^{1/2},
\een
where $\mathbb{P}_{0}$ is the projection on 
$\ker \mathcal{L}$.
It is straightforward to check for $f \in (\ker \tilde{\mathcal{L}}^{\lambda,T})^{\perp}$ that 
\ben \label{functions-in-the-spaces}
w_{f} \in \ker \tilde{\mathcal{L}}^{\lambda,T}
,  \quad \varPhi_{f} =  N_{\lambda}^{-1}\mu^{\f12} f - \mathbb{P}_{0} (N_{\lambda}^{-1}\mu^{\f12}f) \in (\ker \mathcal{L})^{\perp}.
\een
By the above construction, for $f \in (\ker \tilde{\mathcal{L}}^{\lambda,T})^{\perp}$,
it is easy to see
\ben \label{key-relation-btw-quantum-and-classical}
\mathcal{J}_{\lambda} (f) = \mathcal{J}_{\lambda} (f - w_{f}) = \mathcal{J}_{\lambda} (N_{\lambda} \mu^{-1/2} \varPhi_{f})
= \mathcal{H} (\varPhi_{f}).
\een

	By \eqref{key-relation-btw-quantum-and-classical},
\eqref{functions-in-the-spaces},  using \eqref{coercivity-of-classical-operator}, we have
\ben \label{J-rho-lower-bound}
\mathcal{J}_{\lambda} (f) \geq C_0 |\varPhi_{f}|_{L^{2}_{1/2}}^{2} \geq (1-e^{-\lambda})^{2}e^{\lambda} C_{0} |\varPhi_{f} N_{\lambda} \mu^{-1/2}|_{L^{2}_{1/2}}^{2} = (1-e^{-\lambda})^{2}e^{\lambda} C_{0} |f - w_{f}|_{L^{2}_{1/2}}^{2}.
\een
Note that $f \perp w_{f}$, so we have
\ben \label{lower-l2}
|f - w_{f}|_{L^{2}_{1/2}}^{2} \geq |f - w_{f}|_{L^{2}}^{2} = |f|_{L^{2}}^{2} + |w_{f}|_{L^{2}}^{2} \geq |f|_{L^{2}}^{2}.
\een
We also have
\beno
|f - w_{f}|_{L^{2}_{1/2}}^{2} \geq \f{1}{2}|f|_{L^{2}_{1/2}}^{2} - |w_{f}|_{L^{2}_{1/2}}^{2}
\eeno
Note that  $w_{f} = (a_{f} + b_{f} \cdot v + c_{f} |v|^{2}) N_{\lambda}$ where $a_{f}, b_{f}, c_{f}$ are the constants given by
\beno
a_{f} = \f{5}{2} \langle f, N_{\lambda}^{-1} \mu \rangle - \f{1}{2} \langle f, |v|^{2}N_{\lambda}^{-1} \mu \rangle, \quad
b_{f} = \langle f, v N_{\lambda}^{-1} \mu \rangle, \quad
c_{f} =   \f{1}{6} \langle f, |v|^{2}N_{\lambda}^{-1} \mu \rangle - \f{1}{2} \langle f, N_{\lambda}^{-1} \mu \rangle.
\eeno
It is easy to see
\ben \label{up-of-a-b-c}
|a_{f}| + |b_{f}| + |c_{f}| \lesssim e^{\lambda/2} |f|_{L^{2}}.
\een
By Lemma \ref{N-lambda-l2-norm} and the estimate \eqref{up-of-a-b-c},
\beno
|w_{f}|_{L^{2}_{1/2}} \leq (|a_{f}| + |b_{f}| + |c_{f}|) (|(1+|v|^{2})N_{\lambda}|_{L^{2}_{1/2}}) \lesssim (1-e^{-\lambda})^{-1/4}|f|_{L^{2}}.
\eeno
Therefore for some universal constant $C_{2} \geq 1$, we have
\ben \label{lower-l2-weight}
|f - w_{f}|_{L^{2}_{1/2}}^{2} \geq \f{1}{2}|f|_{L^{2}_{1/2}}^{2} - C_{2} (1-e^{-\lambda})^{-1/2}|f|_{L^{2}}^{2}.
\een
Making a suitable combination between \eqref{lower-l2}
and \eqref{lower-l2-weight},
we get
\beno
|f - w_{f}|_{L^{2}_{1/2}}^{2} \geq \f{1}{2(1+ C_{2}(1-e^{-\lambda})^{-1/2})}|f|_{L^{2}_{1/2}}^{2} \geq \f{1}{4 C_{2}(1-e^{-\lambda})^{-1/2}}|f|_{L^{2}_{1/2}}^{2}.
\eeno
Recalling \eqref{inner-product-2-mu-mu} and \eqref{J-rho-lower-bound}, we arrive at
\beno
\langle \tilde{\mathcal{L}}^{\lambda,T}f, f \rangle \geq T^{2} \f{C_{0} e^{-\lambda} (1-e^{-\lambda})^{5/2}}{16C_{2}} |f|_{L^{2}_{1/2}}^{2}
\eeno

By \eqref{key-relation-btw-quantum-and-classical},
\eqref{functions-in-the-spaces}, using the upper bound in \eqref{coercivity-of-classical-operator}, by noting  $N_{\lambda}^{-1} \mu^{1/2} \leq e^{\lambda/2}$,
we can easily get
\ben \label{upper-bound-by-norm}
\mathcal{J}_{\lambda} (f) \leq C_{1} |\varPhi_{f}|_{L^{2}_{1/2}}^{2} 
\lesssim  C_{1} |N_{\lambda}^{-1}\mu^{\f12} f|_{L^{2}_{1/2}}^{2} \lesssim  C_{1} e^{\lambda} | f|_{L^{2}_{1/2}}^{2}.
\een
Recalling \eqref{inner-product-2-mu-mu} and \eqref{J-rho-lower-bound}, we get the upper bound in \eqref{lower-and-upper-bound}.

By relabeling the constants, we finish the proof.
\end{proof}

\subsection{Some preliminary formulas} \label{preliminary}
In this subsection,
we recall some useful formulas for the computation of Boltzmann type integrals involving $B(v-v_{*}, \sigma)=B(|v-v_{*}|, \f{v-v_{*}}{|v-v_{*}|}\cdot \sigma)$.
The change of variable $(v, v_{*}, \sigma) \rightarrow (v_{*}, v, -\sigma)$ has unit Jacobian and thus
\ben \label{change-v-and-v-star}
\int B(v-v_{*}, \sigma) f(v,v_{*},v^{\prime},v^{\prime}_{*}) \mathrm{d}V   =
\int B(v-v_{*}, \sigma) f(v_{*},v,v^{\prime}_{*},v^{\prime}) \mathrm{d}V  ,
\een
where $f$ is a general function such that the integral exists. Thanks to the symmetry of elastic collision formula \eqref{v-prime-v-prime-star}, the change of variable $(v, v_{*}, \sigma) \rightarrow (v^{\prime}, v^{\prime}_{*}, \frac{v-v_{*}}{|v-v_{*}|})$ has unit Jacobian and thus
\ben \label{change-v-and-v-prime}
\int B(v-v_{*}, \sigma) f(v,v_{*},v^{\prime},v^{\prime}_{*}) \mathrm{d}V   =
\int B(v-v_{*}, \sigma) f(v^{\prime},v^{\prime}_{*},v,v_{*}) \mathrm{d}V  .
\een

From now on, we will frequently use the notation \eqref{shorthand-D}.
By \eqref{shorthand-D} and the shorthand $f=f(v), f^{\prime} = f(v^{\prime}), f_{*}=f(v_{*}), f^{\prime}_{*} = f(v^{\prime}_{*})$, it is easy to see
\beno
\mathrm{D}(f) = - \mathrm{D}(f^{\prime}),
\quad
\mathrm{D}(f_{*}) = - \mathrm{D}(f^{\prime}_{*}), \quad
 \mathrm{D}^{2}(f) = \mathrm{D}^{2}(f^{\prime}), \quad \mathrm{D}^{2}(f_{*}) = \mathrm{D}^{2}(f^{\prime}_{*}).
\eeno
Similarly to \eqref{shorthand-D}, we introduce
\ben \label{shorthand-A}
\mathrm{A}(f(v,v_{*},v^{\prime},v^{\prime}_{*})) \colonequals  f(v,v_{*},v^{\prime},v^{\prime}_{*}) + f(v^{\prime},v^{\prime}_{*},v,v_{*}).
\een
The term $\mathrm{A}$ is interpreted as ``addition'' before and after collision.

 Thanks to the symmetry of elastic collision formula \eqref{v-prime-v-prime-star}, we have
\begin{lem}\label{v-square-sum-bounded} Let $v(\kappa) \colonequals     v + \kappa(v^{\prime}-v), v_{*}(\iota) \colonequals     v_{*} + \iota(v^{\prime}_{*}-v_{*})$ for $0 \leq \kappa, \iota \leq 1$. It holds that
\beno
\frac{1}{4}(|v|^{2}+|v_{*}|^{2}) \leq |v(\kappa)|^{2}+|v_{*}(\iota)|^{2} \leq 2(|v|^{2}+|v_{*}|^{2}).
\eeno
As a result, recalling $\mu (v) = e^{-\f{1}{2}|v|^{2}}$, for any $0 \leq \kappa, \iota \leq 1$, it holds that
\ben \label{mu-weight-result}
\mu^{2}(v) \mu^{2}(v_{*}) \leq \mu(v(\kappa)) \mu(v_{*}(\iota)) \leq \mu^{\frac{1}{4}}(v) \mu^{\frac{1}{4}}(v_{*}).
\een
\end{lem}
We omit the proof of Lemma \ref{v-square-sum-bounded} as it is elementary.
We will frequently use \eqref{mu-weight-result} to retain the good negative exponential ($\mu$-type) weight.

We now give a straightforward computation involving the regular change of variable $v \to v^{\prime}$.
\begin{lem} It holds that
\ben
\label{regular-change-of-variable-2}
\int \mathrm{1}_{\cos \theta \geq 0} |v-v_{*}| |g_{*}| f^{2} \mathrm{d}V +
\int \mathrm{1}_{\cos \theta \geq 0} |v-v_{*}| |g_{*}| (f^{2})^{\prime} \mathrm{d}V \lesssim |g|_{L^{1}_{1}}|f|_{L^{2}_{1/2}}^{2}.
\een
\end{lem}
\begin{proof}
For the former integral, as $|v-v_{*}| \leq (1+|v|)(1+|v_{*}|)$, the estimate is obvious.
For the latter integral, for fixed $\sigma, v_{*}$, using the change of variable $v \to v^{\prime}$, denoting
$\cos \alpha = \frac{v^{\prime}-v_{*}}{|v^{\prime}-v_{*}|} \cdot \sigma$,
we have $\theta =2 \alpha, |v-v_{*}| \cos\alpha = |v^{\prime}-v_{*}|,
0 \leq \alpha \leq \pi/4, \f{1}{\sqrt{2}} \leq \cos\alpha \leq 1, |\det(\f{\mathrm{d}v^{\prime}}{\mathrm{d}v})| = \f{\cos^{2}\alpha}{4}$.
Then
\beno
\int \mathrm{1}_{\cos \theta \geq 0}|v-v_{*}| |g_{*}| (f^{2})^{\prime} \mathrm{d}V =
\int \mathrm{1}_{\cos 2\alpha \geq 0} |v^{\prime}-v_{*}| |g_{*}|  (f^{2})^{\prime} \f{4}{\cos^{3}\alpha} \mathrm{d}\sigma \mathrm{d}v_{*} \mathrm{d}v^{\prime}
\lesssim |g|_{L^{1}_{1}} |f|_{L^{2}_{1/2}}^2,
\eeno
which ends the proof.
\end{proof}

Next, we then provide the computation involving the singular change of variable $v_{*} \to v^{\prime}$.
\begin{lem} For any fixed $v \in \mathbb{R}^{3}$, it holds that
\ben \label{singular-change-of-variable}
\int \mathrm{1}_{\cos \theta \geq 0} (\sin \f{\theta}{2})^{3/2} (f^{2})^{\prime} \mathrm{d}\sigma \mathrm{d}v_{*} =  \f{16 \pi}{2^{1/4}} |f|_{L^{2}}^{2}.
\een
\end{lem}
\begin{proof}
For fixed $\sigma, v$, using the change of variable $v_{*} \to v^{\prime}$, denoting
$\cos \alpha = \frac{v^{\prime}-v}{|v^{\prime}-v|} \cdot \sigma$,
we have
\beno
\alpha = \f{\pi}{2} - \f{\theta}{2}, \cos\alpha = \sin \f{\theta}{2},
|v-v_{*}| \cos\alpha = |v^{\prime}-v|,
\\
\pi/4 \leq \alpha \leq \pi/2,  0 \leq \cos\alpha \leq \f{1}{\sqrt{2}}, |\det(\f{\mathrm{d}v^{\prime}}{\mathrm{d}v_{*}})| = \f{\cos^{2}\alpha}{4}
\eeno
Then we get
\beno
\int \mathrm{1}_{\cos \theta \geq 0} 
(\sin \f{\theta}{2})^{3/2} (f^{2})^{\prime} \mathrm{d}\sigma \mathrm{d}v_{*}
= \int \mathrm{1}_{0 \leq \cos\alpha \leq \f{1}{\sqrt{2}}}  (f^{2})^{\prime} \f{4}{\cos^{1/2}\alpha} \mathrm{d}\sigma \mathrm{d}v^{\prime}
\eeno
Using $\frac{v^{\prime}-v}{|v^{\prime}-v|}$ as north pole to represent $\sigma$,
\beno
\int \mathrm{1}_{0 \leq \cos\alpha \leq \f{1}{\sqrt{2}}}  \f{4}{\cos^{1/2}\alpha} \mathrm{d}\sigma
= 2 \pi \int_{\pi/4}^{\pi/2} \f{4}{\cos^{1/2}\alpha} \sin \alpha \mathrm{d} \alpha = 8 \pi \int_{0}^{\f{1}{\sqrt{2}}} \f{1}{t^{1/2}} \mathrm{d}t = \f{16 \pi}{2^{1/4}}.
\eeno
Note that $(\sin \f{\theta}{2})^{3/2}$ is used to cancel the singularity of $\cos^{-2}\alpha$ near $\alpha = \pi/2$ in the change of  variable $v_{*} \to v^{\prime}$.
\end{proof}

\subsection{Non-linear operator estimate}
Recalling \eqref{definition-Gamma-2-epsilon-rho-T}, \eqref{definition-A-2-rho-T}, \eqref{line-1-rho-T}, \eqref{line-2-rho-T} and
\eqref{line-3-rho-T}, we have
\ben \label{tilde-Gamma-2-into-main-and-remaining}
\tilde{\Gamma}_{2}^{\lambda,T}(g,h) = \tilde{\Gamma}_{2,m}^{\lambda,T}(g,h) + \tilde{\Gamma}_{2,r}^{\lambda,T}(g,h), \quad
\tilde{\Gamma}_{2,r}^{\lambda,T}(g,h) \colonequals  \tilde{\Gamma}_{2,r,1}^{\lambda,T}(g,h) + \tilde{\Gamma}_{2,r,2}^{\lambda,T}(g,h) + \tilde{\Gamma}_{2,r,3}^{\lambda,T}(g,h).
\\ \label{tilde-Gamma-2-m}
\tilde{\Gamma}_{2,m}^{\lambda,T}(g,h) \colonequals  T^{2} N_{\lambda}^{-1} \int |v-v_{*}| \mathrm{D} \big((N_{\lambda}g)^{\prime}_{*}(N_{\lambda}h)^{\prime}\big)
\mathrm{d}\sigma \mathrm{d}v_{*}.
\\ \label{tilde-Gamma-2-r1}
\tilde{\Gamma}_{2,r,1}^{\lambda,T}(g,h) \colonequals  T^{2} N_{\lambda}^{-1} \int |v-v_{*}| \mathrm{D}\big((N_{\lambda}g)^{\prime}_{*}(N_{\lambda}h)^{\prime}(M_{\lambda} +  (M_{\lambda})_{*})\big)
\mathrm{d}\sigma \mathrm{d}v_{*}.
\\ \label{tilde-Gamma-2-r2}
\tilde{\Gamma}_{2,r,2}^{\lambda,T}(g,h) \colonequals  T^{2} N_{\lambda}^{-1} \int |v-v_{*}|
\mathrm{D}\big( (N_{\lambda}g)_{*}(N_{\lambda}h)^{\prime}( (M_{\lambda})^{\prime}_{*} -  M_{\lambda}) \big)
\mathrm{d}\sigma \mathrm{d}v_{*}.
\\ \label{tilde-Gamma-2-r3}
\tilde{\Gamma}_{2,r,3}^{\lambda,T}(g,h) \colonequals  T^{2} N_{\lambda}^{-1} \int |v-v_{*}|
\big( (N_{\lambda}g)^{\prime}(N_{\lambda}h)\mathrm{D}((M_{\lambda})^{\prime}_{*}) + (N_{\lambda}g)^{\prime}_{*}(N_{\lambda}h)_{*}\mathrm{D}(M_{\lambda}^{\prime})\big)
\mathrm{d}\sigma \mathrm{d}v_{*}.
\een

In the following Proposition, we derive upper bound estimates for the above bilinear operators.
\begin{prop} \label{bi-linear-term} For the bilinear terms \eqref{tilde-Gamma-2-m}, \eqref{tilde-Gamma-2-r1}, \eqref{tilde-Gamma-2-r2} and \eqref{tilde-Gamma-2-r3}, we have
\ben \label{bound-for-Gamma-2-m}
|\langle \tilde{\Gamma}_{2,m}^{\lambda,T}(g,h), f
\rangle| \lesssim e^{-\lambda/2} (1-e^{-\lambda})^{-2} T^{2} |g|_{L^{2}} |h|_{L^{2}_{1/2}} |f|_{L^{2}_{1/2}}.
\\ \label{bound-for-Gamma-2-r1}
|\langle \tilde{\Gamma}_{2,r,1}^{\lambda,T}(g,h), f
\rangle| \lesssim  e^{-3\lambda/2} (1-e^{-\lambda})^{-3} T^{2} |g|_{L^{2}} |h|_{L^{2}_{1/2}} |f|_{L^{2}_{1/2}}.
\\ \label{bound-for-Gamma-2-r2}
|\langle \tilde{\Gamma}_{2,r,2}^{\lambda,T}(g,h), f
\rangle| \lesssim  e^{-3\lambda/2} (1-e^{-\lambda})^{-3} T^{2} |\mu^{\f{1}{4}}g|_{L^{2}} |h|_{L^{2}_{1/2}} |f|_{L^{2}_{1/2}}.
\\ \label{bound-for-Gamma-2-r3}
|\langle \tilde{\Gamma}_{2,r,3}^{\lambda,T}(g,h), f
\rangle| \lesssim  e^{-3\lambda/2} (1-e^{-\lambda})^{-3}  T^{2} |\mu^{\f{1}{32}}g|_{L^{2}} |\mu^{\f{1}{64}}h|_{L^{2}}|\mu^{\f{1}{64}}f|_{L^{2}}.
\een
As a result, recalling \eqref{tilde-Gamma-2-into-main-and-remaining},
\ben \label{final-result-of-tilde-Gamma2}
|\langle \tilde{\Gamma}_{2}^{\lambda,T}(g,h), f
\rangle| \lesssim  C_{2,\lambda,T}  |g|_{L^{2}} |h|_{L^{2}_{1/2}} |f|_{L^{2}_{1/2}},
\een
where
\ben \label{defintion-of-C-2-lambda-T}
C_{2,\lambda,T} \colonequals  e^{-\lambda/2} (1-e^{-\lambda})^{-3} T^{2}.
\een
\end{prop}

\begin{proof} Note that
\beno \langle \tilde{\Gamma}_{2,m}^{\lambda,T}(g,h), f
\rangle = T^{2} \int
|v-v_{*}| \mathrm{D} \big((N_{\lambda}g)^{\prime}_{*}(N_{\lambda}h)^{\prime}\big) N_{\lambda}^{-1} f
\mathrm{d}V.
\eeno
Recalling \eqref{comparing-NN-with-NN-prime} and \eqref{M-rho-mu-N}, we have
\ben \label{typical-bound-for-a-term}
N_{\lambda}^{-1} |\mathrm{D} \big((N_{\lambda}g)^{\prime}_{*}(N_{\lambda}h)^{\prime}\big)| \lesssim (1-e^{-\lambda})^{-1}  (N_{\lambda})_{*} \mathrm{A} (|g_{*}h|) \lesssim e^{-\lambda/2} (1-e^{-\lambda})^{-2} \mu^{\f{1}{2}}_{*} \mathrm{A} (|g_{*}h|).
\een
Then
\beno |\langle \tilde{\Gamma}_{2,m}^{\lambda,T}(g,h), f
\rangle| \lesssim  T^{2} e^{-\lambda/2} (1-e^{-\lambda})^{-2}
 \int
|v-v_{*}| \mu^{\f{1}{2}}_{*} \mathrm{A} (|g_{*}h|) |f| \mathrm{d}V.
\eeno
Note that the integral equals to
\beno
 \int
|v-v_{*}| \mu^{\f{1}{2}}_{*} |g_{*} h f|
\mathrm{d}V +
 \int
|v-v_{*}| (\mu^{\f{1}{2}})^{\prime}_{*} |g_{*} h f^{\prime}|
\mathrm{d}V.
\eeno
As $\int \mathrm{1}_{\cos \theta \geq 0}\mathrm{d} \sigma = 2 \pi$, it is easy to derive
\beno
\int
|v-v_{*}| \mu^{\f{1}{2}}_{*}|g_{*} h f|
\mathrm{d}V \lesssim |\mu^{\f{1}{4}}g|_{L^{2}} |h|_{L^{2}_{1/2}} |f|_{L^{2}_{1/2}}.
\eeno
Since $|v-v_{*}| = |v^{\prime}-v^{\prime}_{*}| \leq \sqrt{2}|v-v^{\prime}_{*}|$,
 we get
\beno
\int
|v-v_{*}| (\mu^{\f{1}{2}})^{\prime}_{*} |g_{*} h f^{\prime}|
\mathrm{d}V &\leq&
\cs
{\int
|v-v_{*}| (\mu^{\f{1}{2}})^{\prime}_{*} g_{*}^{2} h^{2}
\mathrm{d}V}
{\int
|v-v_{*}| (\mu^{\f{1}{2}})^{\prime}_{*} (f^{2})^{\prime}
\mathrm{d}V},
\\ &\lesssim& \cs
{\int
\langle v \rangle g_{*}^{2} h^{2}
\mathrm{d}V}
{\int
|v-v_{*}| \mu^{\f{1}{2}}_{*} f^{2}
\mathrm{d}V},
\lesssim
|g|_{L^{2}} |h|_{L^{2}_{1/2}} |f|_{L^{2}_{1/2}}.
\eeno
Patching the above formulas together, we get \eqref{bound-for-Gamma-2-m}.

Note that
\beno \langle \tilde{\Gamma}_{2,r,1}^{\lambda,T}(g,h), f
\rangle = T^{2} \int
|v-v_{*}| \mathrm{D}\big((N_{\lambda}g)^{\prime}_{*}(N_{\lambda}h)^{\prime}(M_{\lambda} +  (M_{\lambda})_{*})\big)
 N_{\lambda}^{-1} f \mathrm{d}V.
\eeno
By \eqref{M-rho-mu-N}, $M_{\lambda} \lesssim e^{-\lambda} (1-e^{-\lambda})^{-1}$, from which together with the above argument,
 we obtain \eqref{bound-for-Gamma-2-r1}.
By exactly the same derivation, we also get \eqref{bound-for-Gamma-2-r2}.

Note that
\beno \langle \tilde{\Gamma}_{2,r,3}^{\lambda,T}(g,h), f
\rangle = T^{2} \int
|v-v_{*}| \big( (N_{\lambda}g)^{\prime}(N_{\lambda}h) \mathrm{D}((M_{\lambda})^{\prime}_{*}) + (N_{\lambda}g)^{\prime}_{*}(N_{\lambda}h)_{*}\mathrm{D}(M_{\lambda}^{\prime})\big)
 N_{\lambda}^{-1} f \mathrm{d}V.
\eeno
Recalling \eqref{M-rho-mu-N} and using \eqref{mu-weight-result},  we get
\beno
N_{\lambda}^{-1}|(N_{\lambda})^{\prime}(N_{\lambda}) \mathrm{D}((M_{\lambda})^{\prime}_{*})| \lesssim
e^{-3\lambda/2} (1-e^{-\lambda})^{-2} (\mu^{\f{1}{2}})^{\prime} (\mu_{*} + \mu^{\prime}_{*})
\lesssim
e^{-3\lambda/2} (1-e^{-\lambda})^{-2} \mu^{\f{1}{8}} \mu^{\f{1}{8}}_{*}.
\eeno
Recalling \eqref{comparing-NN-with-NN-prime} and \eqref{M-rho-mu-N}, using \eqref{mu-weight-result},  we can similarly get
\beno
N_{\lambda}^{-1}|(N_{\lambda})^{\prime}_{*}(N_{\lambda})_{*}\mathrm{D}(M_{\lambda}^{\prime})| \lesssim
e^{-3\lambda/2} (1-e^{-\lambda})^{-3} \mu^{\f{1}{8}} \mu^{\f{1}{8}}_{*}.
\eeno
Therefore
\beno
|\langle \tilde{\Gamma}_{2,r,3}^{\lambda,T}(g,h), f
\rangle| &\lesssim& e^{-3\lambda/2} (1-e^{-\lambda})^{-3} T^{2}
 \int
|v-v_{*}| \mu^{\f{1}{8}} \mu^{\f{1}{8}}_{*} \big( |g^{\prime} h f| +
|g^{\prime}_{*} h_{*} f| \big)
   \mathrm{d}V
   \\&\lesssim& e^{-3\lambda/2} (1-e^{-\lambda})^{-3} T^{2} \int
\mu^{\f{1}{16}} \mu^{\f{1}{16}}_{*} \big( |g^{\prime} h f| +
|g^{\prime}_{*} h_{*} f| \big)
   \mathrm{d}V.
\eeno
Note that we retain the good weight. Let us see $\int
\mu^{\f{1}{16}} \mu^{\f{1}{16}}_{*} |g^{\prime} h f|
   \mathrm{d}V$. By the C-S inequality for the integral w.r.t. $\mathrm{d}\sigma \mathrm{d}v_{*}$, by \eqref{singular-change-of-variable},
   we have
\beno
\int
\mu^{\f{1}{16}} \mu^{\f{1}{16}}_{*} |g^{\prime} h f|
   \mathrm{d}V &\leq& \int
\mu^{\f{1}{32}} |h f| \left(\cs{\int \mathrm{1}_{\cos \theta \geq 0} (\sin \f{\theta}{2})^{3/2} (\mu^{\f{1}{16}}g^{2})^{\prime} \mathrm{d}\sigma \mathrm{d}v_{*}}
{\int \mathrm{1}_{\cos \theta \geq 0} (\sin \f{\theta}{2})^{-3/2} \mu^{\f{1}{16}}_{*} \mathrm{d}\sigma \mathrm{d}v_{*}}\right)
   \mathrm{d}v
\\  &\lesssim&  |\mu^{\f{1}{32}}g|_{L^{2}} |\mu^{\f{1}{64}}h|_{L^{2}}|\mu^{\f{1}{64}}f|_{L^{2}}.
\eeno
Here we use
\beno
\int \mathrm{1}_{\cos \theta \geq 0} (\sin \f{\theta}{2})^{-3/2}  \mathrm{d}\sigma = 8 \pi \int_{0}^{\pi/2}  \sin^{-\f{1}{2}} \f{\theta}{2} \mathrm{d} \sin\f{\theta}{2} = 8 \pi \int_{0}^{\f{1}{\sqrt{2}}} \f{1}{t^{1/2}} \mathrm{d}t = \f{16 \pi}{2^{1/4}}.
\eeno
By using \eqref{regular-change-of-variable-2},
it is easy to see $\int
\mu^{\f{1}{16}} \mu^{\f{1}{16}}_{*} |g^{\prime}_{*} h_{*} f|
   \mathrm{d}V = \int
\mu^{\f{1}{16}} \mu^{\f{1}{16}}_{*} |g^{\prime} h f_{*}|
   \mathrm{d}V \lesssim |\mu^{\f{1}{32}}g|_{L^{2}} |\mu^{\f{1}{64}}h|_{L^{2}}|\mu^{\f{1}{64}}f|_{L^{2}}$. Patching together the two estimates, we get \eqref{bound-for-Gamma-2-r3}.
\end{proof}

In the following Proposition, we derive upper bound estimate for the trilinear operator $\tilde{\Gamma}_{3}^{\lambda,T}(\cdot , \cdot, \cdot)$.
\begin{prop} \label{tri-linear-term}
\beno |\langle \tilde{\Gamma}_{3}^{\lambda,T}(g,h,\varrho), f
\rangle| \lesssim  C_{3,\lambda,T} |g|_{L^{2}} |h|_{L^{2}_{1/2}} |\mu^{\f{1}{64}}\varrho|_{L^{2}} |f|_{L^{2}_{1/2}},
\eeno
where
\ben \label{defintion-of-C-3-lambda-T}
C_{3,\lambda,T} \colonequals  e^{-\lambda} (1-e^{-\lambda})^{-3} T^{2}.
\een
\end{prop}

\begin{proof} Recalling \eqref{definition-Gamma-3-epsilon-rho-T}, we have
\beno \langle \tilde{\Gamma}_{3}^{\lambda,T}(g,h,\varrho), f
\rangle =   T^{2} \int
|v-v_{*}| \mathrm{D} \big( (N_{\lambda}g)_{*}^{\prime} (N_{\lambda}h)^{\prime} ((N_{\lambda}\varrho)_{*} + N_{\lambda}\varrho) \big) N_{\lambda}^{-1} f
\mathrm{d}V
.
\eeno
Similarly to \eqref{typical-bound-for-a-term}, we get
\beno
N_{\lambda}^{-1} |\mathrm{D} \big( (N_{\lambda}g)_{*}^{\prime} (N_{\lambda}h)^{\prime} ((N_{\lambda}\varrho)_{*} + N_{\lambda}\varrho) \big)| \lesssim e^{-\lambda} (1-e^{-\lambda})^{-3} \mu^{\f{1}{2}}_{*} \mathrm{A} (|g_{*}h\varrho^{\prime}_{*}|(\mu^{\f{1}{2}})^{\prime}_{*} + |g_{*}h\varrho^{\prime}|(\mu^{\f{1}{2}})^{\prime}).
\eeno
From this together with $\mu \mu_{*} = \mu^{\prime} \mu^{\prime}_{*}$,
we have
\beno |\langle \tilde{\Gamma}_{3}^{\lambda,T}(g,h,\varrho), f
\rangle| \lesssim   T^{2} e^{-\lambda} (1-e^{-\lambda})^{-3} (\mathcal{I}_{1}+\mathcal{I}_{2}+\mathcal{I}_{3}+\mathcal{I}_{4}),
\eeno
where
\beno
\mathcal{I}_{1} \colonequals  \int
|v-v_{*}| \mu^{\f{1}{2}}_{*}  (\mu^{\f{1}{2}})^{\prime}_{*} |g_{*} h \varrho^{\prime}_{*} f|
\mathrm{d}V, \quad
\mathcal{I}_{2} \colonequals  \int
|v-v_{*}| \mu^{\prime}_{*} |g_{*} h \varrho^{\prime}_{*} f^{\prime}|
\mathrm{d}V,
\\
\mathcal{I}_{3} \colonequals  \int
|v-v_{*}| \mu^{\f{1}{2}}_{*} (\mu^{\f{1}{2}})^{\prime} |g_{*} h \varrho^{\prime} f|
\mathrm{d}V, \quad
\mathcal{I}_{4} \colonequals  \int
|v-v_{*}| \mu^{\f{1}{2}}_{*}  \mu^{\f{1}{2}} |g_{*} h \varrho^{\prime} f^{\prime}|
\mathrm{d}V.
\eeno
By the C-S inequality, \eqref{change-v-and-v-prime}
and \eqref{regular-change-of-variable-2},
we have
\beno
|\mathcal{I}_{1}| \leq \cs{\int
|v-v_{*}| \mu_{*}  |g_{*} h|^{2}
\mathrm{d}V}{\int
|v-v_{*}| \mu_{*} |\varrho_{*} f^{\prime}|^{2}
\mathrm{d}V} \lesssim |\mu^{\f{1}{4}}g|_{L^{2}} |h|_{L^{2}_{1/2}} |\mu^{\f{1}{4}}\varrho|_{L^{2}} |f|_{L^{2}_{1/2}}.
\eeno
By the C-S inequality, \eqref{change-v-and-v-prime}, the fact $|v-v_{*}| \mu^{\prime}_{*} \lesssim \langle v \rangle$,
and the estimate \eqref{regular-change-of-variable-2},
we have
\beno
|\mathcal{I}_{2}| \leq \cs{\int
|v-v_{*}| \mu^{\prime}_{*}  |g_{*} h|^{2}
\mathrm{d}V}{\int
|v-v_{*}| \mu_{*} |\varrho_{*} f|^{2}
\mathrm{d}V} \lesssim |g|_{L^{2}} |h|_{L^{2}_{1/2}} |\mu^{\f{1}{4}}\varrho|_{L^{2}} |f|_{L^{2}_{1/2}}.
\eeno
By using \eqref{mu-weight-result}, we have $\mu^{\f{1}{2}}_{*} (\mu^{\f{1}{2}})^{\prime} \leq \mu^{\f{1}{8}}_{*} \mu^{\f{1}{8}}$ and thus
\beno
|\mathcal{I}_{3}| \leq \int
|v-v_{*}| \mu^{\f{1}{8}}_{*} \mu^{\f{1}{8}} |g_{*} h \varrho^{\prime} f|
\mathrm{d}V \lesssim \int
 \mu^{\f{1}{16}}_{*} \mu^{\f{1}{16}} |g_{*} h \varrho^{\prime} f|
\mathrm{d}V
\eeno
Using $\mu \mu_{*} = \mu^{\prime} \mu^{\prime}_{*}$, we get $ \mu^{\f{1}{16}}_{*} \mu^{\f{1}{16}} =  \mu^{\f{1}{32}}_{*} \mu^{\f{1}{32}}  (\mu^{\f{1}{32}})_{*}^{\prime} (\mu^{\f{1}{32}})^{\prime}$ and so
\beno
|\mathcal{I}_{3}| \lesssim  \int
  |(\mu^{\f{1}{64}}g)_{*} (\mu^{\f{1}{64}}h) (\mu^{\f{1}{64}}\varrho)^{\prime} (\mu^{\f{1}{64}}f)|
\mathrm{d}V
\eeno
For any fixed $v$,  we bound the integral over  $\mathrm{d}\sigma \mathrm{d}v_{*}$ by the C-S inequality and \eqref{singular-change-of-variable} to get \beno
\int
  |(\mu^{\f{1}{64}}g)_{*} (\mu^{\f{1}{64}}\varrho)^{\prime}| \mathrm{d}\sigma \mathrm{d}v_{*}
 \leq \cs{\int
   (\mu^{\f{1}{32}}g^{2})_{*} 
   (\sin \f{\theta}{2})^{-3/2} \mathrm{d}\sigma \mathrm{d}v_{*}
 }{\int
 (\mu^{\f{1}{32}}\varrho^{2})^{\prime} (\sin \f{\theta}{2})^{3/2} \mathrm{d}\sigma \mathrm{d}v_{*}
 } \lesssim |\mu^{\f{1}{64}}g|_{L^{2}} |\mu^{\f{1}{64}}\varrho|_{L^{2}}.
\eeno
Then
\beno
|\mathcal{I}_{3}| \lesssim |\mu^{\f{1}{64}}g|_{L^{2}} |\mu^{\f{1}{64}}h|_{L^{2}} |\mu^{\f{1}{64}}\varrho|_{L^{2}} |\mu^{\f{1}{64}}f|_{L^{2}}.
\eeno
Note that
\beno
\mathcal{I}_{4} = \int
|v-v_{*}| \mu^{\f{1}{2}}_{*}  \mu^{\f{1}{2}} |g^{\prime}_{*} h^{\prime} \varrho f|
\mathrm{d}V \lesssim \int
  |(\mu^{\f{1}{16}}g)^{\prime}_{*} (\mu^{\f{1}{16}}h)^{\prime} (\mu^{\f{1}{16}}\varrho) (\mu^{\f{1}{16}}f)|
\mathrm{d}V
\eeno
For any fixed $v$,  by the C-S inequality, using the change of variable $v_{*} \to v^{\prime}_{*}$(in which the Jacobian $\f{4}{\cos^{2}(\theta/2)}$ is bounded) and the estimate \eqref{singular-change-of-variable}, we have
\beno
\int
  |(\mu^{\f{1}{16}}g)^{\prime}_{*} (\mu^{\f{1}{16}}h)^{\prime}| \mathrm{d}\sigma \mathrm{d}v_{*}
 \leq \cs{\int
   (\mu^{\f{1}{8}}g^{2})^{\prime}_{*} (\sin \f{\theta}{2})^{-3/2} \mathrm{d}\sigma \mathrm{d}v_{*}
 }{\int
 (\mu^{\f{1}{8}}h^{2})^{\prime} (\sin \f{\theta}{2})^{3/2} \mathrm{d}\sigma \mathrm{d}v_{*}
 } \lesssim |\mu^{\f{1}{16}}g|_{L^{2}} |\mu^{\f{1}{16}}h|_{L^{2}},
\eeno
which gives
\beno
|\mathcal{I}_{4}| \lesssim |\mu^{\f{1}{16}}g|_{L^{2}} |\mu^{\f{1}{16}}h|_{L^{2}} |\mu^{\f{1}{16}}\varrho|_{L^{2}} |\mu^{\f{1}{16}}f|_{L^{2}}.
\eeno
Patching together the above estimates, we finish the proof.
\end{proof}

\subsection{Energy estimates}

For simplicity, from now on $\mathcal{A} = (a,b,c)$.
By \eqref{L2-of-N} and \eqref{weighted-L2-k-geq-1-of-N},
we have
\ben  \label{f1-weight-equivalent}
|\partial^{\alpha}f_{1}|_{L^{2}_{k}} \lesssim |\partial^{\alpha}f_{1}|_{L^{2}} \leq |\partial^{\alpha}f|_{L^{2}}.
\een
By \eqref{f1-weight-equivalent} and
 recalling \eqref{D-N-lower-bound}, we will frequently use
\ben \label{f1-any-order-bounded-by-energy}
0 \leq j \leq N, \quad 0 \leq |\alpha| \leq N-j \quad \Rightarrow \quad \|\partial^{\alpha}f_1\|_{H^{j}_{x}L^{2}_{1/2}} \lesssim \|f\|_{H^{N}_{x}L^{2}}.
\\
\label{any-order-bounded-by-energy}
0 \leq j \leq N, \quad 0 \leq |\alpha| \leq N-j \quad \Rightarrow \quad \|\partial^{\alpha}f\|_{H^{j}_{x}L^{2}} \leq \|f\|_{H^{N}_{x}L^{2}}.
\\ \label{more-than-1-order-bounded-by-dissipation}
0 \leq k \leq N-1, \quad 1 \leq |\alpha| \leq N-k \quad \Rightarrow \quad \|\partial ^{\alpha}f\|_{H^{k}_{x}L^{2}_{1/2}}^{2} \lesssim \mathcal{D}_{N}(f).
\een

Note that the dissipation $\mathcal{D}_{N}(f)$ lacks $\|f_{1}\|_{L^{2}_{x}L^{2}} \lesssim e^{-\lambda/2} (1-e^{-\lambda})^{-1/4}|\mathcal{A}|_{L^{2}_{x}}$. We will use the following
embedding in dimension 3,
\ben  \label{L6-bounded-by-H-dot-1}
|f|_{L^{6}} \lesssim |\nabla f|_{L^{2}}.
\een
In the 3-dimensional space $\mathbb{R}^{3}_{x}$,
by embedding $L^{\infty}_{x} \hookrightarrow H^{2}_{x}$ or $L^{p}_{x} \hookrightarrow H^{s}_{x}$ with $\frac{s}{3} = \f{1}{2} - \frac{1}{p}$. From these basic embedding results, by \eqref{final-result-of-tilde-Gamma2},
estimates of the bi-linear term $\tilde{\Gamma}_{2}^{\lambda,T}$  in the full space $(x,v)$ are given as
\ben \label{Gamma-2-ub-weighted-full-space}
|(\tilde{\Gamma}_{2}^{\lambda,T}(g,h), f)|
\lesssim C_{2,\lambda,T} \|g\|_{H^{m}_{x}L^{2}}  \|h\|_{H^{n}_{x}L^{2}_{1/2}} \|f\|_{L^{2}_{x}L^{2}_{1/2}},
\\ \label{Gamma-2-ub-weighted-dot-1}
|(\tilde{\Gamma}_{2}^{\lambda,T}(g,h), f)|
\lesssim C_{2,\lambda,T} \|g\|_{H^{1/2}_{x}L^{2}}  \|\nabla_{x}h\|_{L^{2}_{x}L^{2}_{1/2}} \|f\|_{L^{2}_{x}L^{2}_{1/2}},
\een
where $m,n \in \mathbb{N}, m+n=2$.

Based on \eqref{Gamma-2-ub-weighted-full-space} and \eqref{Gamma-2-ub-weighted-dot-1}, by making a suitable choice of parameters $m,n$ to deal with different distribution of derivative order, we get the following estimate.

\begin{thm} \label{Gamma-2-energy-estimate-label}  Let $N \geq 2$. It holds that
\ben  \label{Gamma-2-energy-estimate}
|\sum_{|\alpha|\leq N} ( \partial^{\alpha} \tilde{\Gamma}_{2}^{\lambda,T}(f,f),  \partial^{\alpha}f )| \lesssim C_{2,\lambda,T}\left( \|f\|_{H^{2}_{x}L^{2}} \mathcal{D}^{\f{1}{2}}_{N}(f) + \mathrm{1}_{N \geq 3} C_{N}
\|f\|_{H^{N}_{x}L^{2}}\mathcal{D}^{\f{1}{2}}_{N-1}(f) \right) \|f_{2}\|_{H^{N}_{x}L^{2}_{1/2}}.
\een
\end{thm}
\begin{proof} By binomial formula, we need to consider $(\tilde{\Gamma}_{2}^{\lambda,T}(\partial^{\alpha_{1}}f,\partial^{\alpha_{2}}f), \partial^{\alpha}f)$ for all combinations of $\alpha_{1}+\alpha_{2}=\alpha$ with $|\alpha|\leq N$. We first derive  $\tilde{\Gamma}_{2}^{\lambda,T}(g,h) + \tilde{\Gamma}_{2}^{\lambda,T}(h,g) \in (\tilde{\mathcal{L}}^{\lambda,T})^{\perp} $. Recalling \eqref{definition-Gamma-2-epsilon-rho-T}, we have
\beno
\langle \tilde{\Gamma}_{2}^{\lambda,T}(g,h) + \tilde{\Gamma}_{2}^{\lambda,T}(h,g), f \rangle =  \int
B ( \tilde{\Pi}_{2}(g,h) + \tilde{\Pi}_{2}(h,g) ) N^{-1} f
\mathrm{d}V.
\eeno
Recall \eqref{definition-A-2-rho-T} for the definition of $\tilde{\Pi}_{2}$.
By \eqref{change-v-and-v-star} and \eqref{change-v-and-v-prime}, we have
\beno
\langle \tilde{\Gamma}_{2}^{\lambda,T}(g,h) + \tilde{\Gamma}_{2}^{\lambda,T}(h,g), f \rangle =  \f{1}{4}\int
B ( \tilde{\Pi}_{2}(g,h) + \tilde{\Pi}_{2}(h,g) ) \mathcal{S}(N^{-1} f)
\mathrm{d}V.
\eeno
Therefore,  for any $f \in \tilde{\mathcal{L}}^{\lambda,T}$, we have
\ben \label{symmetry-vanishes-in-kernel}
\langle \tilde{\Gamma}_{2}^{\lambda,T}(g,h) + \tilde{\Gamma}_{2}^{\lambda,T}(h,g), f \rangle = 0.
\een
Note that
\beno
\partial^{\alpha}\tilde{\Gamma}_{2}^{\lambda,T}(f, f) = \sum_{\alpha_{1}+\alpha_{2}=\alpha} C_{\alpha}^{\alpha_{1}}
\tilde{\Gamma}_{2}^{\lambda,T}(\partial^{\alpha_{1}}f, \partial^{\alpha_{2}}f) = \f{1}{2}\sum_{\alpha_{1}+\alpha_{2}=\alpha} C_{\alpha}^{\alpha_{1}}
(\tilde{\Gamma}_{2}^{\lambda,T}(\partial^{\alpha_{1}}f, \partial^{\alpha_{2}}f)+\tilde{\Gamma}_{2}^{\lambda,T}(\partial^{\alpha_{2}}f, \partial^{\alpha_{1}}f).
\eeno
From this together with \eqref{symmetry-vanishes-in-kernel}, for any $\varphi \in \tilde{\mathcal{L}}^{\lambda,T}$, we get
\ben \label{taking-derivative-vanish-with-kernel}
\langle \partial^{\alpha}\tilde{\Gamma}_{2}^{\lambda,T}(f, f), \varphi \rangle = 0.
\een

By \eqref{taking-derivative-vanish-with-kernel}, we have
\beno
(\partial^{\alpha}\tilde{\Gamma}_{2}^{\lambda,T}(f, f), \partial^{\alpha}f) =
(\partial^{\alpha}\tilde{\Gamma}_{2}^{\lambda,T}(f, f), \partial^{\alpha}f_{2}).
\eeno
We will prove that for all combinations of $\alpha_{1}+\alpha_{2}=\alpha$ with $0 \leq |\alpha| \leq N$,
the following inequality holds
\beno
|(\tilde{\Gamma}_{2}^{\lambda,T}(\partial^{\alpha_1}f, \partial^{\alpha_2}f), \partial^{\alpha}f_{2})|
\lesssim  C_{2,\lambda,T}\left( \|f\|_{H^{2}_{x}L^{2}} \mathcal{D}^{\f{1}{2}}_{N}(f) + \mathrm{1}_{N \geq 3}
\|f\|_{H^{N}_{x}L^{2}}\mathcal{D}^{\f{1}{2}}_{N-1}(f) \right) \|f_{2}\|_{H^{N}_{x}L^{2}_{1/2}}.
\eeno

We first deal with the case $|\alpha| = 0$.
By \eqref{Gamma-2-ub-weighted-dot-1} and \eqref{more-than-1-order-bounded-by-dissipation},  we get
\ben \label{alpha=0}
|(\tilde{\Gamma}_{2}^{\lambda,T}(f, f), f_{2})| \lesssim C_{2,\lambda,T}
\|f\|_{H^{1}_{x}L^{2}}  \|\nabla_{x}f\|_{L^{2}_{x}L^{2}_{1/2}} \|f_2\|_{L^{2}_{x}L^{2}_{1/2}}
\lesssim C_{2,\lambda,T} \|f\|_{H^{1}_{x}L^{2}} \mathcal{D}^{\f{1}{2}}_{1}(f)\|f_{2}\|_{L^{2}_{x}L^{2}_{1/2}}.
\een

Now it remains to consider $1 \leq |\alpha| \leq N$.
By \eqref{Gamma-2-ub-weighted-full-space}, it suffices to prove that for all combinations of $\alpha_{1}+\alpha_{2}=\alpha$ with $1 \leq |\alpha| \leq N$,
the following inequality
\beno
\|\partial^{\alpha_{1}}f\|_{H^{m}_{x}L^{2}}  \|\partial^{\alpha_{2}}f\|_{H^{n}_{x}L^{2}_{1/2}} \|\partial^{\alpha}f_{2}\|_{L^{2}_{x}L^{2}_{1/2}}
\lesssim \|f\|_{H^{2}_{x}L^{2}} \mathcal{D}^{\f{1}{2}}_{N}(f) \|f_{2}\|_{H^{N}_{x}L^{2}_{1/2}} + \mathrm{1}_{N \geq 3}
\|f\|_{H^{N}_{x}L^{2}}\mathcal{D}^{\f{1}{2}}_{N-1}(f)\|f_{2}\|_{H^{N}_{x}L^{2}_{1/2}}.
\eeno
holds  for some $m,n$ verifying $m,n \in \mathbb{N}, m+n=2$.

 If $\alpha_{2}=\alpha$,
Taking $(m,n) = (2,0)$ and using \eqref{more-than-1-order-bounded-by-dissipation}, we have
\beno
\|f\|_{H^{2}_{x}L^{2}}  \|\partial^{\alpha}f\|_{H^{0}_{x}L^{2}_{1/2}} \|\partial^{\alpha}f_{2}\|_{L^{2}_{x}L^{2}_{1/2}}
\lesssim  \|f\|_{H^{2}_{x}L^{2}} \mathcal{D}^{\f{1}{2}}_{N}(f) \|f_{2}\|_{H^{N}_{x}L^{2}_{1/2}}
\eeno
Therefore in the following we only need to consider all combinations of $\alpha_{1}+\alpha_{2}=\alpha$ with $1 \leq |\alpha| \leq N, |\alpha_{2}| \leq |\alpha|-1$.

The following is divided into three cases: $|\alpha| = 1; |\alpha| = 2; 3 \leq |\alpha| \leq  N$.

{\it {Case 1: $|\alpha| = 1$.}}
In this case, there remains only one case: $\alpha_{2} = 0, \alpha_{1} = \alpha$. Taking $(m,n) = (0,2)$, using \eqref{f1-any-order-bounded-by-energy}, \eqref{any-order-bounded-by-energy} and \eqref{more-than-1-order-bounded-by-dissipation}, we have
\ben \label{alpha=1}
&&\|\partial^{\alpha}f\|_{H^{0}_{x}L^{2}}  \|f\|_{H^{2}_{x}L^{2}_{1/2}} \|\partial^{\alpha}f_{2}\|_{L^{2}_{x}L^{2}_{1/2}}
\\ \nonumber &\lesssim& \|\partial^{\alpha}f\|_{H^{0}_{x}L^{2}}  (\|f_{1}\|_{H^{2}_{x}L^{2}_{1/2}} + \|f_{2}\|_{H^{2}_{x}L^{2}_{1/2}}) \|\partial^{\alpha}f_{2}\|_{L^{2}_{x}L^{2}_{1/2}}
\\ \nonumber &\lesssim& \mathcal{D}^{\f{1}{2}}_{1}(f) \|f\|_{H^{2}_{x}L^{2}} \|f_{2}\|_{H^{1}_{x}L^{2}_{1/2}} +
\|f\|_{H^{1}_{x}L^{2}}  \mathcal{D}^{\f{1}{2}}_{2}(f) \|f_{2}\|_{H^{1}_{x}L^{2}_{1/2}}
\lesssim \|f\|_{H^{2}_{x}L^{2}} \mathcal{D}^{\f{1}{2}}_{2}(f) \|f_{2}\|_{H^{1}_{x}L^{2}_{1/2}}
\een

{\it {Case 2: $|\alpha| = 2$.}}
In this case, it remains to consider two subcases: $|\alpha_{2}| = 0; |\alpha_{2}| = 1$. In the first subcase $|\alpha_{2}| = 0, \alpha_{1} =\alpha$. Taking $(m,n) = (0,2)$, similarly to \eqref{alpha=1},
we get
\ben \label{first-full-second-0-order-2}
\|\partial^{\alpha}f\|_{H^{0}_{x}L^{2}}  \|f\|_{H^{2}_{x}L^{2}_{1/2}} \|\partial^{\alpha}f_{2}\|_{L^{2}_{x}L^{2}_{1/2}}
\lesssim \|f\|_{H^{2}_{x}L^{2}} \mathcal{D}^{\f{1}{2}}_{2}(f) \|f_{2}\|_{H^{2}_{x}L^{2}_{1/2}}
\een
In the second subcase $|\alpha_{2}| = |\alpha_{1}| = 1$. Taking $(m,n) = (1,1)$ and using \eqref{more-than-1-order-bounded-by-dissipation}, we have
\beno
\|\partial^{\alpha_1}f\|_{H^{1}_{x}L^{2}}  \|\partial^{\alpha_2}f\|_{H^{1}_{x}L^{2}_{1/2}} \|\partial^{\alpha}f_{2}\|_{L^{2}_{x}L^{2}_{1/2}} \lesssim   \|f\|_{H^{2}_{x}L^{2}}\mathcal{D}^{\f{1}{2}}_{2}(f) \|f_{2}\|_{H^{2}_{x}L^{2}_{1/2}}.
\eeno

{\it {Case 3: $3 \leq |\alpha| \leq N$.}}
In this case, it remains to consider three subcases: $|\alpha_{2}| = 0;  |\alpha_{2}| = 1; 2 \leq |\alpha_{2}| \leq |\alpha| -1$. In the first subcase $|\alpha_{2}| = 0, \alpha_{1} =\alpha$. Taking $(m,n) = (0,2)$, similarly to \eqref{alpha=1},
we get
\ben \label{first-full-second-0-order-N}
\|\partial^{\alpha}f\|_{H^{0}_{x}L^{2}}  \|f\|_{H^{2}_{x}L^{2}_{1/2}} \|\partial^{\alpha}f_{2}\|_{L^{2}_{x}L^{2}_{1/2}}
\lesssim \|f\|_{H^{2}_{x}L^{2}} \mathcal{D}^{\f{1}{2}}_{N}(f) \|f_{2}\|_{H^{N}_{x}L^{2}_{1/2}} +
\|f\|_{H^{N}_{x}L^{2}}  \mathcal{D}^{\f{1}{2}}_{2}(f) \|f_{2}\|_{H^{N}_{x}L^{2}_{1/2}}.
\een
In the second subcase $|\alpha_{2}| = 1, |\alpha_{1}| = N-1$. Taking $(m,n) = (1,1)$ and using \eqref{more-than-1-order-bounded-by-dissipation}, we have
\beno
\|\partial^{\alpha_{1}}f\|_{H^{1}_{x}L^{2}}  \|\partial^{\alpha_{2}}f\|_{H^{1}_{x}L^{2}_{1/2}} \|\partial^{\alpha}f_{2}\|_{L^{2}_{x}L^{2}_{1/2}}
\lesssim \|f\|_{H^{N}_{x}L^{2}}  \mathcal{D}^{\f{1}{2}}_{2}(f) \|f_{2}\|_{H^{N}_{x}L^{2}_{1/2}}.
\eeno
In the third subcase $2 \leq |\alpha_{2}| \leq |\alpha| -1 \leq N-1, |\alpha_{1}| \leq N-2$. Taking $(m,n) = (2,0)$ and using \eqref{more-than-1-order-bounded-by-dissipation}, we have
\beno
\|\partial^{\alpha_{1}}f\|_{H^{2}_{x}L^{2}}  \|\partial^{\alpha_{2}}f\|_{H^{0}_{x}L^{2}_{1/2}} \|\partial^{\alpha}f_{2}\|_{L^{2}_{x}L^{2}_{1/2}}
\lesssim \|f\|_{H^{N}_{x}L^{2}}  \mathcal{D}^{\f{1}{2}}_{N-1}(f) \|f_{2}\|_{H^{N}_{x}L^{2}_{1/2}}.
\eeno

Patching together all the above estimates, observing that {\it {Case 3}} happens only if $N \geq 3$,
we obtain  \eqref{Gamma-2-energy-estimate}.
\end{proof}

From Proposition \ref{tri-linear-term}, by H\"{o}lder's inequality and Sobolev embedding inequalities, 
estimates of the trilinear term $\tilde{\Gamma}_{3}^{\lambda,T}$  in the full space $(x,v)$ are given as
\ben \label{Gamma-3-ub-full-space}
|(\tilde{\Gamma}_{3}^{\lambda,T}(g,h,\varrho), f)|
\lesssim  C_{3,\lambda,T} \|g\|_{H^{r_{1}}_{x}L^{2}}  \|h\|_{H^{r_{2}}_{x}L^{2}_{1/2}}  \|\varrho\|_{H^{r_{3}}_{x}L^{2}} \|f\|_{L^{2}_{x}L^{2}_{1/2}},
\\ \label{Gamma-3-ub-weighted-dot-1}
|(\tilde{\Gamma}_{2}^{\lambda,T}(g,h), f)|
\lesssim C_{3,\lambda,T} \|\nabla_{x}g\|_{L^{2}_{x}L^{2}}  \|\nabla_{x}h\|_{L^{2}_{x}L^{2}_{1/2}}  \|\nabla_{x}\varrho\|_{L^{2}_{x}L^{2}} \|f\|_{L^{2}_{x}L^{2}_{1/2}},
\een
where $r_{1}, r_{2}, r_{3} \in \mathbb{N}, 0 \leq  r_{1}, r_{2}, r_{3} \leq 2,  r_{1}+ r_{2}+ r_{3}=4$.

Based on \eqref{Gamma-3-ub-full-space} and \eqref{Gamma-3-ub-weighted-dot-1}, by making a suitable choice of parameters $r_{1}, r_{2}, r_{3}$ to deal with different distribution of derivative order, we get the following estimate.

\begin{thm} \label{Gamma-3-energy-estimate-label}  Let $N \geq 2$. It holds that
\ben  \label{Gamma-3-energy-estimate}
|\sum_{|\alpha|\leq N} ( \partial^{\alpha} \tilde{\Gamma}_{3}^{\lambda,T}(f,f,f),  \partial^{\alpha}f )| \lesssim C_{3,\lambda,T}\left( \|f\|_{H^{2}_{x}L^{2}}^{2}\mathcal{D}_{N}(f) + \mathrm{1}_{N \geq 3}C_{N}
\|f\|_{H^{N-1}_{x}L^{2}}\|f\|_{H^{N}_{x}L^{2}}\mathcal{D}^{\f{1}{2}}_{N-1}(f)\mathcal{D}^{\f{1}{2}}_{N}(f) \right).
\een
\end{thm}

\begin{proof} By binomial formula, we need to consider $(\tilde{\Gamma}_{3}^{\lambda,T}(\partial^{\alpha_{1}}f,\partial^{\alpha_{2}}f,\partial^{\alpha_{3}}f), \partial^{\alpha}f)$ for all combinations of $\alpha_{1}+\alpha_{2}+\alpha_{3}=\alpha$ with $|\alpha|\leq N$. For $|\alpha|=0$, it is easy to check
\beno
(\tilde{\Gamma}_{3}^{\lambda,T}(f, f, f), f) =
(\tilde{\Gamma}_{3}^{\lambda,T}(f, f, f), f_{2}).
\eeno
By taking \eqref{Gamma-3-ub-weighted-dot-1}  and \eqref{more-than-1-order-bounded-by-dissipation},  we get 
\beno
|(\tilde{\Gamma}_{3}^{\lambda,T}(f, f, f), f_{2})|
\lesssim C_{3,\lambda,T} \|\nabla_{x}f\|_{L^{2}_{x}L^{2}}^2  \|\nabla_{x}f\|_{L^{2}_{x}L^{2}_{1/2}}   \|f _2\|_{L^{2}_{x}L^{2}_{1/2}} \lesssim C_{3,\lambda,T} \|f\|_{H^{1}_{x}L^{2}}^{2} \mathcal{D}_{1}(f).
\eeno

Now it remains to consider $1 \leq |\alpha| \leq N$.
By  \eqref{Gamma-3-ub-full-space}, it suffices to prove that for all combinations of $\alpha_{1}+\alpha_{2}+\alpha_{3}=\alpha$ with $1 \leq |\alpha| \leq N$,
the following inequality
\beno
&&\|\partial^{\alpha_{1}}f\|_{H^{r_{1}}_{x}L^{2}}  \|\partial^{\alpha_{2}}f\|_{H^{r_{2}}_{x}L^{2}_{1/2}}  \|\partial^{\alpha_{3}}f\|_{H^{r_{3}}_{x}L^{2}} \|\partial^{\alpha}f\|_{L^{2}_{x}L^{2}_{1/2}}
\\&\lesssim& \|f\|_{H^{2}_{x}L^{2}}^{2}\mathcal{D}_{N}(f) + \mathrm{1}_{N \geq 3}C_{N}
\|f\|_{H^{N-1}_{x}L^{2}}\|f\|_{H^{N}_{x}L^{2}}\mathcal{D}^{\f{1}{2}}_{N-1}(f)\mathcal{D}^{\f{1}{2}}_{N}(f).
\eeno
holds for some   $r_{1}, r_{2}, r_{3} \in \mathbb{N}$ verifying $0 \leq  r_{1}, r_{2}, r_{3} \leq 2,  r_{1}+ r_{2}+ r_{3}=4$.
If $\alpha_{2}=\alpha$,
Taking $(r_{1},r_{2},r_{3}) = (2,0,2)$ and using \eqref{more-than-1-order-bounded-by-dissipation}, we have
\beno
\|f\|_{H^{2}_{x}L^{2}}  \|\partial^{\alpha}f\|_{H^{0}_{x}L^{2}_{1/2}} \|f\|_{H^{2}_{x}L^{2}} \|\partial^{\alpha}f\|_{L^{2}_{x}L^{2}_{1/2}}
\lesssim \|f\|_{H^{2}_{x}L^{2}}^{2}\mathcal{D}_{N}(f)
\eeno
Therefore in the following we only need to consider all combinations of $\alpha_{1}+\alpha_{2}+\alpha_{3}=\alpha$ with $1 \leq |\alpha| \leq N, |\alpha_{2}| \leq |\alpha|-1$. Note that $\alpha_{1}$ and $\alpha_{3}$ play the same role. We can always assume $|\alpha_{1}| \leq |\alpha_{3}|$.

The following is divided into three cases: $|\alpha| = 1; |\alpha| = 2; 3 \leq |\alpha| \leq  N$.

{\it {Case 1: $|\alpha| = 1$.}} In this case, there remains only one case: $\alpha_{1} = \alpha_{2} = 0, \alpha_{3} = \alpha$. Taking $(r_{1},r_{2},r_{3}) = (2,2,0)$, using \eqref{f1-any-order-bounded-by-energy}, \eqref{any-order-bounded-by-energy} and \eqref{more-than-1-order-bounded-by-dissipation}, we have
\ben \label{tri-first-full-second-0}
&& 
\|f\|_{H^{2}_{x}L^{2}}  \|f\|_{H^{2}_{x}L^{2}_{1/2}} \|\partial^{\alpha}f\|_{H^{0}_{x}L^{2}} \|\partial^{\alpha}f\|_{L^{2}_{x}L^{2}_{1/2}}
\\ \nonumber &\lesssim&  
\|f\|_{H^{2}_{x}L^{2}}  (\|f_1\|_{H^{2}_{x}L^{2}_{1/2}} + \|f_2\|_{H^{2}_{x}L^{2}_{1/2}}) \|\partial^{\alpha}f\|_{H^{0}_{x}L^{2}} \|\partial^{\alpha}f\|_{L^{2}_{x}L^{2}_{1/2}}
 \\ \nonumber &\lesssim&   \|f\|_{H^{2}_{x}L^{2}} \|f\|_{H^{2}_{x}L^{2}} \mathcal{D}^{\f{1}{2}}_{1}(f) \mathcal{D}^{\f{1}{2}}_{1}(f) + 
\|f\|_{H^{2}_{x}L^{2}} \mathcal{D}^{\f{1}{2}}_{2}(f) \|f\|_{H^{1}_{x}L^{2}}   \mathcal{D}^{\f{1}{2}}_{1}(f)
\lesssim \|f\|_{H^{2}_{x}L^{2}}^{2} \mathcal{D}_{2}(f).
\een

{\it {Case 2: $|\alpha| = 2$.}}
In this case, it remains to consider two subcases: $|\alpha_{2}| = 0, |\alpha_{1}| \leq 1;  |\alpha_{2}| = 1, |\alpha_{1}|=0$. In the first  subcase $|\alpha_{2}| = 0, |\alpha_{1}| \leq 1 \leq |\alpha_{3}|, |\alpha_{1}| + |\alpha_{3}| \leq 2$. Taking $(r_{1},r_{2},r_{3}) = (2-|\alpha_{1}|,2,|\alpha_{1}|)$ and using the same argument as in \eqref{tri-first-full-second-0}, we will get
\beno
\|\partial^{\alpha_{1}}f\|_{H^{2-|\alpha_{1}|}_{x}L^{2}}  \|f\|_{H^{2}_{x}L^{2}_{1/2}}  \|\partial^{\alpha_{3}}f\|_{H^{|\alpha_{1}|}_{x}L^{2}} \|\partial^{\alpha}f\|_{L^{2}_{x}L^{2}_{1/2}}
\lesssim  \|f\|_{H^{2}_{x}L^{2}}^{2}\mathcal{D}_{2}(f).
\eeno
In the second subcase $|\alpha_{2}| = |\alpha_{3}| = 1, |\alpha_{1}| =0$. Taking $(r_{1},r_{2},r_{3}) = (2,1,1)$ and using \eqref{more-than-1-order-bounded-by-dissipation}, we have
\beno
\|f\|_{H^{2}_{x}L^{2}}  \|\partial^{\alpha_{2}}f\|_{H^{1}_{x}L^{2}_{1/2}}  \|\partial^{\alpha_{3}}f\|_{H^{1}_{x}L^{2}} \|\partial^{\alpha}f\|_{L^{2}_{x}L^{2}_{1/2}}
\lesssim  \|f\|_{H^{2}_{x}L^{2}}^{2}\mathcal{D}_{2}(f).
\eeno

{\it {Case 3: $3 \leq |\alpha| \leq N$.}}
In this case, it remains to consider five subcases: $|\alpha_{2}| = |\alpha|-1; |\alpha_{2}| = |\alpha|-2;
1 \leq |\alpha_{2}| \leq |\alpha|-3;
|\alpha_{2}| = 0, |\alpha_{1}| \leq 2; |\alpha_{2}| = 0, |\alpha_{1}| \geq 3$.
In the first subcase, $|\alpha_{2}| = |\alpha|-1 \geq 2$, then $|\alpha_{1}|=1, |\alpha_{3}| = 0$.
Taking $(r_{1},r_{2},r_{3}) = (1,1,2)$ and using \eqref{more-than-1-order-bounded-by-dissipation}, we get
\beno
\|\partial^{\alpha_{1}}f\|_{H^{1}_{x}L^{2}}  \|\partial^{\alpha_{2}}f\|_{H^{1}_{x}L^{2}_{1/2}}  \|f\|_{H^{2}_{x}L^{2}} \|\partial^{\alpha}f\|_{L^{2}_{x}L^{2}_{1/2}}
\lesssim \|f\|_{H^{2}_{x}L^{2}}^{2}\mathcal{D}_{N}(f).
\eeno
In the second subcase, $|\alpha_{2}| = |\alpha|-2 \geq 1$, then $|\alpha_{1}| + |\alpha_{3}| \leq 2$.
Taking $(r_{1},r_{2},r_{3}) = (2-|\alpha_{1}|,2,|\alpha_{1}|)$ and using \eqref{more-than-1-order-bounded-by-dissipation}, we get
\beno
\|\partial^{\alpha_{1}}f\|_{H^{2-|\alpha_{1}|}_{x}L^{2}}  \|\partial^{\alpha_{2}}f\|_{H^{2}_{x}L^{2}_{1/2}}  \|\partial^{\alpha_{3}}f\|_{H^{|\alpha_{1}|}_{x}L^{2}} \|\partial^{\alpha}f\|_{L^{2}_{x}L^{2}_{1/2}}
\lesssim \|f\|_{H^{2}_{x}L^{2}}^{2}\mathcal{D}_{N}(f).
\eeno
In the third subcase, $1 \leq |\alpha_{2}| \leq |\alpha|-3$, then $N \geq 4, |\alpha_{1}| \leq \f{N}{2} \leq |\alpha_{3}| \leq N-1$. Note that $\f{N}{2} + 2 \leq N$.
Taking $(r_{1},r_{2},r_{3}) = (2,2,0)$ and using \eqref{more-than-1-order-bounded-by-dissipation},
we get
\beno
\|\partial^{\alpha_{1}}f\|_{H^{2}_{x}L^{2}}  \|\partial^{\alpha_{2}}f\|_{H^{2}_{x}L^{2}_{1/2}}  \|\partial^{\alpha_{3}}f\|_{H^{0}_{x}L^{2}} \|\partial^{\alpha}f\|_{L^{2}_{x}L^{2}_{1/2}}
\lesssim \|f\|_{H^{N}_{x}L^{2}}\mathcal{D}^{\f{1}{2}}_{N-1}(f)\|f\|_{H^{N-1}_{x}L^{2}}\mathcal{D}^{\f{1}{2}}_{N}(f).
\eeno

In the fourth subcase $|\alpha_{2}| = 0, |\alpha_{1}| \leq 2,  |\alpha_{3}| = |\alpha|-|\alpha_{1}| \geq 1$. Taking $(r_{1},r_{2},r_{3}) = (2-|\alpha_{1}|,2,|\alpha_{1}|)$, similarly to \eqref{tri-first-full-second-0}, we get
\beno
\|\partial^{\alpha_{1}}f\|_{H^{2-|\alpha_{1}|}_{x}L^{2}}  \|f\|_{H^{2}_{x}L^{2}_{1/2}}  \|\partial^{\alpha_{3}}f\|_{H^{|\alpha_{1}|}_{x}L^{2}} \|\partial^{\alpha}f\|_{L^{2}_{x}L^{2}_{1/2}}
\lesssim \|f\|_{H^{2}_{x}L^{2}}^{2} \mathcal{D}_{N}(f) +
\|f\|_{H^{2}_{x}L^{2}} \|f\|_{H^{N}_{x}L^{2}}  \mathcal{D}^{\f{1}{2}}_{2}(f) \mathcal{D}^{\f{1}{2}}_{N}(f).
\eeno

In the fifth subcase $|\alpha_{2}| = 0, |\alpha_{1}| \geq 3$. Note that $N \geq 6, 3 \leq |\alpha_{1}| \leq \f{N}{2} \leq |\alpha_{3}| \leq N-3$.
Taking $(r_{1},r_{2},r_{3}) = (0,2,2)$, similarly to \eqref{tri-first-full-second-0}, we get
\beno
\|\partial^{\alpha_{1}}f\|_{H^{0}_{x}L^{2}}  \|f\|_{H^{2}_{x}L^{2}_{1/2}}  \|\partial^{\alpha_{3}}f\|_{H^{2}_{x}L^{2}} \|\partial^{\alpha}f\|_{L^{2}_{x}L^{2}_{1/2}}
\lesssim   \|f\|_{H^{N-1}_{x}L^{2}}\|f\|_{H^{N}_{x}L^{2}}\mathcal{D}^{\f{1}{2}}_{N-1}(f)\mathcal{D}^{\f{1}{2}}_{N}(f).
\eeno

Patching together all the above estimates, observing that {\it {Case 3}} happens only if $N \geq 3$,
we obtain  \eqref{Gamma-3-energy-estimate}.
\end{proof}

By revising the above proof of  Theorem \ref{Gamma-2-energy-estimate-label}  and Theorem \ref{Gamma-3-energy-estimate-label}, we can similarly derive the following two results. Let $N \geq 2$. Let $\psi$ be a function of the variable $v$ satisfying $|\psi|_{L^{2}_{1/2}} \lesssim 1$, then
it holds that
\ben  \label{Gamma-2-energy-estimate-with-13}
\sum_{|\alpha|\leq N} |\langle \partial^{\alpha} \tilde{\Gamma}_{2}^{\lambda,T}(f,f),  \psi \rangle|_{L^{2}_{x}}^{2} \lesssim C_{2,\lambda,T}^{2}\left( \|f\|_{H^{2}_{x}L^{2}}^{2}\mathcal{D}_{N}(f) + \mathrm{1}_{N \geq 3}C_{N}
\|f\|_{H^{N}_{x}L^{2}}^{2}\mathcal{D}_{N-1}(f) \right),
\\
 \label{Gamma-3-energy-estimate-with-13}
\sum_{|\alpha|\leq N} |\langle \partial^{\alpha} \tilde{\Gamma}_{3}^{\lambda,T}(f,f,f),  \psi \rangle|_{L^{2}_{x}}^{2} \lesssim
C_{3,\lambda,T}^{2}\left( \|f\|_{H^{2}_{x}L^{2}}^{4}\mathcal{D}_{N}(f) + \mathrm{1}_{N \geq 3}C_{N}
\|f\|_{H^{N-1}_{x}L^{2}}^{2}\|f\|_{H^{N}_{x}L^{2}}^{2}\mathcal{D}_{N-1}(f) \right).
\een

\section{A priori estimate and global well-posedness}   \label{global}
This section is devoted to the proof to Theorem \ref{global-well-posedness}.
For fixed $\lambda,T>0, 0<\epsilon<1$, it is not difficult to derive local existence of the problem \eqref{T-half-scaled-more-consistent}. So we focus on the uniform-in-$\epsilon$ 
{\it a priori} estimate for the equation \eqref{T-half-scaled-more-consistent}. This estimate is given as Theorem \ref{a-priori-estimate-LBE}. Then by a standard continuity argument the global well-posedness result in
Theorem \ref{global-well-posedness} can be established.

\subsection{A priori estimate of a general equation} This subsection is devoted to some uniform-in-$\epsilon$ {\it{a priori}} estimate of the following equation
\ben \label{lBE}
 \partial _t f +  \f{1}{\epsilon}T^{1/2} v \cdot \nabla_{x} f + \f{1}{\epsilon^{2}}\tilde{\mathcal{L}}^{\lambda,T}f = g,  \quad t>0, x \in \mathbb{R}^{3}, v \in \mathbb{R}^{3},
\een
where $g$ is a given function.

Let $f$ be a solution to \eqref{lBE}.
Recalling the formula \eqref{linear-combination-of-basis} of projection operator $\tilde{\mathbb{P}}_{\lambda}$, we denote
\ben \label{definition-f-1} (\tilde{\mathbb{P}}_{\lambda} f)(t,x,v) = (a(t,x) + b(t,x) \cdot v + c(t,x)|v|^{2})N,\een
where $N=N_{\lambda}$ and
\ben \label{a-b-c-for-t-x}
(a(t,x), b(t,x), c(t,x)) \colonequals      (a^{f(t,x,\cdot)}_{\lambda}, b^{f(t,x,\cdot)}_{\lambda}, c^{f(t,x,\cdot)}_{\lambda}).
\een
Here we recall \eqref{explicit-defintion-of-abc} for the definition of $(a^{f(t,x,\cdot)}_{\lambda}, b^{f(t,x,\cdot)}_{\lambda}, c^{f(t,x,\cdot)}_{\lambda})$ for fixed $t,x$.
Note that in \eqref{a-b-c-for-t-x} we omit $f, \lambda, T$ for brevity. However, we should
always keep in mind that $(a, b, c)$ are functions of $(t,x)$ originating from the solution $f$ to \eqref{lBE} for fixed $\lambda, T$.

We set $f_{1} \colonequals      \tilde{\mathbb{P}}_{\lambda} f$ and $f_{2} \colonequals      f - \tilde{\mathbb{P}}_{\lambda} f$.
We first recall some basics of macro-micro decomposition. Plugging the macro-micro decomposition $f = f_{1} + f_{2}$ into \eqref{lBE}
and using the fact $ \tilde{\mathcal{L}}^{\lambda,T} f_{1} =  0$, we get
\ben \label{macro-micro-LBE-2} \epsilon\partial_{t}f_{1} + T^{1/2}v\cdot \nabla_{x} f_{1}  = - \epsilon\partial_{t}f_{2} - T^{1/2}v\cdot \nabla_{x} f_{2} - \f{1}{\epsilon}\tilde{\mathcal{L}}^{\lambda,T}f_{2} + \epsilon g.\een

Recalling \eqref{definition-f-1},
the left-hand of \eqref{macro-micro-LBE-2} reads
\ben \label{left-write-out}
 (\epsilon\partial_{t}a + \sum_{i=1}^{3}\epsilon\partial_{t}b_{i} v_{i} + \epsilon\partial_{t}c |v|^{2} )N +  T^{1/2}(\sum_{i=1}^{3} \partial_{i}a v_{i} + \sum_{i \neq j} \partial_{i}b_{j} v_{i}v_{j} + \sum_{i=1}^{3} \partial_{i}c v_{i}  |v|^{2} )N.
\een
Here $\partial_{i} = \partial_{x_{i}}$ for $i=1,2,3, b=(b_{1},b_{2},b_{3})$ and $v=(v_{1},v_{2},v_{3})$.
We order the 13 functions of $v$ in \eqref{left-write-out} as
\beno e_{1} = N, \quad e_{2} = v_{1}N, \quad e_{3} = v_{2}N, \quad e_{4} = v_{3}N, \quad e_{5} = v_{1}^{2}N, \quad e_{6} = v_{2}^{2}N, \quad e_{7} = v_{3}^{2}N,  \\ e_{8} =v_{1}v_{2}N, \quad e_{9} = v_{2}v_{3}N, \quad e_{10} = v_{3}v_{1}N, \quad  e_{11} = |v|^{2}v_{1}N, \quad e_{12} = |v|^{2}v_{2}N, \quad e_{13} = |v|^{2}v_{3}N.
\eeno
We emphasize that $e_{i}$ depends on $\lambda$ through $N = N_{\lambda}$. We also order the 13 functions of $(t, x)$ in \eqref{left-write-out} as
\beno x_{1} \colonequals     \epsilon\partial_{t}a, \quad x_{2} \colonequals      \epsilon\partial_{t}b_{1}+ T^{1/2}\partial_{1} a, \quad x_{3} \colonequals      \epsilon\partial_{t}b_{2}+ T^{1/2}\partial_{2} a, \quad x_{4} \colonequals      \epsilon\partial_{t}b_{3}+ T^{1/2}\partial_{3} a,
\\ x_{5} \colonequals      \epsilon\partial_{t}c+ T^{1/2}\partial_{1} b_{1}, \quad x_{6} \colonequals      \epsilon\partial_{t}c+ T^{1/2}\partial_{2} b_{2}, \quad x_{7} \colonequals      \epsilon\partial_{t}c+ T^{1/2}\partial_{3} b_{3},
\\ x_{8} \colonequals     T^{1/2}(\partial_{1}b_{2} + \partial_{2}b_{1}), \quad x_{9} \colonequals     T^{1/2}(\partial_{2}b_{3} + \partial_{3}b_{2}), \quad x_{10} \colonequals     T^{1/2}(\partial_{3}b_{1} + \partial_{1}b_{3}), \\ x_{11} \colonequals      T^{1/2}\partial_{1}c, \quad x_{12} \colonequals     T^{1/2}\partial_{2}c, \quad x_{13} \colonequals      T^{1/2}\partial_{3}c.
\eeno
Use $\mathsf{T}$ to denote vector transpose. For simplicity, we define two column vectors
\beno
E \colonequals      (e_{1}, \cdots, e_{13})^{\mathsf{T}}, \quad X \colonequals      (x_{1}, \cdots, x_{13})^{\mathsf{T}}.
\eeno
With these two column vectors, \eqref{left-write-out}  is $\epsilon\partial_{t}f_{1} + T^{1/2}v\cdot \nabla_{x} f_{1} = E^{\mathsf{T}} X$ and thus
 \eqref{macro-micro-LBE-2} is written as
\beno
 E^{\mathsf{T}} X = - \epsilon\partial_{t}f_{2} - T^{1/2}v\cdot \nabla_{x} f_{2} - \f{1}{\epsilon}\tilde{\mathcal{L}}^{\lambda,T}f_{2} + \epsilon g.
\eeno
Taking inner product with $E$
in the space $L^{2}(\mathbb{R}^{3})$, since $X$ depends on $(t,x)$ but not on $v$,
we get
\beno
\langle E, E^{\mathsf{T}} \rangle  X =  \langle E, - \epsilon\partial_{t}f_{2} - T^{1/2}v\cdot \nabla_{x} f_{2} - \f{1}{\epsilon}\tilde{\mathcal{L}}^{\lambda,T}f_{2} + \epsilon g \rangle.
\eeno
We will see soon that the $13 \times 13$ matrix $\langle E, E^{\mathsf{T}} \rangle = (\langle e_{i}, e_{j} \rangle)_{1\leq i \leq 13, 1\leq j \leq 13}$ is invertible and so
\ben  \label{macroscopic-equations}
  X =  (\langle E, E^{\mathsf{T}} \rangle)^{-1} \langle E, - \epsilon\partial_{t}f_{2} - T^{1/2}v\cdot \nabla_{x} f_{2} - \f{1}{\epsilon}\tilde{\mathcal{L}}^{\lambda,T}f_{2} + \epsilon g \rangle.
\een

We now prove that $\langle E, E^{\mathsf{T}} \rangle$ is invertible and give some estimate on its inverse. 

\begin{lem} \label{estimate-of-combination} Let $\lambda>0$.
	The matrix $\langle E, E^{\mathsf{T}} \rangle$ is invertible and
	\ben \label{the-13-vector-l2}
	|(\langle E, E^{\mathsf{T}} \rangle)^{-1} E|_{L^{2}_{2}} \lesssim e^{\lambda/2}.
	\een
\end{lem}
\begin{proof} Recalling \eqref{moment-k-of-N-square}, 
	by directly computing $\langle e_i, e_j \rangle$ for $1 \leq i,j \leq 13$, we obtain
\beno
\langle E, E^{\mathsf{T}} \rangle =
\begin{pmatrix}
m_{0} & 0_{1 \times 3} & \f{m_{2}}{3} 1_{1 \times 3} & 0_{1 \times 3}& 0_{1 \times 3}\\
0_{3 \times 1} & \f{m_{2}}{3}I_{3 \times 3} & 0_{3 \times 3}& 0_{3 \times 3}& \f{m_{4}}{3}I_{3 \times 3}\\
\f{m_{2}}{3}1_{3 \times 1} & 0_{3 \times 3} & \f{m_{4}}{15}A& 0_{3 \times 3}& 0_{3 \times 3}\\
0_{3 \times 1} & 0_{3 \times 3} & 0_{3 \times 3}& \f{m_{4}}{15} I_{3 \times 3}& 0_{3 \times 3}\\
0_{3 \times 1} & \f{m_{4}}{3}I_{3 \times 3} & 0_{3 \times 3}& 0_{3 \times 3}& \f{m_{6}}{3} I_{3 \times 3}
\end{pmatrix},
\eeno
where
\beno
A =
\begin{pmatrix}
3 & 1 & 1\\
1 & 3 & 1\\
1 & 1 & 3
\end{pmatrix}.
\eeno
Note that $\langle E, E^{\mathsf{T}} \rangle$ is represented as a $5 \times 5$ block matrix. We then calculate
the determinate of $\langle E, E^{\mathsf{T}} \rangle$ and end with
\beno
\det(\langle E, E^{\mathsf{T}} \rangle) = \f{4m_{4}^{5}(m_{0}m_{4}-m_{2}^{2})(m_{2}m_{6}-m_{4}^{2})^{3}}{1660753125}.
\eeno
By the Cauchy-Schwarz inequality,  $\det(\langle E, E^{\mathsf{T}} \rangle)>0$ and so  $\langle E, E^{\mathsf{T}} \rangle$ is invertible.
Then we calculate the inverse and find
\beno
(\langle E, E^{\mathsf{T}} \rangle)^{-1} =
\begin{pmatrix}
\f{m_{4}}{m_{0}m_{4}-m_{2}^{2}} & 0_{1 \times 3} & \f{m_{2}}{m_{0}m_{4}-m_{2}^{2}} 1_{1 \times 3} & 0_{1 \times 3}& 0_{1 \times 3}\\
0_{3 \times 1} & \f{3m_{6}}{m_{2}m_{6}-m_{4}^{2}}I_{3 \times 3} & 0_{3 \times 3}& 0_{3 \times 3}& \f{3m_{4}}{m_{2}m_{6}-m_{4}^{2}}I_{3 \times 3}\\
\f{m_{2}}{m_{0}m_{4}-m_{2}^{2}}1_{3 \times 1} & 0_{3 \times 3} & A& 0_{3 \times 3}& 0_{3 \times 3}\\
0_{3 \times 1} & 0_{3 \times 3} & 0_{3 \times 3}& \f{15}{m_{4}} I_{3 \times 3}& 0_{3 \times 3}\\
0_{3 \times 1} & \f{3m_{4}}{m_{2}m_{6}-m_{4}^{2}}I_{3 \times 3} & 0_{3 \times 3}& 0_{3 \times 3}& \f{3m_{2}}{m_{2}m_{6}-m_{4}^{2}} I_{3 \times 3}
\end{pmatrix},
\eeno
where
\beno
A =
\begin{pmatrix}
a & b & b\\
b & a & b\\
b & b & a
\end{pmatrix}, \quad
a = \f{6m_{0}m_{4}-5m_{2}^{2}}{m_{4}(m_{0}m_{4}-m_{2}^{2})}, \quad
b = \f{3m_{0}m_{4}-5m_{2}^{2}}{2m_{4}(m_{0}m_{4}-m_{2}^{2})}.
\eeno
By Lemma \ref{N-lambda-l2-norm},
\ben \label{moment-order-2-4-6}
m_{0} \sim e^{-\lambda} (1-e^{-\lambda})^{-1/2}, \quad  m_{2}, m_{4}, m_{6} \sim e^{-\lambda}.
\een
Then it is easy to derive
\ben \label{cross-order-2-4-6}
 m_{0}m_{4}-m_{2}^{2} \sim e^{-2\lambda} (1-e^{-\lambda})^{-1/2}, \quad  m_{2}m_{6}-m_{4}^{2} \sim e^{-2\lambda}.
\een
Let $(\langle E, E^{\mathsf{T}} \rangle)^{-1} = (a_{ij})_{1\leq i,j \leq 13}$.
By \eqref{moment-order-2-4-6} and \eqref{cross-order-2-4-6}, for the elements in the first column or row, we have
\ben \label{row-column-1}
0 \leq a_{1j}, a_{j1} \lesssim e^{\lambda} (1-e^{-\lambda})^{1/2} \text{ for } 1 \leq j \leq 13.
\een
For the elements except the first column or row, we have
\ben \label{row-column-geq-2}
|a_{ij}|\lesssim e^{\lambda}  \text{ for } 2 \leq i,j \leq 13.
\een
Note that $
(\langle E, E^{\mathsf{T}} \rangle)^{-1} E = (\sum_{j=1}^{13} a_{ij} e_{j})_{1 \leq i \leq  13}
$,
and so for $1 \leq i \leq 13$,
\beno
| \sum_{j=1}^{13} a_{ij} e_{j} |_{L^{2}_{2}} \lesssim
 \sum_{j=1}^{13} |a_{ij} ||e_{j} |_{L^{2}_{2}} = |a_{i1}| |e_{1} |_{L^{2}_{2}} + \sum_{j=2}^{13} |a_{ij} ||e_{j} |_{L^{2}_{2}}.
\eeno
By \eqref{L2-of-N} and \eqref{row-column-1}, we have
\beno
|a_{i1}| |e_{1} |_{L^{2}_{2}} \lesssim e^{\lambda} (1-e^{-\lambda})^{1/2} \times e^{-\lambda/2} (1-e^{-\lambda})^{-1/4} = e^{\lambda/2} (1-e^{-\lambda})^{1/4}.
\eeno
By \eqref{weighted-L2-k-geq-1-of-N} and \eqref{row-column-geq-2}, we have
\beno
\sum_{j=2}^{13} |a_{ij} ||e_{j} |_{L^{2}_{2}} \lesssim e^{\lambda}  \times e^{-\lambda/2} = e^{\lambda/2}.
\eeno
Patching together the previous two estimates, we finish the proof.
\end{proof}

For simplicity, we denote the terms in the right-hand side of
\eqref{macroscopic-equations} by
\beno
\mathcal{U} = (\mathcal{U}^{(0)}, \{\mathcal{U}^{(1)}_{i}\}_{1\leq i \leq 3}, \{\mathcal{U}^{(2)}_{i}\}_{1\leq i \leq 3}, \{\mathcal{U}^{(2)}_{ij} \}_{1\leq i < j  \leq 3}, \{\mathcal{U}^{(3)}_{i}\}_{1\leq i \leq 3})^{\mathsf{T}} \colonequals    (\langle E, E^{\mathsf{T}} \rangle)^{-1} \langle E, f_{2} \rangle
\\
\mathcal{V} = (\mathcal{V}^{(0)}, \{\mathcal{V}^{(1)}_{i}\}_{1\leq i \leq 3}, \{\mathcal{V}^{(2)}_{i}\}_{1\leq i \leq 3}, \{\mathcal{V}^{(2)}_{ij} \}_{1\leq i < j  \leq 3}, \{\mathcal{V}^{(3)}_{i}\}_{1\leq i \leq 3})^{\mathsf{T}} \colonequals    (\langle E, E^{\mathsf{T}} \rangle)^{-1} \langle E, - T^{1/2}v\cdot \nabla_{x} f_{2}  \rangle
\\
\mathcal{W} = (\mathcal{W}^{(0)}, \{\mathcal{W}^{(1)}_{i}\}_{1\leq i \leq 3}, \{\mathcal{W}^{(2)}_{i}\}_{1\leq i \leq 3}, \{\mathcal{W}^{(2)}_{ij} \}_{1\leq i < j  \leq 3}, \{\mathcal{W}^{(3)}_{i}\}_{1\leq i \leq 3})^{\mathsf{T}} \colonequals    (\langle E, E^{\mathsf{T}} \rangle)^{-1}  \langle E, - \f{1}{\epsilon}\tilde{\mathcal{L}}^{\lambda,T}f_{2} \rangle
\\
\mathcal{X} = (\mathcal{X}^{(0)}, \{\mathcal{X}^{(1)}_{i}\}_{1\leq i \leq 3}, \{\mathcal{X}^{(2)}_{i}\}_{1\leq i \leq 3}, \{\mathcal{X}^{(2)}_{ij} \}_{1\leq i < j  \leq 3}, \{\mathcal{X}^{(3)}_{i}\}_{1\leq i \leq 3})^{\mathsf{T}} \colonequals    (\langle E, E^{\mathsf{T}} \rangle)^{-1} \langle E, \epsilon g \rangle
\\
\mathcal{T} = (\mathcal{T}^{(0)}, \{\mathcal{T}^{(1)}_{i}\}_{1\leq i \leq 3}, \{\mathcal{T}^{(2)}_{i}\}_{1\leq i \leq 3}, \{\mathcal{T}^{(2)}_{ij} \}_{1\leq i < j  \leq 3}, \{\mathcal{T}^{(3)}_{i}\}_{1\leq i \leq 3})^{\mathsf{T}} \colonequals    - \epsilon\partial_{t} \mathcal{U} + \mathcal{V} + \mathcal{W} + \mathcal{X}.
\eeno
Now \eqref{macroscopic-equations} is written as
\ben \label{macroscopic-system}
 X =  \mathcal{T} = - \epsilon\partial_{t} \mathcal{U} + \mathcal{V} + \mathcal{W} + \mathcal{X}.
\een

The estimates $\mathcal{U}, \mathcal{V}, \mathcal{W}, \mathcal{X}$  are given in the following lemma.
\begin{lem} \label{estimate-for-UVWX} It holds that
\beno  \sum_{|\alpha|\leq N}|\partial^{\alpha}\mathcal{U}|^{2}_{L^{2}_{x}} \leq e^{\lambda} \|f_{2}\|^{2}_{H^{N}_{x}L^{2}}, \quad \sum_{|\alpha|\leq N-1}|\partial^{\alpha}\mathcal{V}|^{2}_{L^{2}_{x}} \leq e^{\lambda} T\|f_{2}\|^{2}_{H^{N}_{x}L^{2}},  \\
\sum_{|\alpha|\leq N-1}|\partial^{\alpha}\mathcal{W}|^{2}_{L^{2}_{x}} \leq e^{\lambda} \f{1}{\epsilon^{2}} \tilde{C}_{1,\lambda,T}^{2} \|f_{2}\|_{H^{N-1}_{x}L^{2}_{1/2}}^{2}, \quad 
  \sum_{|\alpha|\leq N-1}|\partial^{\alpha}\mathcal{X}|^{2}_{L^{2}_{x}} \leq \epsilon^{2} \mathcal{C}_{\lambda, N-1}(g), \eeno
where
\beno
\mathcal{C}_{\lambda, n}(g)\colonequals \sum_{|\alpha| \leq n}
|\langle  (\langle E, E^{\mathsf{T}} \rangle)^{-1} E, \partial^{\alpha}g \rangle|^{2}_{L^{2}_{x}}.
\eeno
\end{lem}
\begin{proof}
Note that
$
\partial^{\alpha}\mathcal{U} = (\langle E, E^{\mathsf{T}} \rangle)^{-1} \langle E, \partial^{\alpha}f_{2} \rangle
$. By \eqref{the-13-vector-l2}, 
\ben \label{partial-alpha-mathcal-U}
|\partial^{\alpha}\mathcal{U}| \lesssim e^{\lambda/2}  |\partial^{\alpha}f_{2}|_{L^{2}},
\een
which gives the first inequality on $\mathcal{U}$.
Note that
$
\partial^{\alpha}\mathcal{V} = (\langle E, E^{\mathsf{T}} \rangle)^{-1} \langle E, - T^{1/2}v\cdot \nabla_{x} \partial^{\alpha} f_{2}  \rangle$.
Then a similar argument yields the second inequality on $\mathcal{V}$.

By \eqref{lower-and-upper-bound}, we have
\beno
  \langle \tilde{\mathcal{L}}^{\lambda,T}g, h \rangle = \langle \tilde{\mathcal{L}}^{\lambda,T}g_{2}, h_{2} \rangle \leq  (\langle \tilde{\mathcal{L}}^{\lambda,T}g_{2}, g_{2} \rangle)^{1/2} (\langle \tilde{\mathcal{L}}^{\lambda,T}h_{2}, h_{2} \rangle)^{1/2} \leq \tilde{C}_{1,\lambda,T}|g_{2}|_{L^{2}_{1/2}} |h_{2}|_{L^{2}_{1/2}}.
\eeno
For any $g$ and $k \geq 0$, one can derive that $|g_{2}|_{L^{2}_{k}} \lesssim |g|_{L^{2}_{k}}, |g_{1}|_{L^{2}_{k}} \lesssim |g|_{L^{2}}$.
From this together with \eqref{the-13-vector-l2}, we have
\beno
|\partial^{\alpha}\mathcal{W}| = |(\langle E, E^{\mathsf{T}} \rangle)^{-1}  \langle E, - \f{1}{\epsilon}\tilde{\mathcal{L}}^{\lambda,T}\partial^{\alpha}f_{2} \rangle| \lesssim \f{1}{\epsilon} e^{\lambda/2}\tilde{C}_{1,\lambda,T}
|\partial^{\alpha}f_{2}|_{L^{2}_{1/2}},
\eeno
which gives
 the third inequality on $\mathcal{W}$.
The last result on $\mathcal{X}$ is trivial.
\end{proof}

We now estimate the dynamics of $(a,b,c)$ in the following lemma.
\begin{lem}\label{estimate-for-ptabc} Recall \eqref{defintion-of-l-i} and \eqref{explicit-defintion-of-abc}.
	It holds that
\ben \label{simple-estimate-of-pt-a}
\epsilon^{2}\sum_{|\alpha|\leq N-1}|\partial^{\alpha}\partial_{t}a|^{2}_{L^{2}_{x}} \lesssim T e^{\lambda}\|f_{2}\|_{H^{N}_{x}L^{2}}^{2}
+ \epsilon^{2}
\mathcal{P}_{\lambda, N-1}(g),
\\ \label{estimate-of-pt-b-c}
\epsilon^{2}\sum_{|\alpha|\leq N-1}|\partial^{\alpha}\partial_{t} (b,c)|^{2}_{L^{2}_{x}} \lesssim T |\nabla_{x}(a,b,c)|^{2}_{H^{N-1}_{x}}
+ T e^{\lambda}\|f_{2}\|_{H^{N}_{x}L^{2}}^{2} + \epsilon^{2} \mathcal{P}_{\lambda, N-1}(g),
\een
where
\beno
\mathcal{P}_{\lambda, n}(g)\colonequals \sum_{|\alpha| \leq n}
(|\langle  l_{1}N -  l_{2}N|v|^{2}, \partial^{\alpha}g \rangle|^{2}_{L^{2}_{x}} + |\langle  l_{3}N v, \partial^{\alpha}g \rangle|^{2}_{L^{2}_{x}} + |\langle  l_{4} N|v|^{2} - l_{2} N, \partial^{\alpha}g \rangle|^{2}_{L^{2}_{x}}).
\eeno
\end{lem}
\begin{proof} 
Taking inner products between equation \eqref{lBE} and the functions $l_{1}N -  l_{2}N|v|^{2}, l_{3}N v, l_{4} N|v|^{2} - l_{2} N$, using $
\langle  \tilde{\mathcal{L}}^{\lambda,T}f, N \rangle = \langle  \tilde{\mathcal{L}}^{\lambda,T}f, N v \rangle = \langle  \tilde{\mathcal{L}}^{\lambda,T}f, N|v|^{2} \rangle =0,
$
we get
\beno \epsilon\partial_{t} a +  \langle T^{1/2}v\cdot \nabla_{x} f, l_{1}N -  l_{2}N|v|^{2} \rangle = \langle \epsilon g, l_{1}N -  l_{2}N|v|^{2}\rangle, \\ \epsilon\partial_{t} b +  \langle T^{1/2}v\cdot \nabla_{x} f, l_{3}N v  \rangle = \langle \epsilon g,  l_{3}N v \rangle,
\\ \epsilon\partial_{t} c + \langle T^{1/2}v\cdot \nabla_{x} f, l_{4} N|v|^{2} - l_{2} N  \rangle = \langle \epsilon g,  l_{4} N|v|^{2} - l_{2} N \rangle.
\eeno
Since $\langle v_{i}N, N\rangle = \langle v_{i} v_{j}N, v_{k}N\rangle = \langle v_{i}|v|^{2}N, |v|^{2}N\rangle = 0$ for $i, j, k \in \{1,2,3\}$.
Recalling the definition of $l_{1}, l_{2}$ in \eqref{defintion-of-l-i}, it is straightforward to see
\beno
\langle v\cdot \nabla_{x} f_{1}, l_{1}N -  l_{2}N|v|^{2} \rangle = \f{1}{3} \langle N|v|^{2}, l_{1}N -  l_{2}N|v|^{2}\rangle
\nabla_{x} \cdot b =  \f{1}{3} (m_{2}l_{1} - m_{4}l_{2})
\nabla_{x} \cdot b =
0,
\\
\quad  \langle v\cdot \nabla_{x} f_{1}, l_{3}N v  \rangle = \langle v \cdot \nabla_{x} (a N + c N |v|^{2}), l_{3}N v  \rangle = \langle  Nv_{i}, l_{3}N v_{i}  \rangle \nabla_{x} a + \langle  Nv_{i}, l_{3}N|v|^{2} v_{i}  \rangle \nabla_{x} c
\\= \f{m_{2}}{3} l_{3}\nabla_{x} a + \f{m_{4}}{3} l_{3} \nabla_{x} c = \nabla_{x} a + \f{m_{4}}{m_{2}} \nabla_{x} c,
\\ \langle v\cdot \nabla_{x} f_{1}, l_{4} N|v|^{2} - l_{2} N  \rangle = \langle v \cdot \nabla_{x} (b \cdot v N), l_{4} N|v|^{2} - l_{2} N  \rangle
\\= \langle N v_{i}^{2}, l_{4} N|v|^{2} - l_{2} N  \rangle \nabla_{x} \cdot b = \f{1}{3}(m_{4}l_{4}-m_{2}l_{2}) \nabla_{x} \cdot b = \f{1}{3} \nabla_{x} \cdot b,
\eeno
which gives the local conservation laws
\ben \label{local-conservation-laws} \left\{ \begin{aligned}
 & \epsilon\partial_{t} a = \langle \epsilon g - T^{1/2}v\cdot \nabla_{x} f_{2}, l_{1}N -  l_{2}N|v|^{2}\rangle,
\\&
\epsilon\partial_{t} b +  T^{1/2} \nabla_{x} a + T^{1/2} \f{m_{4}}{m_{2}} \nabla_{x} c
 = \langle \epsilon g - T^{1/2}v\cdot \nabla_{x} f_{2},  l_{3}N v \rangle,
\\&
\epsilon\partial_{t} c + T^{1/2} \f{1}{3} \nabla_{x} \cdot b = \langle \epsilon g - T^{1/2}v\cdot \nabla_{x} f_{2},  l_{4} N|v|^{2} - l_{2} N \rangle.
\end{aligned} \right.
\een
Note that $\f{m_{4}}{m_{2}} \leq 5$.
By \eqref{moment-order-2-4-6} and \eqref{cross-order-2-4-6},  we have
\ben \label{bounds-of-l-i}
l_{1}\sim   e^{\lambda} (1-e^{-\lambda})^{1/2}, \quad
l_{2}\sim   e^{\lambda} (1-e^{-\lambda})^{1/2}, \quad
l_{3}\sim   e^{\lambda}, \quad
l_{4}\sim   e^{\lambda}.
\een
From this together with \eqref{L2-of-N} and \eqref{weighted-L2-k-geq-1-of-N},
\ben \label{l22-fives-functions}
|l_{1}N -  l_{2}N|v|^{2}|_{L^{2}_{2}} \lesssim e^{\lambda/2}, \quad  |l_{3}N v|_{L^{2}_{2}} \lesssim e^{\lambda/2}, \quad
 |l_{4} N|v|^{2} - l_{2} N|_{L^{2}_{2}} \lesssim e^{\lambda/2},
\een
which gives
\beno
|\langle v\cdot \nabla_{x} f_{2}, l_{1}N -  l_{2}N|v|^{2}\rangle|+ |\langle v\cdot \nabla_{x} f_{2}, l_{3}N v \rangle| +|\langle v\cdot \nabla_{x} f_{2}, l_{4} N|v|^{2} - l_{2} N \rangle| \lesssim  e^{\lambda/2}  |\nabla_{x} f_{2}|_{L^{2}}.\eeno
Patching together the above results, we get
 \eqref{simple-estimate-of-pt-a} and \eqref{estimate-of-pt-b-c}.
\end{proof}

By \eqref{the-13-vector-l2} and \eqref{l22-fives-functions}, $\mathcal{C}_{\lambda, n}(g)$ and $\mathcal{P}_{\lambda, n}(g)$ are bounded by some linear combination of the following quantities,
\ben \label{non-linear-functional-general-source}
|\langle  \psi(v), \partial^{\alpha}g \rangle|^{2}_{L^{2}_{x}} \text{ with } |\psi|_{L^{2}_{2}} \lesssim e^{\lambda/2}, |\alpha| \leq n.
\een

With Lemma \ref{estimate-for-UVWX} and Lemma \ref{estimate-for-ptabc}, based on the macroscopic system \eqref{macroscopic-system}, using integration by parts to balance derivative, we derive 
macroscopic dissipation in the following lemma.
\begin{lem}\label{estimate-for-highorder-abc}  Let $N \geq 1, T>0$. Let $f \in L^{\infty}([0,T]; H^{N}_{x}L^{2})$ be a solution to \eqref{lBE},
then there exists a universal constant $C > 0$ such that
\ben \label{solution-property-part2} \epsilon\frac{\mathrm{d}}{\mathrm{d}t}\mathcal{I}_{N}(f) + \f{1}{2}T^{1/2}|\nabla_{x}(a,b,c)|^{2}_{H^{N-1}_{x}} &\leq&
C T^{1/2} e^{\lambda} \|f_{2}\|_{H^{N}_{x}L^{2}}^{2}
\\ \nonumber && + C T^{-1/2} e^{\lambda} \epsilon^{-2}
 \tilde{C}_{1,\lambda,T}^{2} \|f_{2}\|_{H^{N}_{x}L^{2}_{1/2}}^{2} + C T^{-1/2}\epsilon^{2}\mathcal{Q}_{\lambda, N-1}(g),\een
where $\mathcal{I}_{N}(f)$ is defined in \eqref{interactive-INf} and satisfying
\ben \label{temporal-bounded-by-norm}
|\mathcal{I}_{N}(f)| \leq C e^{\lambda} \|f\|^{2}_{H^{N}_{x}L^{2}}.
\een
Here $\mathcal{Q}_{\lambda, n}(g)$ is bounded by some linear combination of the quantities in \eqref{non-linear-functional-general-source}.
\end{lem}
\begin{proof}  Note that \eqref{macroscopic-system} is equivalent to
\ben \label{equation-a}\epsilon\partial_{t} a  &=& \mathcal{T}^{(0)} = -\epsilon\partial_{t}\mathcal{U}^{(0)} + \mathcal{V}^{(0)} + \mathcal{W}^{(0)} + \mathcal{X}^{(0)} ,
\\ \label{equation-ba}\epsilon\partial_{t}b_{i}+ T^{1/2} \partial_{i} a  &=& \mathcal{T}^{(1)}_{i} = -\epsilon\partial_{t}\mathcal{U}^{(1)}_{i} + \mathcal{V}^{(1)}_{i} +\mathcal{W}^{(1)}_{i} +\mathcal{X}^{(1)}_{i}, ~~~~~~~~ 1\leq i \leq 3,
\\ \label{equation-cb}\epsilon\partial_{t}c+ T^{1/2} \partial_{i} b_{i}  &=& \mathcal{T}^{(2)}_{i} = -\epsilon\partial_{t}\mathcal{U}^{(2)}_{i} + \mathcal{V}^{(2)}_{i} +\mathcal{W}^{(2)}_{i} +\mathcal{X}^{(2)}_{i}, ~~~~~~~~ 1\leq i \leq 3,
\\ \label{equation-b}T^{1/2}(\partial_{i}b_{j}+ \partial_{j} b_{i})  &=& \mathcal{T}^{(2)}_{ij} = -\epsilon\partial_{t}\mathcal{U}^{(2)}_{ij} + \mathcal{V}^{(2)}_{ij} +\mathcal{W}^{(2)}_{ij} +\mathcal{X}^{(2)}_{ij},~~~~~~~~1 \leq i < j \leq 3,
\\ \label{equation-c}T^{1/2} \partial_{i}c  &=& \mathcal{T}^{(3)}_{i} = -\epsilon\partial_{t}\mathcal{U}^{(3)}_{i} + \mathcal{V}^{(3)}_{i} +\mathcal{W}^{(3)}_{i} +\mathcal{X}^{(3)}_{i}, ~~~~~~~~ 1 \leq i \leq 3.\een

The derivation of $|\nabla_{x}c|^{2}_{H^{N-1}_{x}}$ is the easiest one and so we begin with $c$.
By \eqref{equation-c}, we have
\ben \label{equation-c-itself-laplace}
-T^{1/2}\Delta_{x}c  = -\sum_{j} \partial_{j} \mathcal{T}^{(3)}_{j}.
\een
Applying $\pa^{\alpha}$ with $|\alpha| \leq N-1$ to equation \eqref{equation-c-itself-laplace}, then by taking inner product with $\pa^{\alpha} c$,   one has
\beno T^{1/2}|\nabla_{x} \pa^{\alpha}c|^{2}_{L^{2}_{x}} = \langle -\sum_{j} \partial_{j} \pa^{\alpha}\mathcal{T}^{(3)}_{j}, \pa^{\alpha} c \rangle_{x}.
\eeno
Recalling \eqref{macroscopic-system}, $\mathcal{T} = - \epsilon\partial_{t} \mathcal{U} + \mathcal{V} + \mathcal{W} + \mathcal{X}$. Then
\beno
\langle -\sum_{j} \partial_{j} \pa^{\alpha}\mathcal{T}^{(3)}_{j}, \pa^{\alpha} c \rangle_{x} = \langle \epsilon\partial_{t}\sum_{j} \partial_{j} \pa^{\alpha}\mathcal{U}^{(3)}_{j}, \pa^{\alpha} c \rangle_{x} - \langle \sum_{j} \partial_{j} \pa^{\alpha}(\mathcal{V}^{(3)}_{j}+\mathcal{W}^{(3)}_{j}+\mathcal{X}^{(3)}_{j}), \pa^{\alpha} c \rangle_{x}.
\eeno
For the term involving $\mathcal{U}^{(3)}_{j}$, interchanging the $t$-derivative and $L^{2}_{x}$ inner product,
we have
\beno
\langle \epsilon\partial_{t} \sum_{j} \partial_{j} \pa^{\alpha}\mathcal{U}^{(3)}_{j}, \pa^{\alpha} c\rangle_{x}
= - \epsilon \frac{\mathrm{d}}{\mathrm{d}t}  \mathcal{I}^{c}_{\alpha}(f) +
\mathfrak{U}^{c}_{\alpha},
\eeno
where
\ben \label{c-temperal-energy}
\mathcal{I}^{c}_{\alpha}(f) \colonequals - \sum_{j} \langle \partial_{j} \pa^{\alpha}\mathcal{U}^{(3)}_{j}, \pa^{\alpha} c\rangle_{x}, \quad
\mathfrak{U}^{c}_{\alpha} \colonequals   - \sum_{j} \langle \partial_{j} \pa^{\alpha}\mathcal{U}^{(3)}_{j}, \epsilon\partial_{t}\pa^{\alpha} c\rangle_{x}.
\een
For the term involving $\mathcal{V}^{(3)}_{j} + \mathcal{W}^{(3)}_{j} + \mathcal{X}^{(3)}_{j}$, via integrating by parts
we have
\beno
- \langle \sum_{j} \partial_{j} \pa^{\alpha}(\mathcal{V}^{(3)}_{j} + \mathcal{W}^{(3)}_{j} + \mathcal{X}^{(3)}_{j}), \pa^{\alpha} c \rangle_{x}
= \sum_{j} \langle   \pa^{\alpha}(\mathcal{V}^{(3)}_{j} + \mathcal{W}^{(3)}_{j} + \mathcal{X}^{(3)}_{j}), \partial_{j}\pa^{\alpha} c\rangle_{x}
\colonequals    \mathfrak{O}^{c}_{\alpha}.
\eeno
Patching together the above formulas, we get
\beno
T^{1/2}|\nabla_{x} \pa^{\alpha}c|^{2}_{L^{2}_{x}} + \epsilon \frac{\mathrm{d}}{\mathrm{d}t} \mathcal{I}^{c}_{\alpha}(f) = \mathfrak{U}^{c}_{\alpha} + \mathfrak{O}^{c}_{\alpha}.
\eeno
By the Cauchy-Schwartz inequality, for any $0<\eta<1$, one has
\beno  |\mathfrak{U}^{c}_{\alpha}| \leq \f{\eta}{T^{1/2}} \epsilon^{2}|\partial^{\alpha}\partial_{t}c|^{2}_{L^{2}_{x}} + \frac{CT^{1/2}}{\eta} |\nabla_{x}\partial^{\alpha}\mathcal{U}|^{2}_{L^{2}_{x}}, 
\\ |\mathfrak{O}^{c}_{\alpha}| \leq \eta T^{1/2} |\nabla_{x} \partial^{\alpha} c|^{2}_{L^{2}_{x}} + \frac{C}{\eta T^{1/2}} (|\partial^{\alpha}\mathcal{V}|^{2}_{L^{2}_{x}} + |\partial^{\alpha}\mathcal{W}|^{2}_{L^{2}_{x}} + |\partial^{\alpha}\mathcal{X}|^{2}_{L^{2}_{x}}).
\eeno
Taking sum over $|\alpha| \leq N-1$,
by Lemma \ref{estimate-for-ptabc} and \eqref{estimate-of-pt-b-c}, we get
\ben \label{dissipation-c-H-N-minus-1}
T^{1/2}|\nabla_{x}c|^{2}_{H^{N-1}_{x}} +\epsilon \frac{\mathrm{d}}{\mathrm{d}t} \sum_{|\alpha|\leq N-1} \mathcal{I}^{c}_{\alpha}(f)
\leq \mathfrak{H}_{\eta, C},
\een
where for simplicity,
\ben \label{definition-H-eta}
\mathfrak{H}_{\eta, C} &\colonequals &   \eta T^{1/2}|\nabla_{x}(a,b,c)|^{2}_{H^{N-1}_{x}} +
\eta^{-1}C T^{1/2} e^{\lambda} \|f_{2}\|_{H^{N}_{x}L^{2}}^{2}
\\ \nonumber && + \eta^{-1}C T^{-1/2} e^{\lambda} \epsilon^{-2}
 \tilde{C}_{1,\lambda,T}^{2} \|f_{2}\|_{H^{N}_{x}L^{2}_{1/2}}^{2} + \eta^{-1}C T^{-1/2}\epsilon^{2}\mathcal{Q}_{\lambda, N-1}(g).
 \een
Here and in the rest of this proof, $0<\eta<1$ is an arbitrary constant and $C$ is a universal constant that could change across different lines.

We now derive $|\nabla_{x}b|^{2}_{H^{N-1}_{x}}$.
Based on equations \eqref{equation-cb} and \eqref{equation-b}, we have
\ben \nonumber
-T^{1/2}\Delta_{x}b_{j}-T^{1/2}\partial^{2}_{j}b_{j} &=&
\sum_{i\neq j} \partial_{j} (\epsilon\partial_{t}c+ T^{1/2}\partial_{i} b_{i}) - \sum_{i \neq j}  \partial_{i} (T^{1/2}\partial_{i}b_{j}+ T^{1/2}\partial_{j} b_{i})
- 2 \partial_{j} (\epsilon\partial_{t}c+ T^{1/2}\partial_{j} b_{j})
\\ \label{equation-b-itself-2} &=&
\sum_{i\neq j} \partial_{j}\mathcal{T}^{(2)}_{i} - \sum_{i \neq j} \partial_{i} \mathcal{T}^{(2)}_{ij} - 2 \partial_{j}\mathcal{T}^{(2)}_{j}.  \een
For $|\alpha| \leq N-1$,
applying $\pa^{\alpha}$ to equation \eqref{equation-b-itself-2} for $b_{j}$, then taking inner product with $\pa^{\alpha} b_{j}$,   one has
\beno T^{1/2}|\nabla_{x} \pa^{\alpha}b_{j}|^{2}_{L^{2}_{x}} + T^{1/2}|\pa_{j} \pa^{\alpha}b_{j}|^{2}_{L^{2}_{x}} = \langle \sum_{i\neq j} \partial_{j}\pa^{\alpha}\mathcal{T}^{(2)}_{i} - \sum_{i \neq j} \partial_{i} \pa^{\alpha}\mathcal{T}^{(2)}_{ij} - 2 \partial_{j}\pa^{\alpha}\mathcal{T}^{(2)}_{j}, \pa^{\alpha} b_{j}\rangle_{x}.
\eeno
Recalling \eqref{macroscopic-system}, $\mathcal{T} = - \epsilon\partial_{t} \mathcal{U} + \mathcal{V} + \mathcal{W} + \mathcal{X}$. For the term $- \epsilon\partial_{t} \mathcal{U}$, interchanging the $t$-derivative and $L^{2}_{x}$ inner product,
we have
\beno
&&-\langle \epsilon\partial_{t}(\sum_{i\neq j} \partial_{j}\pa^{\alpha}\mathcal{U}^{(2)}_{i} - \sum_{i \neq j} \partial_{i} \pa^{\alpha}\mathcal{U}^{(2)}_{ij} - 2 \partial_{j}\pa^{\alpha}\mathcal{U}^{(2)}_{j}), \pa^{\alpha} b_{j}\rangle_{x}
= -\epsilon \frac{\mathrm{d}}{\mathrm{d}t} \mathcal{I}^{b_{j}}_{\alpha}(f)
+
\mathfrak{U}^{b_{j}}_{\alpha},
\eeno
where
\ben \label{bj-temperal-energy}
\mathcal{I}^{b_{j}}_{\alpha}(f) \colonequals    \langle \sum_{i\neq j} \partial_{j}\pa^{\alpha}\mathcal{U}^{(2)}_{i} - \sum_{i \neq j} \partial_{i} \pa^{\alpha}\mathcal{U}^{(2)}_{ij} - 2 \partial_{j}\pa^{\alpha}\mathcal{U}^{(2)}_{j}, \pa^{\alpha} b_{j}\rangle_{x},
\\ \nonumber
\mathfrak{U}^{b_{j}}_{\alpha} \colonequals    \langle \sum_{i\neq j} \partial_{j}\pa^{\alpha}\mathcal{U}^{(2)}_{i} - \sum_{i \neq j} \partial_{i} \pa^{\alpha}\mathcal{U}^{(2)}_{ij} - 2 \partial_{j}\pa^{\alpha}\mathcal{U}^{(2)}_{j}, \epsilon\partial_{t}\pa^{\alpha} b_{j}\rangle_{x}.
\een
For the term $\mathcal{V} + \mathcal{W} + \mathcal{X}$, via integrating by parts,
we have
\beno
&&\langle \sum_{i\neq j} \partial_{j}\pa^{\alpha}(\mathcal{V}^{(2)}_{i} + \mathcal{W}^{(2)}_{i} + \mathcal{X}^{(2)}_{i}) - \sum_{i \neq j} \partial_{i} \pa^{\alpha}(\mathcal{V}^{(2)}_{ij} + \mathcal{W}^{(2)}_{ij} + \mathcal{X}^{(2)}_{ij}) - 2 \partial_{j}\pa^{\alpha}(\mathcal{V}^{(2)}_{j} + \mathcal{W}^{(2)}_{j} + \mathcal{X}^{(2)}_{j}), \pa^{\alpha} b_{j}\rangle_{x}
\\&=& -\langle \sum_{i\neq j} \pa^{\alpha}(\mathcal{V}^{(2)}_{i} + \mathcal{W}^{(2)}_{i} + \mathcal{X}^{(2)}_{i}), \partial_{j}\pa^{\alpha} b_{j}\rangle_{x}
+ \sum_{i \neq j} \langle \pa^{\alpha}(\mathcal{V}^{(2)}_{ij} + \mathcal{W}^{(2)}_{ij} + \mathcal{X}^{(2)}_{ij}), \partial_{i} \pa^{\alpha} b_{j}\rangle_{x}
\\&&+ 2\langle \pa^{\alpha}(\mathcal{V}^{(2)}_{j} + \mathcal{W}^{(2)}_{j} + \mathcal{X}^{(2)}_{j}), \partial_{j} \pa^{\alpha} b_{j}\rangle_{x} \colonequals  \mathfrak{O}^{b_{j}}_{\alpha}.
\eeno
Patching together the above formulas, we get
\beno
T^{1/2}|\nabla_{x} \pa^{\alpha}b_{j}|^{2}_{L^{2}_{x}} + T^{1/2}|\pa_{j} \pa^{\alpha}b_{j}|^{2}_{L^{2}_{x}} + \epsilon \frac{\mathrm{d}}{\mathrm{d}t} \mathcal{I}^{b_{j}}_{\alpha}(f) = \mathfrak{U}^{b_{j}}_{\alpha} + \mathfrak{O}^{b_{j}}_{\alpha}.
\eeno
By the Cauchy-Schwartz inequality, one has
\beno  |\mathfrak{U}^{b_{j}}_{\alpha}| \leq \f{\eta}{T^{1/2}} \epsilon^{2}|\partial^{\alpha}\partial_{t}b_{j}|^{2}_{L^{2}_{x}} + \frac{CT^{1/2}}{\eta} |\nabla_{x}\partial^{\alpha}\mathcal{U}|^{2}_{L^{2}_{x}},
\\ |\mathfrak{O}^{b_{j}}_{\alpha}| \leq \eta T^{1/2} |\nabla_{x} \partial^{\alpha} b_{j}|^{2}_{L^{2}_{x}} + \frac{C}{\eta T^{1/2}} (|\partial^{\alpha}\mathcal{V}|^{2}_{L^{2}_{x}} + |\partial^{\alpha}\mathcal{W}|^{2}_{L^{2}_{x}} + |\partial^{\alpha}\mathcal{X}|^{2}_{L^{2}_{x}}).
\eeno
Taking sum over $1 \leq j \leq 3$, taking sum over $|\alpha| \leq N-1$,
by Lemma \ref{estimate-for-ptabc} and \eqref{estimate-of-pt-b-c}, we get
\ben \label{dissipation-b-H-N-minus-1}
 T^{1/2}|\nabla_{x}b|^{2}_{H^{N-1}_{x}} +\epsilon \frac{\mathrm{d}}{\mathrm{d}t} \sum_{|\alpha|\leq N-1}\sum_{j=1}^{3}\mathcal{I}^{b_{j}}_{\alpha}(f)
\leq \mathfrak{H}_{\eta, C}.\een

We now derive $|\nabla_{x}a|^{2}_{H^{N-1}_{x}}$.
By \eqref{equation-ba}, we have
\ben \label{equation-a-itself-laplace}
-T^{1/2}\Delta_{x}a  = \sum_{j} \partial_{j} \epsilon\partial_{t}b_{j} -\sum_{j} \partial_{j} \mathcal{T}^{(1)}_{j}.
\een
Applying $\pa^{\alpha}$ to equation \eqref{equation-a-itself-laplace} for, by taking inner product with $\pa^{\alpha} a$,   one has
\beno T^{1/2}|\nabla_{x} \pa^{\alpha}a|^{2}_{L^{2}_{x}} &=& \langle \sum_{j} \partial_{j} \pa^{\alpha} \epsilon\partial_{t}  b_{j} -\sum_{j} \partial_{j} \pa^{\alpha}\mathcal{T}^{(1)}_{j}, \pa^{\alpha} a \rangle_{x}
\\&=& \langle \sum_{j} \partial_{j} \pa^{\alpha} \epsilon\partial_{t}  (b_{j} + \mathcal{U}^{(1)}_{j}) -\sum_{j} \partial_{j} \pa^{\alpha}(\mathcal{V}^{(1)}_{j} + \mathcal{W}^{(1)}_{j} + \mathcal{X}^{(1)}_{j}), \pa^{\alpha} a \rangle_{x}.
\eeno
where we recall $\mathcal{T} = - \epsilon\partial_{t} \mathcal{U} + \mathcal{V} + \mathcal{W} + \mathcal{X}$ from \eqref{macroscopic-system}. For the term involving $\epsilon\partial_{t}  (b_{j} + \mathcal{U}^{(1)}_{j})$, interchanging the $t$-derivative and $L^{2}_{x}$ inner product,
we have
\beno
\langle \epsilon\partial_{t} \sum_{j} \partial_{j} \pa^{\alpha} (b_{j} + \mathcal{U}^{(1)}_{j}), \pa^{\alpha} a\rangle_{x}
= - \epsilon \frac{\mathrm{d}}{\mathrm{d}t} \mathcal{I}^{a}_{\alpha}(f) + \mathfrak{U}^{a}_{\alpha},
\eeno
where
\ben \label{a-temperal-energy}
\mathcal{I}^{a}_{\alpha}(f)\colonequals    -\sum_{j} \langle \partial_{j} \pa^{\alpha}b_{j}, \pa^{\alpha} a\rangle_{x} - \sum_{j} \langle \partial_{j} \pa^{\alpha} \mathcal{U}^{(1)}_{j}, \pa^{\alpha} a\rangle_{x},
\\ \nonumber
\mathfrak{U}^{a}_{\alpha} \colonequals    -
\langle \sum_{j} \partial_{j} \pa^{\alpha}b_{j}, \epsilon\partial_{t}\pa^{\alpha} a\rangle_{x} -
\langle \sum_{j} \partial_{j} \pa^{\alpha} \mathcal{U}^{(1)}_{j}, \epsilon\partial_{t}\pa^{\alpha} a\rangle_{x}.
\een
For the term involving $\mathcal{V} + \mathcal{W} + \mathcal{X}$, via integrating by parts,
we  have
\beno
-\langle \sum_{j} \partial_{j} \pa^{\alpha}(\mathcal{V}^{(1)}_{j} + \mathcal{W}^{(1)}_{j} + \mathcal{X}^{(1)}_{j}), \pa^{\alpha} a \rangle_{x} = \sum_{j} \langle  \pa^{\alpha}(\mathcal{V}^{(1)}_{j} + \mathcal{W}^{(1)}_{j} + \mathcal{X}^{(1)}_{j}), \partial_{j} \pa^{\alpha} a \rangle_{x} \colonequals \mathfrak{O}^{a}_{\alpha}.
\eeno
Patching together the above formulas, we get
\beno
|\nabla_{x} \pa^{\alpha}a|^{2}_{L^{2}_{x}} + \epsilon \frac{\mathrm{d}}{\mathrm{d}t} \mathcal{I}^{a}_{\alpha}(f) = \mathfrak{U}^{a}_{\alpha} + \mathfrak{O}^{a}_{\alpha}.
\eeno
By the Cauchy-Schwartz inequality, one has
\beno  |\mathfrak{U}^{a}_{\alpha}| \leq \eta T^{1/2} |\nabla_{x} \partial^{\alpha} b|^{2}_{L^{2}_{x}} + \frac{C}{\eta T^{1/2}}
|\epsilon\partial_{t}\pa^{\alpha} a|^{2}_{L^{2}_{x}}  + \frac{CT^{1/2}}{\eta} |\nabla_{x}\partial^{\alpha}\mathcal{U}|^{2}_{L^{2}_{x}},
\\ |\mathfrak{O}^{a}_{\alpha}| \leq \eta T^{1/2} |\nabla_{x} \partial^{\alpha} a|^{2}_{L^{2}_{x}} + \frac{C}{\eta T^{1/2}} (|\partial^{\alpha}\mathcal{V}|^{2}_{L^{2}_{x}} + |\partial^{\alpha}\mathcal{W}|^{2}_{L^{2}_{x}} + |\partial^{\alpha}\mathcal{X}|^{2}_{L^{2}_{x}}).
\eeno
Taking sum over $|\alpha| \leq N-1$,
by Lemma \ref{estimate-for-ptabc} and \eqref{simple-estimate-of-pt-a}, we get
\ben \label{dissipation-a-H-N-minus-1} |\nabla_{x}a|^{2}_{H^{N-1}_{x}} +\epsilon \frac{\mathrm{d}}{\mathrm{d}t} \sum_{|\alpha|\leq N-1} \mathcal{I}^{a}_{\alpha}(f)
\leq \mathfrak{H}_{\eta, C}. \een

For brevity,
let us define the temporal  energy functional  $\mathcal{I}_{N}(f)$ as in \cite{duan2008cauchy} by
\ben \label{interactive-INf} \mathcal{I}_{N}(f) \colonequals    \sum_{|\alpha|\leq N-1} (\mathcal{I}^{a}_{\alpha}(f) +
\sum_{j=1}^{3}\mathcal{I}^{b_{j}}_{\alpha}(f)+\mathcal{I}^{c}_{\alpha}(f)), \een
where $\mathcal{I}^{a}_{\alpha}(f), \mathcal{I}^{b_{j}}_{\alpha}(f), \mathcal{I}^{c}_{\alpha}(f)$ are defined in \eqref{a-temperal-energy}, \eqref{bj-temperal-energy}, \eqref{c-temperal-energy}.
Patching together \eqref{dissipation-a-H-N-minus-1}, \eqref{dissipation-b-H-N-minus-1} and \eqref{dissipation-c-H-N-minus-1}, recalling  \eqref{interactive-INf} and \eqref{definition-H-eta}, we have
\beno  \epsilon\frac{\mathrm{d}}{\mathrm{d}t}\mathcal{I}_{N}(f) + T^{1/2} |\nabla_{x}(a,b,c)|^{2}_{H^{N-1}_{x}}
\leq \mathfrak{H}_{\eta, C}.\eeno
Taking $\eta=\frac{1}{2}$, we arrive at \eqref{solution-property-part2}. We emphasize that $C$ is a universal constant. In particular, it is even independent of the integer $N$ by the above derivation.

Recall $f=f_{1}+f_{2}$, then $\partial^{\alpha}f= (\partial^{\alpha}f)_{1}+(\partial^{\alpha}f)_{2} = \partial^{\alpha}f_{1}+\partial^{\alpha}f_{2}$. That is, the projection operator $\tilde{\mathbb{P}}_{\lambda}$ and the spatial derivative $\partial^{\alpha}$ is commutative. So we have
\beno
|\partial^{\alpha}f|_{L^{2}_{x}L^{2}}^{2} = |\partial^{\alpha}f_{1}|_{L^{2}_{x}L^{2}}^{2} + |\partial^{\alpha}f_{2}|_{L^{2}_{x}L^{2}}^{2}.
\eeno
Recalling the definition of $\mathcal{I}_{N}(f)$ in \eqref{interactive-INf}, back to  \eqref{a-temperal-energy}, \eqref{bj-temperal-energy} and \eqref{c-temperal-energy},
by the C-S inequality inequality, recalling \eqref{explicit-defintion-of-abc}, using \eqref{l22-fives-functions} and
\eqref{partial-alpha-mathcal-U}, we get
\beno
|\mathcal{I}_{N}(f)| \lesssim |(a,b,c)|_{H^{N}_{x}}^{2} + |\mathcal{U}|_{H^{N}_{x}}^{2}  \lesssim e^{\lambda}|f_{1}|_{H^{N}_{x}L^{2}}^{2} + e^{\lambda}|f_{2}|_{H^{N}_{x}L^{2}}^{2} = e^{\lambda}|f|_{H^{N}_{x}L^{2}}^{2}.
\eeno
That is, \eqref{temporal-bounded-by-norm} holds for some universal constant $C$.
\end{proof}

We derive the following {\it a priori} estimate for equation \eqref{lBE}.

\begin{prop}\label{essential-estimate-of-micro-macro} Let $N \geq 1; \lambda, T>0, T_{*} >0$. Let $f \in L^{\infty}([0,T_{*}]; H^{N}_{x}L^{2})$ be a solution to \eqref{lBE},
then for any $0 \leq t \leq T_{*}$, it holds that
\ben \label{essential-micro-macro-result-final}
&&\frac{\mathrm{d}}{\mathrm{d}t}\Xi^{\lambda,T}_{N}(f)+ \f{1}{2} \mathcal{D}_{N}(f) + \mathrm{C}_{1}(\lambda,T) \f{1}{\epsilon^{2}}
\|f_{2}\|^{2}_{H^{N}_{x}L^{2}_{1/2}}
\\ \nonumber &\lesssim&    \mathrm{C}_{2}(\lambda,T)\sum_{|\alpha| \leq N}
|(\pa^{\alpha}g, \pa^{\alpha}f)| + \mathrm{C}_{3}(\lambda,T) \epsilon^{2} \mathcal{Q}_{\lambda, N-1}(g), \een
where
\ben \label{defintion-of-Xi-lambda-T}
\Xi^{\lambda,T}_{N}(f)\colonequals  \epsilon T^{-1/2} e^{-\lambda}(1-e^{-\lambda})^{-\f{1}{2}} C_0 \mathcal{I}_{N}(f) + K(\lambda,T) \|f\|^{2}_{H^{N}_{x}L^{2}},
\\
\label{defintion-of-K-lambda-T}
K(\lambda,T) \colonequals  C_{*} e^{2\lambda} (1-e^{-\lambda})^{-\f{21}{2}} \max\{T, T^{-2}\},
\\ \label{defintion-of-C1-lambda-T}
\mathrm{C}_{1}(\lambda,T)  \colonequals e^{\lambda} (1-e^{-\lambda})^{-\f{21}{2}} \max\{T, T^{-2}\} T^{2},
\\ \label{defintion-of-C2-lambda-T}
\mathrm{C}_{2}(\lambda,T)  \colonequals  e^{2\lambda} (1-e^{-\lambda})^{-\f{21}{2}} \max\{T, T^{-2}\},
\\ \label{defintion-of-C3-lambda-T}
\mathrm{C}_{3}(\lambda,T)  \colonequals  e^{-\lambda}(1-e^{-\lambda})^{-\f{1}{2}} T^{-1},
\een
for some large universal constant $C_{*}>0$. Moreover, it holds that
\ben \label{equivalence-between-energy}
\frac{1}{2} K(\lambda,T) \|f\|^{2}_{H^{N}_{x}L^{2}} \leq \Xi^{\lambda,T}_{N}(f) \leq \frac{3}{2} K(\lambda,T) \|f\|^{2}_{H^{N}_{x}L^{2}}.
\een
\end{prop}

\begin{proof}
Note that $\Xi^{\lambda,T}_{N}(f)$ is a combination of $\mathcal{I}_{N}(f)$ and $\|f\|^{2}_{H^{N}_{x}L^{2}}$. We already have $\mathcal{I}_{N}(f)$ from Lemma \ref{estimate-for-highorder-abc}. That is, the solution $f$ verifies \eqref{solution-property-part2}. By \eqref{solution-property-part2} and recalling the constant $\tilde{C}_{1,\lambda,T}$ from \eqref{constant-of-C-1},
\ben  \label{inter-energy-estimate}
T^{-1/2}\epsilon\frac{\mathrm{d}}{\mathrm{d}t}\mathcal{I}_{N}(f) + \f{1}{2}|\nabla_{x}(a,b,c)|^{2}_{H^{N-1}_{x}}
\leq
 C T^{3} e^{\lambda} (1-e^{-\lambda})^{-8} \epsilon^{-2}
 \|f_{2}\|_{H^{N}_{x}L^{2}_{1/2}}^{2} + C T^{-1}\epsilon^{2}\mathcal{Q}_{\lambda, N-1}(g),\een
where $C$ is a universal constant.

Applying $\partial^{\alpha}$ to equation \eqref{lBE}, taking inner product with $\partial^{\alpha}f$, recalling $f_{2} = (\mathbb{I}-\tilde{\mathbb{P}}_{\lambda})f$ and $(\partial^{\alpha}f)_{2} = \partial^{\alpha}f_{2}$, by the lower bound in \eqref{lower-and-upper-bound},
 taking sum over $|\alpha|\leq N$, we have
\ben \label{solution-property-part-g}\f{1}{2}\frac{\mathrm{d}}{\mathrm{d}t}\|f\|^{2}_{H^{N}_{x}L^{2}} + C_{1,\lambda,T}  \f{1}{\epsilon^{2}} \|f_{2}\|^{2}_{H^{N}_{x}L^{2}_{1/2}} \leq \sum_{|\alpha| \leq N}
|(\pa^{\alpha}g, \pa^{\alpha}f)|.\een
We now use the term $\|f_{2}\|^{2}_{H^{N}_{x}L^{2}_{1/2}}$ in
\eqref{solution-property-part-g} to control the right-hand of  \eqref{inter-energy-estimate}. The combination $\eqref{solution-property-part-g} \times 2K(\lambda,T) + \eqref{inter-energy-estimate} \times C_0 e^{-\lambda}(1-e^{-\lambda})^{-\f{1}{2}}$ gives
\beno % \label{essential-micro-macro-result}
&&\frac{\mathrm{d}}{\mathrm{d}t}(\epsilon T^{-1/2} C_0 e^{-\lambda}(1-e^{-\lambda})^{-\f{1}{2}} \mathcal{I}_{N}(f) + K(\lambda,T) \|f\|^{2}_{H^{N}_{x}L^{2}})
\\&&+ (\f{1}{2} C_0 e^{-\lambda}(1-e^{-\lambda})^{-\f{1}{2}}|\nabla_{x}(a,b,c)|^{2}_{H^{N-1}_{x}}+ K(\lambda,T) C_{1,\lambda,T}   \f{1}{\epsilon^{2}}
\|f_{2}\|^{2}_{H^{N}_{x}L^{2}_{1/2}})
\\ % \nonumber
&\leq& C e^{-\lambda}(1-e^{-\lambda})^{-\f{1}{2}} T^{-1}\epsilon^{2} \mathcal{Q}_{\lambda, N-1}(g) + 2K(\lambda,T)\sum_{|\alpha| \leq N}
|(\pa^{\alpha}g, \pa^{\alpha}f)|,
\eeno
by taking $C_*$ large enough in
\eqref{defintion-of-K-lambda-T}
such that
\beno
K(\lambda,T) C_{1,\lambda,T} \geq \f{1}{4} C T^{3} e^{\lambda} (1-e^{-\lambda})^{-8}.
\eeno
Thanks to \eqref{temporal-bounded-by-norm}, we can also ask $C_{*}$ large enough in \eqref{defintion-of-K-lambda-T} such that \eqref{equivalence-between-energy} holds.
Recalling from  \eqref{defintion-of-dissipation-after-scaling} that $\mathcal{D}_{N}(f) = C_0 e^{-\lambda}(1-e^{-\lambda})^{-\f{1}{2}}|\nabla_{x}(a,b,c)|^{2}_{H^{N-1}_{x}}+\|f_{2}\|^{2}_{H^{N}_{x}L^{2}_{1/2}}$, by taking $C_*$ large enough in
\eqref{defintion-of-K-lambda-T}
such that
$
K(\lambda,T) C_{1,\lambda,T} \geq \f{1}{8},
$
we get the dissipation $\f{1}{2}\mathcal{D}_{N}(f)$
and finish the proof.
\end{proof}

\subsection{A priori estimate of the quantum Boltzmann equation.}
In this subsection, we derive the following uniform-in-$\epsilon$ { \it a priori} estimate of the Cauchy problem \eqref{T-half-scaled-more-consistent}.
\begin{thm}\label{a-priori-estimate-LBE} Let $\lambda, T > 0$. Recall the constants $\tilde{C}_{*}(\lambda,T), K(\lambda,T), \mathrm{C}_{1}(\lambda,T), \mathrm{C}_{2}(\lambda,T), \mathrm{C}_{3}(\lambda,T)$ defined in \eqref{tilde-C-star}, \eqref{defintion-of-K-lambda-T}, \eqref{defintion-of-C1-lambda-T}, \eqref{defintion-of-C2-lambda-T}, \eqref{defintion-of-C3-lambda-T}.
There exists a constant  $\delta_{2}>0$ which is independent of $\epsilon,\lambda$ and $T$, such that if
a solution $f$ to the Cauchy problem \eqref{T-half-scaled-more-consistent}  satisfies
\ben \label{smalness-on-h2l2}
\sup_{t}\|f(t)\|_{H^{2}_{x}L^{2}}^{2} \leq \delta_{2} \tilde{C}_{*}(\lambda,T),
\een
then $f$ verifies  for any $N \geq 2$,
\ben \label{uniform-estimate-propagation}
\sup_{t}\|f(t)\|_{H^{N}_{x}L^{2}}^{2}+ \f{1}{K(\lambda,T)}  \int_{0}^{\infty} \mathcal{D}_{N}(f) \mathrm{d}\tau + \f{\mathrm{C}_{1}(\lambda,T)}{K(\lambda,T)}  \f{1}{\epsilon^{2}} \int_{0}^{\infty} \|f_{2}\|^{2}_{H^{N}_{x}L^{2}_{1/2}} \mathrm{d}\tau \leq  P_{N}(f_{0})\|f(t)\|_{H^{N}_{x}L^{2}}^{2}.
\een
where
\ben \label{definition-of-pn}
P_{2}(f_{0}) \equiv 6, \quad P_{N}(f_{0}) \colonequals  12 \exp\left(Q_{3}(\lambda,T,N,f_{0}) P_{N-1}(f_{0})\|f_{0}\|_{H^{N-1}_{x}L^{2}}^{2} \right)   \text{ for } N \geq 3,
\\ \label{definition of -Q-3-f0}
Q_{3}(\lambda,T,N,f_{0}) \colonequals  2(Q_{1}(\lambda,T,N)  +
Q_{2}(\lambda,T,N) P_{N-1}(f_{0})\|f_{0}\|_{H^{N-1}_{x}L^{2}}^{2}).
\een
where the constants $Q_{1}(\lambda,T,N)$ and $Q_{2}(\lambda,T,N)$ are defined in \eqref{defintion-Q-1-constant} and \eqref{defintion-Q-2-constant} respectively.
\end{thm}
%Let $N \geq 2$. There is a constant  $\delta_{N}>0$ which is independent of $\lambda$ and $T$, such that if
%a solution $f$ to the Cauchy problem \eqref{T-half-scaled-more-consistent}  satisfies
%\ben  \label{smalness-on-hnl2}
%\sup_{t}\|f(t)\|_{H^{N}_{x}L^{2}}^{2} \leq \delta_{N} \tilde{C}_{*}(\lambda,T),\een
%then $f$ verifies
%\ben \label{uniform-estimate-propagation-case-N}
%\sup_{t}\|f(t)\|_{H^{N}_{x}L^{2}}^{2}+ \f{1}{K(\lambda,T)}  \int_{0}^{\infty} \mathcal{D}_{N}(f) \mathrm{d}\tau + \f{\mathrm{C}_{1}(\lambda,T)}{K(\lambda,T)}  \f{1}{\epsilon^{2}} \int_{0}^{\infty} \|f_{2}\|^{2}_{H^{N}_{x}L^{2}_{1/2}} \mathrm{d}\tau \leq 12 \|f_{0}\|_{H^{N}_{x}L^{2}}^{2}.
%\een
\begin{proof} Observe that $f$ solves \eqref{lBE} with $g= \f{1}{\epsilon}\tilde{\Gamma}_{2}^{\lambda,T}(f,f) + \tilde{\Gamma}_{3}^{\lambda,T}(f,f,f)$.
By  Proposition  \ref{essential-estimate-of-micro-macro}, recalling that $\mathcal{Q}_{\lambda, n}(g)$ is bounded by some linear combination of the quantities in \eqref{non-linear-functional-general-source},
we have
\ben \label{essential-micro-macro-result-final-ap}
\frac{\mathrm{d}}{\mathrm{d}t}\Xi^{\lambda,T}_{N}(f)+ \f{1}{2} \mathcal{D}_{N}(f) + \mathrm{C}_{1}(\lambda,T) \f{1}{\epsilon^{2}}\|f_{2}\|^{2}_{H^{N}_{x}L^{2}_{1/2}}
\lesssim   \mathrm{C}_{2}(\lambda,T) \mathcal{I}_{1}(f) + \mathrm{C}_{3}(\lambda,T) \epsilon^{2} \mathcal{I}_{2}(f). \een
Here the two terms on the right-hand of \eqref{essential-micro-macro-result-final-ap} are
\beno
\mathcal{I}_{1}(f) &\colonequals & \sum_{|\alpha| \leq N}
|(\f{1}{\epsilon}\pa^{\alpha}\tilde{\Gamma}_{2}^{\lambda,T}(f,f) + \pa^{\alpha} \Gamma_{3}^{\lambda,T}(f,f,f), \pa^{\alpha}f)|,
\\
\mathcal{I}_{2}(f) &\colonequals &  \sum_{|\alpha| \leq N-1}  |\langle  \f{1}{\epsilon}\pa^{\alpha}\tilde{\Gamma}_{2}^{\lambda,T}(f,f) + \pa^{\alpha} \Gamma_{3}^{\lambda,T}(f,f,f), \psi \rangle|_{L^{2}_{x}}^{2}.
\eeno
By \eqref{Gamma-2-energy-estimate} and \eqref{Gamma-3-energy-estimate}, we have
\beno
\mathcal{I}_{1}(f) &\lesssim& \f{1}{\epsilon} C_{2,\lambda,T}\left( \|f\|_{H^{2}_{x}L^{2}} \mathcal{D}^{\f{1}{2}}_{N}(f) + \mathrm{1}_{N \geq 3}C_{N}
\|f\|_{H^{N}_{x}L^{2}}\mathcal{D}^{\f{1}{2}}_{N-1}(f) \right) \|f_{2}\|_{H^{N}_{x}L^{2}_{1/2}}
\\&&+ C_{3,\lambda,T}\left( \|f\|_{H^{2}_{x}L^{2}}^{2}\mathcal{D}_{N}(f) + \mathrm{1}_{N \geq 3}C_{N}
\|f\|_{H^{N-1}_{x}L^{2}}\|f\|_{H^{N}_{x}L^{2}}\mathcal{D}^{\f{1}{2}}_{N-1}(f)\mathcal{D}^{\f{1}{2}}_{N}(f) \right).
\eeno
By \eqref{Gamma-2-energy-estimate-with-13} and \eqref{Gamma-3-energy-estimate-with-13}, we have
\beno
\mathcal{I}_{2}(f) &\lesssim& \f{1}{\epsilon^{2}} C_{2,\lambda,T}^{2}\left( \|f\|_{H^{2}_{x}L^{2}}^{2}\mathcal{D}_{N}(f) + \mathrm{1}_{N \geq 3}C_{N}
\|f\|_{H^{N}_{x}L^{2}}^{2}\mathcal{D}_{N-1}(f) \right)
\\&&+ C_{3,\lambda,T}^{2}\left( \|f\|_{H^{2}_{x}L^{2}}^{4}\mathcal{D}_{N}(f) + \mathrm{1}_{N \geq 3}C_{N}
\|f\|_{H^{N-1}_{x}L^{2}}^{2}\|f\|_{H^{N}_{x}L^{2}}^{2}\mathcal{D}_{N-1}(f) \right).
\eeno
Plugging the previous two inequalities into \eqref{essential-micro-macro-result-final-ap}, we have
\ben \label{essential-micro-macro-all-bounds} &&\frac{\mathrm{d}}{\mathrm{d}t}\Xi^{\lambda,T}_{N}(f)+ \f{1}{2} \mathcal{D}_{N}(f) + \mathrm{C}_{1}(\lambda,T) \f{1}{\epsilon^{2}}\|f_{2}\|^{2}_{H^{N}_{x}L^{2}_{1/2}}
\\&\leq& \label{most-dangerous-term} C \mathrm{C}_{2}(\lambda,T) \f{1}{\epsilon} C_{2,\lambda,T}\left( \|f\|_{H^{2}_{x}L^{2}} \mathcal{D}^{\f{1}{2}}_{N}(f) + \mathrm{1}_{N \geq 3}C_{N}
\|f\|_{H^{N}_{x}L^{2}}\mathcal{D}^{\f{1}{2}}_{N-1}(f) \right) \|f_{2}\|_{H^{N}_{x}L^{2}_{1/2}}
\\ &&+ C \mathrm{C}_{2}(\lambda,T)
C_{3,\lambda,T}\left( \|f\|_{H^{2}_{x}L^{2}}^{2}\mathcal{D}_{N}(f) + \mathrm{1}_{N \geq 3}C_{N}
\|f\|_{H^{N-1}_{x}L^{2}}\|f\|_{H^{N}_{x}L^{2}}\mathcal{D}^{\f{1}{2}}_{N-1}(f)\mathcal{D}^{\f{1}{2}}_{N}(f) \right)
\\ &&+
C \mathrm{C}_{3}(\lambda,T) C_{2,\lambda,T}^{2}\left( \|f\|_{H^{2}_{x}L^{2}}^{2}\mathcal{D}_{N}(f) + \mathrm{1}_{N \geq 3}C_{N}
\|f\|_{H^{N}_{x}L^{2}}^{2}\mathcal{D}_{N-1}(f) \right)
\\ &&+
C \mathrm{C}_{3}(\lambda,T) \epsilon^{2} C_{3,\lambda,T}^{2}\left( \|f\|_{H^{2}_{x}L^{2}}^{4}\mathcal{D}_{N}(f) + \mathrm{1}_{N \geq 3}C_{N}
\|f\|_{H^{N-1}_{x}L^{2}}^{2}\|f\|_{H^{N}_{x}L^{2}}^{2}\mathcal{D}_{N-1}(f) \right). \een
Here $C \geq 2$ is a universal constant.

Recalling \eqref{defintion-of-C-2-lambda-T} for the definition of $C_{2,\lambda,T}$. Recalling \eqref{defintion-of-C-3-lambda-T} for the definition of $C_{3,\lambda,T}$.
Note that the line \eqref{most-dangerous-term} is bounded by
\beno
C \mathrm{C}_{2}(\lambda,T) \f{1}{\epsilon} C_{2,\lambda,T} \|f\|_{H^{2}_{x}L^{2}} \mathcal{D}^{\f{1}{2}}_{N}(f) \|f_{2}\|_{H^{N}_{x}L^{2}_{1/2}}
 \leq \f{1}{16}\mathcal{D}_{N}(f) + 4 C^{2} \mathrm{C}_{2}^{2}(\lambda,T)  C_{2,\lambda,T}^{2} \|f\|_{H^{2}_{x}L^{2}}^{2} \f{1}{\epsilon^{2}} \|f_{2}\|_{H^{N}_{x}L^{2}_{1/2}}^{2}.
\eeno
Under the following conditions
\ben \label{condition-1}
4 C^{2} \mathrm{C}_{2}^{2}(\lambda,T)  C_{2,\lambda,T}^{2} \|f\|_{H^{2}_{x}L^{2}}^{2} \leq \f{1}{2}\mathrm{C}_{1}(\lambda,T),
\\ \label{condition-2}
C \mathrm{C}_{2}(\lambda,T)
C_{3,\lambda,T} \|f\|_{H^{2}_{x}L^{2}}^{2} \leq \f{1}{16},
\\ \label{condition-3}
C \mathrm{C}_{3}(\lambda,T) C_{2,\lambda,T}^{2} \|f\|_{H^{2}_{x}L^{2}}^{2} \leq \f{1}{16},
\\ \label{condition-4}
C \mathrm{C}_{3}(\lambda,T) C_{3,\lambda,T}^{2} \|f\|_{H^{2}_{x}L^{2}}^{4} \leq \f{1}{16},
\een
 we get
\ben \label{case-N-2}
N=2 && \frac{\mathrm{d}}{\mathrm{d}t}\Xi^{\lambda,T}_{N}(f)+ \f{1}{4} \mathcal{D}_{N}(f) + \f{1}{2} \mathrm{C}_{1}(\lambda,T) \f{1}{\epsilon^{2}}\|f_{2}\|^{2}_{H^{N}_{x}L^{2}_{1/2}}  \leq 0,
 \\ \label{case-N-geq-3}
N \geq 3 &&\frac{\mathrm{d}}{\mathrm{d}t}\Xi^{\lambda,T}_{N}(f)+ \f{1}{4} \mathcal{D}_{N}(f) + \f{1}{2} \mathrm{C}_{1}(\lambda,T) \f{1}{\epsilon^{2}}\|f_{2}\|^{2}_{H^{N}_{x}L^{2}_{1/2}}
\\ \nonumber &\leq&C_{N} \mathrm{C}_{2}(\lambda,T) \f{1}{\epsilon} C_{2,\lambda,T}
\|f\|_{H^{N}_{x}L^{2}}\mathcal{D}^{\f{1}{2}}_{N-1}(f) \|f_{2}\|_{H^{N}_{x}L^{2}_{1/2}}
\\  \nonumber &&+ C_{N} \mathrm{C}_{2}(\lambda,T)
C_{3,\lambda,T}
\|f\|_{H^{N-1}_{x}L^{2}}\|f\|_{H^{N}_{x}L^{2}}\mathcal{D}^{\f{1}{2}}_{N-1}(f)\mathcal{D}^{\f{1}{2}}_{N}(f)
\\  \nonumber &&+
C_{N} \mathrm{C}_{3}(\lambda,T) C_{2,\lambda,T}^{2}
\|f\|_{H^{N}_{x}L^{2}}^{2}\mathcal{D}_{N-1}(f)
\\  \nonumber &&+
C_{N} \mathrm{C}_{3}(\lambda,T) \epsilon^{2} C_{3,\lambda,T}^{2}
\|f\|_{H^{N-1}_{x}L^{2}}^{2}\|f\|_{H^{N}_{x}L^{2}}^{2}\mathcal{D}_{N-1}(f).
\een
Here $C_{N}$ is a large constant that depends only on 
$N$ and could change from line to line.

Note that \eqref{condition-1} is stronger than \eqref{condition-2}, \eqref{condition-3}, \eqref{condition-4} since $C \geq 2$. By
taking $\delta_{2} = \f{1}{8 C^{2}}$, \eqref{smalness-on-h2l2} yields
the condition \eqref{condition-1}.

We now prove \eqref{uniform-estimate-propagation} for $N \geq 2$.
If $N=2$, integrating \eqref{case-N-2} w.r.t. time, we get
\beno
\Xi^{\lambda,T}_{N}(f(t))+ \f{1}{4} \int_{0}^{t} \mathcal{D}_{N}(f) \mathrm{d}\tau + \f{1}{2} \mathrm{C}_{1}(\lambda,T) \f{1}{\epsilon^{2}} \int_{0}^{t} \|f_{2}\|^{2}_{H^{N}_{x}L^{2}_{1/2}} \mathrm{d}\tau \leq \Xi^{\lambda,T}_{N}(f_{0})
\eeno
Recalling \eqref{equivalence-between-energy}, we get
\ben \label{result-n-equals-2}
\sup_{t}\|f(t)\|_{H^{2}_{x}L^{2}}^{2}+ \f{1}{K(\lambda,T)} \int_{0}^{\infty} \mathcal{D}_{2}(f) \mathrm{d}\tau + \f{\mathrm{C}_{1}(\lambda,T)}{K(\lambda,T)}  \f{1}{\epsilon^{2}} \int_{0}^{t} \|f_{2}\|^{2}_{H^{2}_{x}L^{2}_{1/2}} \mathrm{d}\tau \leq 6 \|f_{0}\|_{H^{2}_{x}L^{2}}^{2}.
\een

Now we use mathematical induction to establish \eqref{uniform-estimate-propagation} for $N \geq 2$. Suppose \eqref{uniform-estimate-propagation} is valid for $N-1 \geq 2$, we now prove it is also valid for $N > 3$.
By \eqref{case-N-geq-3}, using $ab \leq \eta a^{2} + \f{1}{4\eta} b^{2}$, we get
\beno
&&\frac{\mathrm{d}}{\mathrm{d}t}\Xi^{\lambda,T}_{N}(f)+ \f{1}{8} \mathcal{D}_{N}(f) + \f{1}{4} \mathrm{C}_{1}(\lambda,T) \f{1}{\epsilon^{2}}\|f_{2}\|^{2}_{H^{N}_{x}L^{2}_{1/2}}
\\ \nonumber &\leq& \mathrm{C}_{1}^{-1}(\lambda,T) C_{N}^{2} \mathrm{C}_{2}^{2}(\lambda,T)  C_{2,\lambda,T}^{2}
\|f\|_{H^{N}_{x}L^{2}}^{2}\mathcal{D}_{N-1}(f)
\\  \nonumber &&+ 2C_{N}^{2} \mathrm{C}_{2}^{2}(\lambda,T)
C^{2}_{3,\lambda,T}
\|f\|_{H^{N-1}_{x}L^{2}}^{2}\|f\|_{H^{N}_{x}L^{2}}^{2}\mathcal{D}_{N-1}(f)
\\  \nonumber &&+
C_{N} \mathrm{C}_{3}(\lambda,T) C_{2,\lambda,T}^{2}
\|f\|_{H^{N}_{x}L^{2}}^{2}\mathcal{D}_{N-1}(f)
\\  \nonumber &&+
C_{N} \mathrm{C}_{3}(\lambda,T) \epsilon^{2} C_{3,\lambda,T}^{2}
\|f\|_{H^{N-1}_{x}L^{2}}^{2}\|f\|_{H^{N}_{x}L^{2}}^{2}\mathcal{D}_{N-1}(f)
\\ \nonumber &\leq& Q_{1}(\lambda,T) \|f\|_{H^{N}_{x}L^{2}}^{2}\mathcal{D}_{N-1}(f) +
Q_{2}(\lambda,T) \|f\|_{H^{N-1}_{x}L^{2}}^{2} \|f\|_{H^{N}_{x}L^{2}}^{2}\mathcal{D}_{N-1}(f),
\eeno
where we define for simplicity
\ben \label{defintion-Q-1-constant}
Q_{1}(\lambda,T,N) \colonequals  \mathrm{C}_{1}^{-1}(\lambda,T) C_{N}^{2} \mathrm{C}_{2}^{2}(\lambda,T)  C_{2,\lambda,T}^{2} + C_{N} \mathrm{C}_{3}(\lambda,T) C_{2,\lambda,T}^{2},
\\ \label{defintion-Q-2-constant}
Q_{2}(\lambda,T,N) \colonequals  2C_{N}^{2} \mathrm{C}_{2}^{2}(\lambda,T) C^{2}_{3,\lambda,T} + C_{N} \mathrm{C}_{3}(\lambda,T) C_{3,\lambda,T}^{2}.
\een

By the induction assumption,
\beno
\sup_{t}\|f(t)\|_{H^{N-1}_{x}L^{2}}^{2}+ \f{1}{K(\lambda,T)}  \int_{0}^{\infty} \mathcal{D}_{N-1}(f) \mathrm{d}\tau + \f{\mathrm{C}_{1}(\lambda,T)}{K(\lambda,T)}  \f{1}{\epsilon^{2}} \int_{0}^{t} \|f_{2}\|^{2}_{H^{N-1}_{x}L^{2}_{1/2}} \mathrm{d}\tau \leq  P_{N-1}(f_{0})\|f_{0}\|_{H^{N-1}_{x}L^{2}}^{2}.
\eeno
Then the energy inequality becomes
\beno
&&\frac{\mathrm{d}}{\mathrm{d}t}\Xi^{\lambda,T}_{N}(f)+ \f{1}{8} \mathcal{D}_{N}(f) + \f{1}{4} \mathrm{C}_{1}(\lambda,T)
\f{1}{\epsilon^{2}}\|f_{2}\|^{2}_{H^{N}_{x}L^{2}_{1/2}}
\\&\leq& (Q_{1}(\lambda,T,N)  +
Q_{2}(\lambda,T,N) P_{N-1}(f_{0})\|f_{0}\|_{H^{N-1}_{x}L^{2}}^{2}) \|f\|_{H^{N}_{x}L^{2}}^{2} \mathcal{D}_{N-1}(f)
\\&\leq& Q_{3}(\lambda,T,N,f_{0}) \Xi^{\lambda,T}_{N}(f) \f{\mathcal{D}_{N-1}(f)}{K(\lambda,T)},
\eeno
where we recall  \eqref{equivalence-between-energy} and \eqref{definition of -Q-3-f0}.
 Then by Gr\"{o}nwall's inequality and the induction assumption, we get
\beno
&&\Xi^{\lambda,T}_{N}(f(t))+ \int_{0}^{t} \big(\f{1}{8} \mathcal{D}_{N}(f) + \f{1}{4} \mathrm{C}_{1}(\lambda,T)
\f{1}{\epsilon^{2}}\|f_{2}\|^{2}_{H^{N}_{x}L^{2}_{1/2}}\big)  \mathrm{d}\tau
\\ \nonumber &\leq& \exp\left(Q_{3}(\lambda,T,N,f_{0}) \int_{0}^{\infty} \f{\mathcal{D}_{N-1}(f(t))}{K(\lambda,T)} \mathrm{d}t \right) \Xi^{\lambda,T}_{N}(f_{0}) \leq \exp\left(Q_{3}(\lambda,T,N,f_{0}) P_{N-1}(f_{0})\|f_{0}\|_{H^{N-1}_{x}L^{2}}^{2} \right) \Xi^{\lambda,T}_{N}(f_{0}).
\eeno
Recalling \eqref{equivalence-between-energy}, we get
\beno
&&\sup_{t}\|f(t)\|_{H^{N}_{x}L^{2}}^{2} + \f{1}{K(\lambda,T)} \int_{0}^{t} \big( \mathcal{D}_{N}(f) +  \mathrm{C}_{1}(\lambda,T)
\f{1}{\epsilon^{2}}\|f_{2}\|^{2}_{H^{N}_{x}L^{2}_{1/2}}\big)  \mathrm{d}\tau
\\&\leq& 12 \exp\left(Q_{3}(\lambda,T,N,f_{0}) P_{N-1}(f_{0})\|f_{0}\|_{H^{N-1}_{x}L^{2}}^{2} \right) \|f_{0}\|_{H^{N}_{x}L^{2}}^{2}.
\eeno
By recalling \eqref{definition-of-pn},
we finish the proof of \eqref{uniform-estimate-propagation}.
\end{proof}

%
%Suppose we have the smallness assumption  \eqref{smalness-on-hnl2} on the norm $H^{N}_{x}L^{2}$, we get directly from \eqref{essential-micro-macro-all-bounds} that
%\ben \label{small-on-hnl2-result} \frac{\mathrm{d}}{\mathrm{d}t}\Xi^{\lambda,T}_{N}(f)+ \f{1}{4} \mathcal{D}_{N}(f) + \f{1}{2}\mathrm{C}_{1}(\lambda,T) \f{1}{\epsilon^{2}}\|f_{2}\|^{2}_{H^{N}_{x}L^{2}_{1/2}}
%\leq 0. \een
%From this together with \eqref{equivalence-between-energy}, as the same as
%\eqref{result-n-equals-2}, we get \eqref{uniform-estimate-propagation-case-N} and finish the proof.

\section{Hydrodynamic limits} \label{hydrodynamic-limit}

This whole section is devoted to prove Theorem \ref{thm-hydrodynamical-limit}. We first derive some basic formulas involving Bose-Einstein distribution.
Recall from \eqref{definition-of-M-lambda} that $M_{\lambda}$ is a radial function.
Sometimes, we also write for $r \geq 0$,
\ben
M_{\lambda}(r)\colonequals  \f{1}{ \exp(\f{r^{2}}{2}+\lambda) -1}.
\een
It is easy to check that
\ben \label{derivatives-of-M}
\f{\mathrm{d} M_{\lambda}}{\mathrm{d} r} = - r M_{\lambda}(1+M_{\lambda}), \quad
\f{\mathrm{d} [M_{\lambda}(1+M_{\lambda})]}{\mathrm{d} r} = - r M_{\lambda}(1+M_{\lambda})(1+2M_{\lambda}).
\een
Then by \eqref{derivatives-of-M}, polar coordinates, integration by parts formula, it is easy to derive
\ben \label{2-to-1}
\int_{\mathbb{R}^{3}} f(|v|)|v|^{k} M_{\lambda}(v) (1+M_{\lambda}(v)) \mathrm{d}v = (k+1) \int_{\mathbb{R}^{3}} f(|v|)|v|^{k-2} M_{\lambda}(v) \mathrm{d}v,
\\ \label{3-to-2}
\int_{\mathbb{R}^{3}} f(|v|)|v|^{k} M_{\lambda}(v) (1+M_{\lambda}(v)) (1+2M_{\lambda}(v)) \mathrm{d}v = (k+1) \int_{\mathbb{R}^{3}} f(|v|)|v|^{k-2} M_{\lambda}(v) (1+M_{\lambda}(v)) \mathrm{d}v.
\een

Recall \eqref{moment-k-of-N-square} for $m_{k}$.
Let
\ben \label{defintion-of-KA}
K_{A} = \f{m_{4}}{2m_{2}}, \quad
K_{\lambda} = K_{A} -1 = \f{m_{4}}{2m_{2}} - 1, \quad C_{A}  = \f{m_{4}}{4 m_{2}^{2}} (m_{4} m_{0} - m_{2}^{2}).
\een
Then it is elementary to check
\ben \label{why-the-constant-KA}
\int_{\mathbb{R}^{3}} (\f{|v|^{2}}{2} - K_{A}) |v|^{2}M_{\lambda}(v) (1+M_{\lambda}(v)) \mathrm{d}v = \f{1}{2} m_{4} - K_{A} m_{2} = 0,
\\ \label{why-the-constant-CA}
\f{1}{3} \int_{\mathbb{R}^{3}} (\f{|v|^{2}}{2} - K_{A})^{2} |v|^{2} M_{\lambda}(v) (1+M_{\lambda}(v)) (1+ 2M_{\lambda}(v)) \mathrm{d}v = C_{A}.
\een

Let us introduce
\ben \label{definition-of-vector-and-matrix-bos}
A (v)= N_{\lambda}(v) (\f{|v|^{2}}{2} - K_{A}) v, \quad B(v) =  N_{\lambda}(v)(v \otimes v - \f{|v|^{2}}{3}I_{3}).
\een
Note that $A$ is  a vector, and $B$ is a symmetric matrix. With the help of \eqref{why-the-constant-KA}, $A_{i}, B_{ij} \in (\ker\tilde{\mathcal{L}}^{\lambda,T})^{\perp}$. More precisely, for any $f \in \ker\tilde{\mathcal{L}}^{\lambda,T}$, it holds that
\ben \label{perpendicular-to-kernel-space}
\langle A_{i}, f \rangle = \langle B_{ij}, f \rangle = 0.
\een

By rotational invariance, it is standard to derive
\begin{thm} \label{form-of-inverse-lambda-T}
There exist unique radial functions $\alpha_{\lambda,T}(|v|), \beta_{\lambda,T}(|v|)$  such that
\ben \label{inverse-function-definition-lambda-T}
\tilde{\mathcal{L}}^{\lambda,T} \left(\alpha_{\lambda,T} A\right) = A, \quad \tilde{\mathcal{L}}^{\lambda,T}\left(\beta_{\lambda,T}B\right) = B.
\een
For later reference, let us denote
\ben \label{inverse-function-form-lambda-T}
\hat{A}(v) := \alpha_{\lambda,T} (|v|) A(v), \quad \hat{B}(v) := \beta_{\lambda,T}(|v|) B(v).
\een
\end{thm}
See Theorem \ref{form-of-inverse} for why $\hat{A}, \hat{B}$ must be of the form in \eqref{inverse-function-form-lambda-T} in order to satisfy $\tilde{\mathcal{L}}^{\lambda,T} ( \hat{A}) = A, \tilde{\mathcal{L}}^{\lambda,T} ( \hat{B}) = B$.

\begin{proof}[Proof of Theorem \ref{thm-hydrodynamical-limit}]
By the well-posedness theory in Theorem \ref{global-well-posedness}, the family $\{f_{\epsilon}\}_{0<\epsilon<1}$ verifies
\ben \label{L-infty-t-HNx-L2-cb}
M_{\infty} \colonequals  \sup_{0<\epsilon<1} \sup_{t \geq 0}\|f^{\epsilon}(t)\|_{H^{N}_{x}L^{2}} \leq C(M_{0}),
\\ \label{L-2-t-HNx-L2-cb}
M_{2} \colonequals   \sup_{0<\epsilon<1} \left( \f{1}{\epsilon^{2}}
\int_{0}^{\infty} (\|f^{\epsilon}(t) - \mathbb{P} f^{\epsilon}(t)\|_{H^{N}_{x}L^{2}}^{2}
 \mathrm{d}t \right)\leq C(M_{0}).
\een
Here $C(M_0)$ is a constant depending the constant $M_0$ given in  
\eqref{small-condition}. Note that for brevity, we drop the dependence on $\lambda,T,N$. 

By \eqref{L-infty-t-HNx-L2-cb}, there is a subsequence of $\{f^{\epsilon}\}$ still denoting it by $\{f^{\epsilon}\}$ such that
 \ben \label{weak-in-L-infty}
 f^{\epsilon} \to f^{0} \text{ as }
\epsilon \to 0, \text{ weakly-* in } L^{\infty}(\mathbb{R}_{+}; H^{N}_{x}L^{2}),
 \een
for some $f^{0} \in L^{\infty}(\mathbb{R}_{+}; H^{N}_{x}L^{2})$.
By \eqref{L-2-t-HNx-L2-cb}, $\{f^{\epsilon} - \mathbb{P} f^{\epsilon} \}_{0<\epsilon<1}$ converges to $0$ in $L^{2}_{t}H^{N}_{x}L^{2}$.
\ben \label{f2-strong-in-L-2}
f^{\epsilon} - \mathbb{P} f^{\epsilon} \to 0 \text{ as } \epsilon \to 0,
\text{ strongly in } L^{2}(\mathbb{R}_{+}; H^{N}_{x}L^{2}).
\een
By \eqref{weak-in-L-infty} and \eqref{f2-strong-in-L-2}, we have $f^{0} \in \ker \tilde{\mathcal{L}}^{\lambda,T}$. Then there exists $(\rho, u, \theta) \in L^{\infty}(\mathbb{R}_{+}; H^{N}_{x})$ such that
\ben \label{limit-form}
f^{0}(t,x,v) = (\rho(t,x) + u(t,x) \cdot v + \theta(t,x) (\f{|v|^{2}}{2} - K_{\lambda}) )N_{\lambda}(v).
\een
Note that we already get \eqref{weakly-star-convergence-to-some-special-form}.

We now prove that $(\rho, u, \theta) \in C(\mathbb{R}_{+}; H^{N-1}_{x})$ satisfies the system \eqref{NSF-problem} and the moment convergence \eqref{momoent-1} and \eqref{momoent-2}. This is done by looking at $(\rho^{\epsilon}, u^{\epsilon}, \theta^{\epsilon})$, the macroscopic components of $f^{\epsilon}$. Recalling \eqref{linear-combination-of-basis}, \eqref{defintion-of-l-i} and \eqref{explicit-defintion-of-abc}, we have
\beno
\mathbb{P}_{\lambda} f^{\epsilon} = \left(\rho^{\epsilon} + u^{\epsilon} \cdot v + \theta^{\epsilon} (\f{|v|^{2}}{2} - K_{\lambda}) \right) N_{\lambda},
\eeno
where
\beno
\rho^{\epsilon} = \f{K_{A}}{C_{A}} \langle f^{\epsilon}, N_{\lambda}\rangle + \f{K_{A}}{C_{A}}  (\f{K_{\lambda}m_{0}}{m_{2}} - \f{1}{2})\langle f^{\epsilon}, |v|^{2} N_{\lambda}\rangle, \quad
u^{\epsilon} = \f{3}{m_{2}}\langle f^{\epsilon}, v N_{\lambda} \rangle,
\\
\theta^{\epsilon} = \f{K_{A}}{C_{A}} \f{m_{0}}{m_{2}} \langle f^{\epsilon}, |v|^{2} N_{\lambda}\rangle
-\f{K_{A}}{C_{A}} \langle f^{\epsilon}, N_{\lambda}\rangle.
\eeno
By \eqref{weak-in-L-infty}, \eqref{f2-strong-in-L-2} and \eqref{limit-form}, we have
\ben \label{rho-u-theta-convergence}
(\rho^{\epsilon}, u^{\epsilon}, \theta^{\epsilon}) \to (\rho, u, \theta) \text{ as }
\epsilon \to 0, \text{ weakly-* in } L^{\infty}(\mathbb{R}_{+}; H^{N}_{x}).
\een

Taking inner products
between
\eqref{T-half-scaled-more-consistent} and the following functions
\beno
\f{K_{A}}{C_{A}} N_{\lambda} (1 + (\f{K_{\lambda}m_{0}}{m_{2}} - \f{1}{2})|v|^{2}),
\quad \f{3}{m_{2}}  N_{\lambda} v,
\quad
\f{K_{A}}{C_{A}} \f{m_{0}}{m_{2}}N_{\lambda}(|v|^{2}-\f{m_{2}}{m_{0}}),
\eeno
we get
\ben \label{eq-rho-epsilon-more-consistent}
\partial_{t} \rho^{\epsilon} +  \f{T^{\f{1}{2}}}{\epsilon} \langle  v \cdot \nabla_{x} f^{\epsilon}, \f{K_{A}}{C_{A}} N_{\lambda} (1 + (\f{K_{\lambda}m_{0}}{m_{2}} - \f{1}{2})|v|^{2})\rangle &=& 0,
\\ \label{eq-u-epsilon-more-consistent}
\partial_{t} u^{\epsilon} + \f{T^{\f{1}{2}}}{\epsilon} \langle  v \cdot \nabla_{x} f^{\epsilon}, \f{3}{m_{2}}  N_{\lambda} v\rangle &=& 0,
\\ \label{eq-theta-epsilon-more-consistent}
\partial_{t} \theta^{\epsilon} + \f{T^{\f{1}{2}}}{\epsilon} \langle  v \cdot \nabla_{x} f^{\epsilon}, \f{K_{A}}{C_{A}} \f{m_{0}}{m_{2}}N_{\lambda}(|v|^{2}-\f{m_{2}}{m_{0}})\rangle &=& 0.
\een

Now we need to compute the inner products in \eqref{eq-rho-epsilon-more-consistent}, \eqref{eq-u-epsilon-more-consistent} and
\eqref{eq-theta-epsilon-more-consistent}.
Recall the macro-micro decomposition $f^{\epsilon} = f^{\epsilon}_{1} + f^{\epsilon}_{2}$ where
$
f^{\epsilon}_{1} = \mathbb{P}_{\lambda} f^{\epsilon}, \quad f^{\epsilon}_{2} = f^{\epsilon} - \mathbb{P}_{\lambda} f^{\epsilon}
$.
Using
\beno
\int v_{i} \phi(|v|) \mathrm{d}v = 0, \int v_{i}v_{j} \phi(|v|) \mathrm{d}v = \delta_{ij} \int \f{|v|^{2}}{3} \phi(|v|) \mathrm{d}v, \int v_{i}v_{j}v_{k} \phi(|v|) \mathrm{d}v = 0,
\eeno
it is easy to see
\beno
\langle v \cdot \nabla_{x} f^{\epsilon}_{1}, N_{\lambda}\rangle =
\langle v \cdot \nabla_{x}  (\rho^{\epsilon} + u^{\epsilon} \cdot v + \theta^{\epsilon} (\f{|v|^{2}}{2}-K_{\lambda}))N_{\lambda}, N_{\lambda}\rangle = \langle v \cdot \nabla_{x}  (u^{\epsilon} \cdot v), N_{\lambda}^{2}\rangle = \f{m_{2}}{3}\nabla_{x} \cdot u^{\epsilon},
\\
\langle v \cdot \nabla_{x} f^{\epsilon}_{1}, |v|^{2}N_{\lambda}\rangle =
\langle v \cdot \nabla_{x}  (\rho^{\epsilon} + u^{\epsilon} \cdot v + \theta^{\epsilon} (\f{|v|^{2}}{2}-K_{\lambda}))N_{\lambda}, |v|^{2}N_{\lambda}\rangle = \langle v \cdot \nabla_{x}  (u^{\epsilon} \cdot v), |v|^{2}N_{\lambda}^{2}\rangle = \f{m_{4}}{3}\nabla_{x} \cdot u^{\epsilon}.
\eeno
As a result, by recalling \eqref{defintion-of-KA},
we have
\ben \label{f-1-order-0}
\langle v \cdot \nabla_{x} f^{\epsilon}_{1}, \f{K_{A}}{C_{A}} N_{\lambda} (1 + (\f{K_{\lambda}m_{0}}{m_{2}} - \f{1}{2})|v|^{2})\rangle
= \f{K_{A}}{C_{A}} (\f{m_{2}}{3} + (\f{K_{\lambda}m_{0}}{m_{2}} - \f{1}{2}) \f{m_{4}}{3}) \nabla_{x} \cdot u^{\epsilon} = \f{2}{3} K_{\lambda} \nabla_{x} \cdot u^{\epsilon},
\\ \label{f-1-order-2}
\langle v \cdot \nabla_{x} f^{\epsilon}_{1}, \f{K_{A}}{C_{A}} \f{m_{0}}{m_{2}}N_{\lambda}(|v|^{2}-\f{m_{2}}{m_{0}}) \rangle
= \f{K_{A}}{C_{A}} (\f{m_{2}}{3} + (\f{K_{\lambda}m_{0}}{m_{2}} - \f{1}{2}) \f{m_{4}}{3}) \nabla_{x} \cdot u^{\epsilon} = \f{2}{3} \nabla_{x} \cdot u^{\epsilon}.
\een

Note that
\beno
\langle v \cdot \nabla_{x} f^{\epsilon}_{1}, v N_{\lambda}\rangle &=&
\langle v \cdot \nabla_{x}  (\rho^{\epsilon} + u^{\epsilon} \cdot v + \theta^{\epsilon} (\f{|v|^{2}}{2}-K_{\lambda}))N_{\lambda}, v N_{\lambda}\rangle
\\
&=& \langle v \cdot \nabla_{x}  (\rho^{\epsilon} + \theta^{\epsilon} (\f{|v|^{2}}{2}-K_{\lambda}))N_{\lambda}, vN_{\lambda}\rangle
\\
&=& \langle v \cdot \nabla_{x}  (\rho^{\epsilon} + \theta^{\epsilon}), v N_{\lambda}^{2} \rangle
= \f{m_{2}}{3} \nabla_{x}  (\rho^{\epsilon} + \theta^{\epsilon}),
\eeno
where we recall $K_{\lambda} = K_{A} -1$ and use \eqref{perpendicular-to-kernel-space}.
 As a result
\ben \label{f-1-order-1}
\langle  v \cdot \nabla_{x} f^{\epsilon}_{1}, \f{3}{m_{2}}  N_{\lambda} v\rangle = \nabla_{x}  (\rho^{\epsilon} + \theta^{\epsilon}).
\een

Since $f^{\epsilon}_{2} \in (\ker \tilde{\mathcal{L}}^{\lambda,T})^{\perp}$, we have
\ben \label{f-2-order-0}
\langle v \cdot \nabla_{x} f^{\epsilon}_{2}, N_{\lambda}\rangle = \nabla_{x} \cdot \langle  f^{\epsilon}_{2}, v N_{\lambda}\rangle =0,
\\ \nonumber
\langle v \cdot \nabla_{x} f^{\epsilon}_{2}, v N_{\lambda}\rangle = \nabla_{x} \cdot \langle  f^{\epsilon}_{2}, v \otimes v N_{\lambda}\rangle =  \nabla_{x} \cdot \langle  f^{\epsilon}_{2}, (v \otimes v - \f{1}{3}|v^{2}| I_{3}) N_{\lambda}\rangle,
\\ \nonumber
\langle v \cdot \nabla_{x} f^{\epsilon}_{2}, \f{|v|^{2}}{2} N_{\lambda}\rangle = \nabla_{x} \cdot \langle  f^{\epsilon}_{2}, \f{|v|^{2}}{2} v N_{\lambda}\rangle = \nabla_{x} \cdot \langle  f^{\epsilon}_{2}, (\f{|v|^{2}}{2}-K_{A})v N_{\lambda}\rangle.
\een
Recalling \eqref{definition-of-vector-and-matrix-bos} and \eqref{inverse-function-definition-lambda-T}, since $\tilde{\mathcal{L}}^{\lambda,T}$ is self-adjoint and
 $f^{\epsilon}_{2} \in (\ker \tilde{\mathcal{L}}^{\lambda,T})^{\perp}$, we have
\ben \label{f-2-order-1}
\langle v \cdot \nabla_{x} f^{\epsilon}_{2}, v N_{\lambda}\rangle  =  \nabla_{x} \cdot \langle  f^{\epsilon}_{2}, B \rangle =  \nabla_{x} \cdot \langle  f^{\epsilon}_{2}, \tilde{\mathcal{L}}^{\lambda,T}\hat{B} \rangle = \nabla_{x} \cdot \langle  \tilde{\mathcal{L}}^{\lambda,T}f^{\epsilon}_{2}, \hat{B} \rangle = \nabla_{x} \cdot \langle  \tilde{\mathcal{L}}^{\lambda,T}f^{\epsilon}, \hat{B} \rangle,
\\ \label{f-2-order-2}
\langle v \cdot \nabla_{x} f^{\epsilon}_{2}, \f{|v|^{2}}{2} N_{\lambda}\rangle = \nabla_{x} \cdot \langle  f^{\epsilon}_{2}, A \rangle = \nabla_{x} \cdot \langle  f^{\epsilon}_{2}, \tilde{\mathcal{L}}^{\lambda,T}\hat{A} \rangle = \nabla_{x} \cdot \langle \tilde{\mathcal{L}}^{\lambda,T} f^{\epsilon}_{2}, \hat{A} \rangle = \nabla_{x} \cdot \langle \tilde{\mathcal{L}}^{\lambda,T} f^{\epsilon}, \hat{A} \rangle.
\een

Plugging \eqref{f-1-order-0}-\eqref{f-2-order-2} into
\eqref{eq-rho-epsilon-more-consistent}-\eqref{eq-theta-epsilon-more-consistent}, we get
\ben \label{eq-rho-epsilon-2}
\partial_{t} \rho^{\epsilon} + \f{T^{\f{1}{2}}}{\epsilon} \f{2}{3} K_{\lambda} \nabla_{x} \cdot u^{\epsilon} + \f{T^{\f{1}{2}}}{\epsilon}
\f{K_{A}}{C_{A}}\f{2K_{\lambda}m_{0}-m_{2}}{m_{2}}\nabla_{x} \cdot \langle \hat{A}, \tilde{\mathcal{L}}^{\lambda,T} f^{\epsilon} \rangle&=& 0,
\\ \label{eq-u-epsilon-2}
\partial_{t} u^{\epsilon} +  \f{T^{\f{1}{2}}}{\epsilon} \nabla_{x}  (\rho^{\epsilon} + \theta^{\epsilon}) + \f{T^{\f{1}{2}}}{\epsilon} \f{3}{m_{2}}
\nabla_{x} \cdot \langle  \hat{B}, \tilde{\mathcal{L}}^{\lambda,T}f^{\epsilon} \rangle&=& 0,
\\ \label{eq-theta-epsilon-2}
\partial_{t} \theta^{\epsilon} + \f{T^{\f{1}{2}}}{\epsilon}  \f{2}{3} \nabla_{x} \cdot u^{\epsilon} + \f{T^{\f{1}{2}}}{\epsilon} \f{2K_{A}}{C_{A}} \f{m_{0}}{m_{2}} \nabla_{x} \cdot \langle \hat{A}, \tilde{\mathcal{L}}^{\lambda,T} f^{\epsilon} \rangle &=& 0.
\een

By \eqref{eq-rho-epsilon-2} and recalling \eqref{f-2-order-2}, by
\eqref{L-infty-t-HNx-L2-cb} and \eqref{L-2-t-HNx-L2-cb},
in the distributional sense,
\ben \label{nabla-u-epsilon-to-0}
T^{\f{1}{2}} \f{2}{3} K_{\lambda} \nabla_{x} \cdot u^{\epsilon} = - \epsilon \partial_{t} \rho^{\epsilon} - T^{\f{1}{2}} \f{K_{A}}{C_{A}}\f{2K_{\lambda}m_{0}-m_{2}}{m_{2}}\nabla_{x} \cdot \langle f^{\epsilon}_{2}, A \rangle
\to 0.
\een
By \eqref{nabla-u-epsilon-to-0} and recalling \eqref{rho-u-theta-convergence}, we get
\ben \label{incompressible-condition}
\nabla_{x} \cdot u = 0.
\een

By \eqref{eq-u-epsilon-2} and recalling \eqref{f-2-order-1}, by
\eqref{L-infty-t-HNx-L2-cb} and \eqref{L-2-t-HNx-L2-cb}, in the distributional sense,
\ben \label{nabla-rho-plus-theta-epsilon-to-0}
T^{\f{1}{2}} \nabla_{x}  (\rho^{\epsilon} + \theta^{\epsilon}) = - \epsilon \partial_{t} u - T^{\f{1}{2}} \f{3}{m_{2}}
\nabla_{x} \cdot \langle  f^{\epsilon}_{2}, B \rangle \to 0.
\een
By \eqref{nabla-rho-plus-theta-epsilon-to-0} and recalling \eqref{rho-u-theta-convergence}, we get
\ben \label{Boussinesq-relation}
\nabla_{x}  (\rho + \theta) = 0, \quad \Rightarrow \quad \rho + \theta = 0.
\een

{\it Convergence of $\f{K_{\lambda}\theta^{\epsilon} - \rho^{\epsilon}}{K_{\lambda}+1}$ and regularity of $\theta$.}
 Making a suitable combination of \eqref{eq-rho-epsilon-2} and \eqref{eq-theta-epsilon-2}, we get
\ben \label{combination-rho-theta}
\partial_{t} (\f{K_{\lambda}\theta^{\epsilon} - \rho^{\epsilon}}{K_{\lambda}+1}) + \f{T^{\f{1}{2}}}{\epsilon}
\f{K_{A}}{C_{A}(K_{\lambda}+1)} \nabla_{x} \cdot \langle \hat{A}, \tilde{\mathcal{L}}^{\lambda,T} f^{\epsilon} \rangle&=& 0
\een
Then by
\eqref{L-infty-t-HNx-L2-cb} and \eqref{L-2-t-HNx-L2-cb},
\ben
\partial_{t} (\f{K_{\lambda}\theta^{\epsilon} - \rho^{\epsilon}}{K_{\lambda}+1}) \in L^{2}(\mathbb{R}_{+}; H^{N-1}_{x}), \quad \f{K_{\lambda}\theta^{\epsilon} - \rho^{\epsilon}}{K_{\lambda}+1} \in L^{\infty}(\mathbb{R}_{+}; H^{N}_{x}).
\een
Then by Aubin-Lions-Simon Theorem, \eqref{rho-u-theta-convergence} and \eqref{Boussinesq-relation}, we get $\rho, \theta \in L^{\infty}(\mathbb{R}_{+}; H^{N}_{x}) \cap C(\mathbb{R}_{+}; H^{N-1}_{x})$ and
\ben \label{strong-convergence-moment-2}
\f{K_{\lambda}\theta^{\epsilon} - \rho^{\epsilon}}{K_{\lambda}+1} \to \theta \text{ strongly in } C(\mathbb{R}_{+}; H^{N-1}_{x}).
\een

Recall that the
Leray projection $\mathcal{P}$ is defined by
$
\mathcal{P}u = \mathcal{I} - \nabla \Delta \nabla \cdot
$.
Note that
\ben \label{Leray-projection-property-1}
\mathcal{P}u = u  \quad \Leftrightarrow \quad \nabla \cdot u = 0.
\\ \label{Leray-projection-property-2}
\mathcal{P} \nabla  \phi = 0.
\een

{\it{Convergence of $\mathcal{P}u^{\epsilon}$ and regularity of $u$.}} Applying $\mathcal{P}$ to \eqref{eq-u-epsilon-2}, we get
\ben \label{Leray-projection-u-epsilon-2}
\partial_{t} \mathcal{P}u^{\epsilon} + \f{T^{\f{1}{2}}}{\epsilon} \f{3}{m_{2}}
 \mathcal{P} \nabla_{x} \cdot \langle  \hat{B}, \tilde{\mathcal{L}}^{\lambda,T}f^{\epsilon} \rangle= 0.
\een
Then by
\eqref{L-infty-t-HNx-L2-cb} and \eqref{L-2-t-HNx-L2-cb},
\ben
\partial_{t} \mathcal{P}u^{\epsilon} \in L^{2}(\mathbb{R}_{+}; H^{N-1}_{x}), \quad \mathcal{P}u^{\epsilon} \in L^{\infty}(\mathbb{R}_{+}; H^{N}_{x}).
\een
Then by Aubin-Lions-Simon Theorem, \eqref{rho-u-theta-convergence} and \eqref{incompressible-condition}, we get $u \in L^{\infty}(\mathbb{R}_{+}; H^{N}_{x}) \cap C(\mathbb{R}_{+}; H^{N-1}_{x})$ and
\ben \label{Pu-epsilon-to-u}
\mathcal{P}u^{\epsilon} \to u = \mathcal{P}u\text{ strongly in } C(\mathbb{R}_{+}; H^{N-1}_{x}).
\een 
By \eqref{rho-u-theta-convergence} and \eqref{Pu-epsilon-to-u}, in the distributional sense, $\mathcal{P}^{\perp}u^{\epsilon} \to 0$. By \eqref{L-infty-t-HNx-L2-cb}, $\|\mathcal{P}^{\perp}u^{\epsilon}\|_{L^{\infty}(\mathbb{R}_{+}; H^{N}_{x})} \lesssim C(M_{0})$. As a result,
\ben \label{p-perp-u-epsilon-2-0}
\mathcal{P}^{\perp}u^{\epsilon} \to 0 \text{ weakly-* in } L^{\infty}(\mathbb{R}_{+}; H^{N}_{x}).
\een

Note that we already proved $(\rho, u, \theta) \in C(\mathbb{R}_{+}; H^{N-1}_{x})$ and the moment convergence \eqref{momoent-1}(given by \eqref{Pu-epsilon-to-u}) and \eqref{momoent-2}(given by \eqref{strong-convergence-moment-2}). We also proved $(\rho, u, \theta)$ satisfies ${\eqref{NSF-problem}}_{1}$. Now it remains to show $(u, \theta)$ satisfies ${\eqref{NSF-problem}}_{2}, {\eqref{NSF-problem}}_{3}$ with the initial conditions in ${\eqref{NSF-problem}}_{4}$. To this end, we will come back to \eqref{T-half-scaled-more-consistent} to evaluate $\langle  \hat{B}, \tilde{\mathcal{L}}^{\lambda,T}f^{\epsilon} \rangle$ and $\langle  \hat{A}, \tilde{\mathcal{L}}^{\lambda,T}f^{\epsilon} \rangle$. By \eqref{T-half-scaled-more-consistent}, we have
\ben \label{equation-for-Lf}
\f{1}{\epsilon} \tilde{\mathcal{L}}^{\lambda,T}f^{\epsilon} = -  T^{1/2} v \cdot \nabla_{x} f^{\epsilon}_{1} + \tilde{\Gamma}_{2}^{\lambda,T}(f^{\epsilon}_{1},f^{\epsilon}_{1}) + R_{\epsilon},
\een
where
\beno
R_{\epsilon} \colonequals  -\epsilon \partial _t f^{\epsilon} + \epsilon \tilde{\Gamma}_{3}^{\lambda,T}(f^{\epsilon}, f^{\epsilon}, f^{\epsilon})
- T^{1/2} v \cdot \nabla_{x} f^{\epsilon}_{2} + \tilde{\Gamma}_{2}^{\lambda,T}(f^{\epsilon}_{1},f^{\epsilon}_{2}) + \tilde{\Gamma}_{2}^{\lambda,T}(f^{\epsilon}_{2},f^{\epsilon}_{1}) + \tilde{\Gamma}_{2}^{\lambda,T}(f^{\epsilon}_{2},f^{\epsilon}_{2}).
\eeno
Note that $R_{\epsilon}$ is bounded by $O(\epsilon)+O(f^{\epsilon}_{2})$. Then it is easy to check
\ben \label{error-term}
R_{\epsilon} \to 0 \text{ weakly-* in } L^{\infty}(\mathbb{R}_{+}; H^{N-1}_{x}L^{2}).
\een

 To derive more equations,
we need to consider the following four quantities
\beno
\langle  \hat{A}, v \cdot \nabla_{x} f^{\epsilon}_{1} \rangle, \quad \langle  \hat{B}, v \cdot \nabla_{x} f^{\epsilon}_{1} \rangle, \quad \langle  \hat{A}, \tilde{\Gamma}_{2}^{\lambda,T}(f^{\epsilon}_{1},f^{\epsilon}_{1}) \rangle, \quad \langle  \hat{B}, \tilde{\Gamma}_{2}^{\lambda,T}(f^{\epsilon}_{1},f^{\epsilon}_{1}) \rangle,
\eeno
We give three lemmas to derive them.
\begin{lem}
Let
$
g = (\rho + u \cdot v + \theta (\f{|v|^{2}}{2} - K_{\lambda})) N_{\lambda}
$.
Then
\ben \label{A-with-macro}
 \langle \hat{A}, v \cdot \nabla_{x} g \rangle = \kappa_{1,\lambda, T} \nabla_{x} \theta + \kappa_{2,\lambda, T} \nabla_{x} (\rho+\theta),
\\ \label{B-with-macro}
 \langle \hat{B}, v \cdot \nabla_{x} g \rangle = \nu_{\lambda, T} \mathcal{T}(u).
\een
where
\beno
\kappa_{1,\lambda, T} \colonequals  \f{1}{3} \int \alpha_{\lambda, T}(|v|) (\f{|v|^{2}}{2} - K_{A})^{2}|v|^{2} N_{\lambda}^{2} \mathrm{d}v
,\quad
\kappa_{2,\lambda, T} \colonequals  \f{1}{3} \int \alpha_{\lambda, T}(|v|) (\f{|v|^{2}}{2} - K_{A}) |v|^{2} N_{\lambda}^{2} \mathrm{d}v
,
\\
\nu_{\lambda, T}\colonequals  \f{1}{15} \int \beta_{\lambda, T}(|v|) |v|^{4} N_{\lambda}^{2} \mathrm{d}v,  \quad  \mathcal{T}(u) \colonequals   \nabla_{x} u + (\nabla_{x} u)^{\mathrm{T}} - \f{2}{3} (\nabla_{x} \cdot u) I_{3}.
\eeno
\end{lem}
\begin{proof}Note that
\beno
\langle \hat{A}, v \cdot \nabla_{x} g \rangle = \nabla_{x} \cdot \langle \hat{A} \otimes v,  g \rangle
= \nabla_{x} \cdot \langle \hat{A} \otimes v,  (\rho + \theta + \theta (\f{|v|^{2}}{2} - K_{A})) N_{\lambda} \rangle
= \nabla_{x} \cdot \langle \hat{A} \otimes A, \theta \rangle + \nabla_{x} \cdot \langle \hat{A} \otimes N_{\lambda}v,  \rho + \theta \rangle.
\eeno
We then get \eqref{A-with-macro} by observing
\beno
 \int \hat{A} \otimes A \mathrm{d}v = \int \alpha_{\lambda, T}(|v|) A \otimes A \mathrm{d}v =
 \f{1}{3} \int \alpha_{\lambda, T}(|v|) (\f{|v|^{2}}{2} - K_{A})^{2}|v|^{2} N_{\lambda}^{2} \mathrm{d}v I_{3} = \kappa_{1,\lambda, T} I_{3},
 \\
 \int \hat{A} \otimes N_{\lambda}v \mathrm{d}v = \int \alpha_{\lambda, T}(|v|) A \otimes N_{\lambda}v \mathrm{d}v =
 \f{1}{3} \int \alpha_{\lambda, T}(|v|) (\f{|v|^{2}}{2} - K_{A}) |v|^{2} N_{\lambda}^{2} \mathrm{d}v I_{3} = \kappa_{1,\lambda, T} I_{3}.
\eeno

Note that
\beno
\langle \hat{B}, v \cdot \nabla_{x} g \rangle = \langle \hat{B}, N_{\lambda} (v \cdot \nabla_{x})(u \cdot v) \rangle
= \langle \hat{B}, N_{\lambda} \sum_{i,j} \partial_{i}u_{j} v_{i}v_{j} \rangle.
\eeno
For the $(1,1)$ element, we have
\beno
&&\langle \beta_{\lambda, T}(|v|) N_{\lambda}^{2}, (v_{1}^{2} - \f{|v|^{2}}{3})  \sum_{i,j} \partial_{i}u_{j} v_{i}v_{j} \rangle  \\&=&
\langle \beta_{\lambda, T}(|v|)N_{\lambda}^{2}, (v_{1}^{2} - \f{|v|^{2}}{3})  v_{1}^{2} \rangle \partial_{1}u_{1}
+ \langle \beta_{\lambda, T}(|v|)N_{\lambda}^{2}, (v_{1}^{2} - \f{|v|^{2}}{3})  v_{2}^{2} \rangle \partial_{2}u_{2}
+ \langle \beta_{\lambda, T}(|v|)N_{\lambda}^{2}, (v_{1}^{2} - \f{|v|^{2}}{3})  v_{3}^{2} \rangle \partial_{3}u_{3}
\\&=& \langle \beta_{\lambda, T}(|v|)N_{\lambda}^{2}, (v_{1}^{4} - v_{1}^{2}v_{2}^{2})   \rangle \partial_{1}u_{1} + \langle \beta_{\lambda, T}(|v|)N_{\lambda}^{2}, (v_{1}^{2} - \f{|v|^{2}}{3})  v_{2}^{2} \rangle \nabla_{x} \cdot u
\colonequals  \lambda_{11} \partial_{1}u_{1} + \lambda_{*} \nabla_{x} \cdot u.
\eeno
For the $(1,2)$ element, we have
\beno
\langle \beta_{\lambda, T}(|v|) N_{\lambda}^{2}, v_{1}v_{2}  \sum_{i,j} \partial_{i}u_{j} v_{i}v_{j} \rangle  =
\langle \beta_{\lambda, T}(|v|)N_{\lambda}^{2}, v_{1}^{2}v_{2}^{2}  \rangle (\partial_{1}u_{2} + \partial_{2}u_{1})
\colonequals  \lambda_{12} (\partial_{1}u_{2} + \partial_{2}u_{1}).
\eeno
Observe that for a general radial function $f(v)=f(|v|)$,
\beno
\int f(|v|) |v|^{n}  \mathrm{d}v = (n+1)\int f(|v|) v_{1}^{n}  \mathrm{d}v.
\eeno
Since $\int \beta(|v|) N_{\lambda}^{2} |v|^{4} \mathrm{d}v = 5 \int \beta(|v|) N_{\lambda}^{2} v_{1}^{4} \mathrm{d}v$, then $\lambda_{11} = 2 \lambda_{12}, \lambda_{*} = -\f{2}{3} \lambda_{12}$.
As a result, we get
\beno
\langle \hat{B}, v \cdot \nabla_{x} g \rangle = \lambda_{12} (\nabla_{x} u + (\nabla_{x} u)^{\mathrm{T}}) + \lambda_{*} (\nabla_{x} \cdot u) I_{3} = \lambda_{12} (\nabla_{x} u + (\nabla_{x} u)^{\mathrm{T}} - \f{2}{3} (\nabla_{x} \cdot u) I_{3}).
\eeno
It is easy to check that $\lambda_{12} = \langle \beta_{\lambda, T}(|v|)N_{\lambda}^{2}, v_{1}^{2}v_{2}^{2}  \rangle = \f{1}{15}\langle \beta_{\lambda, T}(|v|)N_{\lambda}^{2}, |v|^{4}  \rangle = \nu_{\lambda, T}.$
\end{proof}

\begin{lem} \label{Gamma-2-gg-to-L}
For $g \in \ker \tilde{\mathcal{L}}^{\lambda,T}$, we have
\ben \label{Gamma-g-g-to-L-g-square-bosons}
\tilde{\Gamma}_{2}^{\lambda,T} (g,g) = \f{1}{2} \tilde{\mathcal{L}}^{\lambda,T}((1+2M_{\lambda})N_{\lambda}^{-1}g^{2}).
\een
\end{lem}
\begin{proof}
Recall from \eqref{scaled-linearized-operator-more-consistent} and
\eqref{definition-Gamma-2-epsilon-rho-T} the definitions of $\tilde{\mathcal{L}}^{\lambda,T}$ and  $\Gamma_{2}^{\lambda,T}$. Let $K_{\lambda}\colonequals  N_{\lambda} (N_{\lambda})_{*} N_{\lambda}^{\prime} (N_{\lambda})^{\prime}_{*}$, then
\beno
(\tilde{\mathcal{L}}^{\lambda,T}f)(v)  \colonequals   N_{\lambda}^{-1} \int  B_{T} K_{\lambda} \mathrm{S}(N_{\lambda}^{-1}f)
\mathrm{d}\sigma \mathrm{d}v_{*},
\\
\tilde{\Gamma}_{2}^{\lambda,T}(g,h) \colonequals    N_{\lambda}^{-1}\int
B_{T} K_{\lambda} \Theta_{2} (g,h)
\mathrm{d}\sigma \mathrm{d}v_{*},
\eeno
where
\beno   %\label{definition-theta-2}
\Theta_{2}(g,h) &\colonequals   & \mathrm{D} \big((N_{\lambda}^{-1}g)^{\prime}_{*}(N_{\lambda}^{-1}h)^{\prime}\big)
\\ %\label{line-1-another} 
&&+
\mathrm{D}\big((N_{\lambda}^{-1}g)^{\prime}_{*}(N_{\lambda}^{-1}h)^{\prime}(M_{\lambda}^{\prime} +  (M_{\lambda})^{\prime}_{*})\big)
\\ %\label{line-2-another} 
&&+ \mathrm{D}\big( (N_{\lambda}^{-1}g)^{\prime}_{*}(N_{\lambda}^{-1}h)( M_{\lambda} -  (M_{\lambda})^{\prime}_{*}) \big)
\\ %\label{line-3-another} 
&&+ \mathrm{D}((N_{\lambda}^{-1}g)^{\prime}_{*}(N_{\lambda}^{-1}g)_{*}(M_{\lambda})_{*}) + \mathrm{D}((N_{\lambda}^{-1}h)^{\prime}(N_{\lambda}^{-1}h)M_{\lambda}).
\eeno
Note that
\beno  % \label{definition-theta-2-another}
\Theta_{2}(g,g) &\colonequals   & \mathrm{D} \big((N_{\lambda}^{-1}g)^{\prime}_{*}(N_{\lambda}^{-1}g)^{\prime}\big)
\\ %\label{line-1-another-2}
 &&+
\mathrm{D}\left((M_{\lambda}N_{\lambda}^{-1}g)^{\prime}_{*}\big((N_{\lambda}^{-1}g)^{\prime} - (N_{\lambda}^{-1}g) -(N_{\lambda}^{-1}g)_{*}\big)\right)
\\ %\label{line-2-another-2} 
&&+ \mathrm{D}\left((M_{\lambda}N_{\lambda}^{-1}g)^{\prime}\big((N_{\lambda}^{-1}g)^{\prime}_{*} - (N_{\lambda}^{-1}g) -(N_{\lambda}^{-1}g)_{*}\big)\right).
\eeno
For $g \in \ker \tilde{\mathcal{L}}^{\lambda,T}$,
setting $g = N_{\lambda} h$, then $h$ is collision invariant. Then
\beno %  \label{definition-theta-2-another-h-h}
\Theta_{2}(g,g) &=& \mathrm{D} \big(h^{\prime}_{*}h^{\prime}\big)
\\ % \label{line-1-another-2-h-h}
 &&+
\mathrm{D}\left((M_{\lambda}h)^{\prime}_{*}\big(h^{\prime} - h -h_{*}\big)\right)
\\ % \label{line-2-another-2-h-h} 
&&+ \mathrm{D}\left((M_{\lambda}h)^{\prime}\big(h^{\prime}_{*} - h -h_{*}\big)\right)
\\ &=& -\mathrm{D}(h_{*}h) - \mathrm{D}\left((M_{\lambda}h)_{*}\big(h- h^{\prime} -h^{\prime}_{*}\big)\right)
- \mathrm{D}\left((M_{\lambda}h)\big(h_{*} - h^{\prime} -h^{\prime}_{*}\big)\right).
\eeno
Note that
\beno
\f{1}{2} \mathrm{S}(N_{\lambda}^{-1}(1+2M_{\lambda})N_{\lambda}^{-1}g^{2}) &=&
\f{1}{2} \mathrm{S}((1+2M_{\lambda})h^{2}) = \f{1}{2} \mathrm{S}(h^{2}) + \mathrm{S}(M_{\lambda}h^{2})
\\ &=& \f{1}{2} \mathrm{D}(h^{2} + h^{2}_{*}) + \mathrm{D}(M_{\lambda}h^{2} + (M_{\lambda}h^{2})_{*}).
\eeno
Recalling $h$ is collision invariant, we get
\beno
\f{1}{2} \mathrm{S}(N_{\lambda}^{-1}(1+2M_{\lambda})N_{\lambda}^{-1}g^{2}) - \Theta_{2}(g,g)
= \f{1}{2} \mathrm{D}((h + h_{*})^{2}) + \mathrm{D}\left(M_{\lambda}h \mathrm{S}h + (M_{\lambda}h)_{*}\mathrm{S}h\right) =0,
\eeno
which gives \eqref{Gamma-g-g-to-L-g-square-bosons}.
\end{proof}

\begin{lem} Let
$
g = (\rho + u \cdot v + \theta (\f{|v|^{2}}{2} - K_{\lambda})) N_{\lambda}
$.
Then
\ben \label{A-with-Gamma11}
 \langle \hat{A}, \tilde{\Gamma}_{2}^{\lambda,T}(g, g) \rangle = \langle \hat{A}, \f{1}{2}\mathcal{L}(\mu^{-\f{1}{2}}g^{2}) \rangle = C_{A} \theta u + C_{*} (\rho + \theta) u,
\\ \label{B-with-Gamma11}
 \langle \hat{B}, \tilde{\Gamma}_{2}^{\lambda,T}(g, g) \rangle = \langle \hat{B}, \f{1}{2}\mathcal{L}(\mu^{-\f{1}{2}}g^{2}) \rangle = \f{m_{2}}{3} \mathcal{B}(u),
\een
where
\beno C_{*} = \f{m_{2}^{2} - m_{0}m_{4}}{2m_{2}} , \quad \mathcal{B}(u)\colonequals  u \otimes u - \f{|u|^{2}}{3}I_{3}
\eeno
\end{lem}
\begin{proof}  The formula \eqref{Gamma-g-g-to-L-g-square-bosons} will significantly simplify the computation of $ \langle \hat{A}, \tilde{\Gamma}_{2}^{\lambda,T}(g, g) \rangle,  \langle \hat{B}, \tilde{\Gamma}_{2}^{\lambda,T}(g, g) \rangle$ for
$g \in \ker \tilde{\mathcal{L}}^{\lambda,T}$. 
By \eqref{Gamma-g-g-to-L-g-square-bosons},
since $\tilde{\mathcal{L}}^{\lambda,T}$ is self-adjoint,
we get
\beno
\langle \hat{A}, \tilde{\Gamma}_{2}^{\lambda,T}(g, g) \rangle
&=& \langle \hat{A}, \f{1}{2} \tilde{\mathcal{L}}^{\lambda,T}((1+2M_{\lambda})N_{\lambda}^{-1}g^{2}) \rangle
\\&=& \f{1}{2}\langle A, (1+2M_{\lambda})N_{\lambda}^{-1}g^{2} \rangle
\\&=&
\f{1}{2}\langle N_{\lambda}(\f{|v|^{2}}{2} - K_{A})v, (1+2M_{\lambda})N_{\lambda}(\rho + u \cdot v + \theta (\f{|v|^{2}}{2} - K_{\lambda}))^{2} \rangle
\\ &=& \langle N_{\lambda}(\f{|v|^{2}}{2} - K_{A})v, (1+2M_{\lambda})N_{\lambda} (\rho + \theta (\f{|v|^{2}}{2} - K_{\lambda}) )u \cdot v \rangle
\\ &=& \langle M_{\lambda}(1+M_{\lambda})(1+2M_{\lambda})(\f{|v|^{2}}{2} - K_{A})v, (\rho +
\theta + \theta (\f{|v|^{2}}{2} - K_{A})) u \cdot v \rangle
\\ &=&  \f{1}{3}
\langle M_{\lambda}(1+M_{\lambda})(1+2M_{\lambda}), (\f{|v|^{2}}{2} - K_{A})^{2} |v|^{2} \rangle \theta u
\\&&+ \f{1}{3}
\langle M_{\lambda}(1+M_{\lambda})(1+2M_{\lambda}), (\f{|v|^{2}}{2} - K_{A}) |v|^{2} \rangle (\rho + \theta) u
\\ &=& C_{A} \theta u + C_{*} (\rho + \theta) u,
\eeno
where we use \eqref{3-to-2} to get
 \beno \f{1}{3}
\langle M_{\lambda}(1+M_{\lambda})(1+2M_{\lambda}), (\f{|v|^{2}}{2} - K_{A}) |v|^{2} \rangle  =
\langle M_{\lambda}(1+M_{\lambda}), \f{|v|^{2}}{2} - K_{A}  \rangle = \f{m_{2}^{2} - m_{0}m_{4}}{2m_{2}} = C_{*}. \eeno

Similarly,
\beno
\langle \hat{B}, \tilde{\Gamma}_{2}^{\lambda,T}(g, g) \rangle
&=& \langle \hat{B}, \f{1}{2} \tilde{\mathcal{L}}^{\lambda,T}((1+2M_{\lambda})N_{\lambda}^{-1}g^{2}) \rangle
\\&=& \f{1}{2}\langle B, (1+2M_{\lambda})N_{\lambda}^{-1}g^{2} \rangle
\\&=&
\f{1}{2}\langle N_{\lambda}(v \otimes v - \f{|v|^{2}}{3}I_{3}), (1+2M_{\lambda})N_{\lambda}(\rho + u \cdot v + \theta (\f{|v|^{2}}{2} - K_{\lambda}))^{2} \rangle
\\ &=&
\f{1}{2}\langle M_{\lambda}(1+M_{\lambda})(1+2M_{\lambda})(v \otimes v - \f{|v|^{2}}{3}I_{3}), (u \cdot v)^{2} \rangle
\\ &=& \f{1}{15}   \langle \tilde{\mu}, |v|^{4} \rangle (u \otimes u - \f{|u|^{2}}{3}I_{3}).
\eeno
 For simplicity, let $\tilde{\mu} \colonequals  M_{\lambda}(1+M_{\lambda})(1+2M_{\lambda})$.
For the $(1,1)$ element,
\beno
\langle \tilde{\mu}, (v_{1}^{2} - \f{|v|^{2}}{3})  \sum_{i,j} u_{i}u_{j} v_{i}v_{j} \rangle  &=&
\langle \tilde{\mu}, (v_{1}^{2} - \f{|v|^{2}}{3})  v_{1}^{2} \rangle u_{1}^{2}
+ \langle \tilde{\mu}, (v_{1}^{2} - \f{|v|^{2}}{3})  v_{2}^{2} \rangle u_{2}^{2}
+ \langle \tilde{\mu}, (v_{1}^{2} - \f{|v|^{2}}{3})  v_{3}^{2} \rangle u_{3}^{2}
\\&=&  \langle \tilde{\mu}, v_{1}^{2} v_{2}^{2} \rangle (\f{4}{3} u_{1}^{2} - \f{2}{3}  u_{2}^{2} - \f{2}{3}  u_{3}^{2}) = \f{2}{15}  \langle \tilde{\mu}, |v|^{4} \rangle (u_{1}^{2} - \f{|u|^{2}}{3}).
\eeno
For the $(1,2)$ element,
\beno
\langle \tilde{\mu}, v_{1}v_{2}  \sum_{i,j} u_{i}u_{j} v_{i}v_{j} \rangle  =
2\langle \tilde{\mu}, v_{1}^{2}v_{2}^{2}  \rangle u_{1}u_{2}
= \f{2}{15}   \langle \tilde{\mu}, |v|^{4} \rangle  u_{1}u_{2}.
\eeno
As a result, we obtain \eqref{B-with-Gamma11} by observing $\langle \tilde{\mu}, |v|^{4} \rangle = 5 m_{2}$.
\end{proof}

Plugging \eqref{equation-for-Lf} into \eqref{Leray-projection-u-epsilon-2}, using \eqref{B-with-macro} and \eqref{B-with-Gamma11}, we get
\beno
\partial_{t} \mathcal{P}u^{\epsilon} - T \f{3 \nu_{\lambda, T}}{m_{2}}
 \mathcal{P}
 \nabla_{x} \cdot \mathcal{T}(u^{\epsilon})
  + T^{\f{1}{2}}
 \mathcal{P}
 \nabla_{x} \cdot  \mathcal{B}(u^{\epsilon})
 + T^{\f{1}{2}} \f{3}{m_{2}}
 \mathcal{P}
 \nabla_{x} \cdot \langle  \hat{B}, R_{\epsilon} \rangle
= 0.
\eeno
Note that
\beno
\nabla_{x} \cdot \mathcal{T}(u^{\epsilon}) = \Delta_{x}u^{\epsilon} + \f{1}{3} \nabla_{x}(\nabla_{x}\cdot u^{\epsilon}), \quad
\nabla_{x} \cdot \mathcal{B}(u^{\epsilon}) = \nabla_{x} \cdot (u^{\epsilon} \otimes u^{\epsilon}) - \f{1}{3} \nabla_{x} (|u^{\epsilon}|^{2}),
\eeno
then
\beno
\mathcal{P} \nabla_{x} \cdot \mathcal{T}(u^{\epsilon}) = \mathcal{P} \Delta_{x}u^{\epsilon} = \Delta_{x} \mathcal{P}u^{\epsilon}, \quad
\mathcal{P} \nabla_{x} \cdot  \mathcal{B}(u^{\epsilon}) = \mathcal{P}\nabla_{x} \cdot (u^{\epsilon} \otimes u^{\epsilon}).
\eeno
Then we have
\beno
\partial_{t} \mathcal{P}u^{\epsilon} - T \f{3 \nu_{\lambda, T}}{m_{2}} \Delta_{x} \mathcal{P}u^{\epsilon}
  + T^{\f{1}{2}} \mathcal{P}\nabla_{x} \cdot (u^{\epsilon} \otimes u^{\epsilon})
 + T^{\f{1}{2}} \f{3}{m_{2}}
 \mathcal{P}
 \nabla_{x} \cdot \langle  \hat{B}, R_{\epsilon} \rangle
= 0.
\eeno
With the decomposition $u^{\epsilon} = \mathcal{P}u^{\epsilon} + \mathcal{P}^{\perp}u^{\epsilon}$, defining $\mu_{\lambda, T} \colonequals T \f{3 \nu_{\lambda, T}}{m_{2}}$,
we get
\beno
\partial_{t} \mathcal{P}u^{\epsilon} -  \mu_{\lambda, T} \Delta_{x} \mathcal{P}u^{\epsilon}
  + T^{\f{1}{2}} \mathcal{P} \nabla_{x} \cdot ( \mathcal{P} u^{\epsilon} \otimes \mathcal{P} u^{\epsilon})
= R_{\epsilon, u},
\eeno
where
\beno
R_{\epsilon, u} \colonequals  - T^{\f{1}{2}} \mathcal{P} \nabla_{x} \cdot ( \mathcal{P}^{\perp} u^{\epsilon} \otimes \mathcal{P} u^{\epsilon})
- T^{\f{1}{2}} \mathcal{P} \nabla_{x} \cdot ( \mathcal{P} u^{\epsilon} \otimes \mathcal{P}^{\perp} u^{\epsilon})
- T^{\f{1}{2}} \mathcal{P} \nabla_{x} \cdot ( \mathcal{P}^{\perp} u^{\epsilon} \otimes \mathcal{P}^{\perp} u^{\epsilon})
- T^{\f{1}{2}} \f{3}{m_{2}}
 \mathcal{P}
 \nabla_{x} \cdot \langle  \hat{B}, R_{\epsilon} \rangle.
\eeno
For any $T_{*}>0$, and any test function $\psi(t,x) \in C^{1}(0,T_{*}; C_{c}^{\infty}(\mathbb{R}^{3}))$ with $\nabla_{x} \cdot \psi =0, \psi(T_{*},x) =0$, we first have
\beno
\int_{0}^{T_{*}}\int_{\mathbb{R}^{3}} \partial_{t} \mathcal{P}u^{\epsilon} \cdot \psi \mathrm{d}t \mathrm{d}x
= - \int_{\mathbb{R}^{3}} \mathcal{P}u^{\epsilon}(0,x) \cdot \psi(0,x) \mathrm{d}x - \int_{0}^{T_{*}}\int_{\mathbb{R}^{3}}  \mathcal{P}u^{\epsilon} \cdot \partial_{t} \psi \mathrm{d}t \mathrm{d}x.
\eeno
By convergence of the initial data  \eqref{initial-convergence},
\beno
\int_{\mathbb{R}^{3}} \mathcal{P}u^{\epsilon}(0,x) \cdot \psi(0,x) \mathrm{d}x \to
\int_{\mathbb{R}^{3}} \mathcal{P}u_{0} \cdot \psi(0,x) \mathrm{d}x.
\eeno
By the convergence \eqref{Pu-epsilon-to-u},
\beno
\int_{0}^{T_{*}}\int_{\mathbb{R}^{3}}  \mathcal{P}u^{\epsilon} \cdot \partial_{t} \psi \mathrm{d}t \mathrm{d}x \to
\int_{0}^{T_{*}}\int_{\mathbb{R}^{3}}  u \cdot \partial_{t} \psi \mathrm{d}t \mathrm{d}x. 
\eeno
As a result, we get 
\ben \label{partial-t-u-epsilon}
\int_{0}^{T_{*}}\int_{\mathbb{R}^{3}} \partial_{t} \mathcal{P}u^{\epsilon} \cdot \psi \mathrm{d}t \mathrm{d}x
\to - \int_{\mathbb{R}^{3}} \mathcal{P}u_{0} \cdot \psi(0,x) \mathrm{d}x - \int_{0}^{T_{*}}\int_{\mathbb{R}^{3}}  u \cdot \partial_{t} \psi \mathrm{d}t \mathrm{d}x.
\een
Recalling \eqref{rho-u-theta-convergence},  \eqref{Pu-epsilon-to-u}, and \eqref{p-perp-u-epsilon-2-0},
\ben \label{condition=on-convergence-of-Pu}
u^{\epsilon} \to u \text{ and } \mathcal{P}^{\perp}u^{\epsilon} \to 0 \text{ weakly-* in } L^{\infty}(\mathbb{R}_{+}; H^{2}_{x}),  \quad
\mathcal{P} u^{\epsilon} \to u \text{ strongly in } C(\mathbb{R}_{+}; H^{1}_{x}). \een
From this,  
we first have
\ben \label{Delta-p-u-epsilon-2-Delta-u}
\Delta_{x} \mathcal{P}u^{\epsilon} \to \Delta_{x} u \text{ in the sense of distribution}.
\een 
It is easy to see
\beno
\|\mathcal{P} \nabla_{x} \cdot ( \mathcal{P} u^{\epsilon} \otimes \mathcal{P} u^{\epsilon}) -
\mathcal{P} \nabla_{x} \cdot ( u \otimes u)\|_{L^{2}_{x}}
\lesssim \|\mathcal{P} u^{\epsilon} -u \|_{H^{1}_{x}} \|\mathcal{P} u^{\epsilon}\|_{H^{2}_{x}}
+ \|u\|_{H^{2}_{x}} \|\mathcal{P} u^{\epsilon} -u\|_{H^{1}_{x}},
\eeno
which gives
\beno
\|\mathcal{P} \nabla_{x} \cdot ( \mathcal{P} u^{\epsilon} \otimes \mathcal{P} u^{\epsilon}) -
\mathcal{P} \nabla_{x} \cdot ( u \otimes u)\|_{C(\mathbb{R}_{+}; L^{2}_{x})}
\leq \|\mathcal{P} u^{\epsilon} -u \|_{C(\mathbb{R}_{+}; H^{1}_{x})} (\|u^{\epsilon}\|_{L^{\infty}(\mathbb{R}_{+}; H^{2}_{x})} + \|u\|_{L^{\infty}(\mathbb{R}_{+}; H^{2}_{x})}).
\eeno
Then we have
\ben 
\label{P-u-epsilon-tensor-2-u-tensor}
\mathcal{P} \nabla_{x} \cdot ( \mathcal{P} u^{\epsilon} \otimes \mathcal{P} u^{\epsilon}) \to
\mathcal{P} \nabla_{x} \cdot ( u \otimes u) \text{ strongly in } C(\mathbb{R}_{+}; L^{2}_{x}).
\een

By \eqref{condition=on-convergence-of-Pu}, it is standard to derive 
\ben \label{cross-product-term}
\mathcal{P} \nabla_{x} \cdot (\mathcal{P}^{\perp} u^{\epsilon} \otimes \mathcal{P} u^{\epsilon} + \mathcal{P} u^{\epsilon} \otimes \mathcal{P}^{\perp} u^{\epsilon}) \to 0 \text{ in the distributional sense.}
\een
First note that
\beno
\|\mathcal{P} \nabla_{x} \cdot (\mathcal{P}^{\perp} u^{\epsilon} \otimes \mathcal{P} u^{\epsilon} + \mathcal{P} u^{\epsilon} \otimes \mathcal{P}^{\perp} u^{\epsilon})\|_{L^{2}_{x}}
\leq \|\mathcal{P}^{\perp} u^{\epsilon} \otimes \mathcal{P} u^{\epsilon} + \mathcal{P} u^{\epsilon} \otimes \mathcal{P}^{\perp} u^{\epsilon}\|_{H^{1}_{x}}
\lesssim \|u^{\epsilon}\|_{H^{1}_{x}} \|u^{\epsilon}\|_{H^{2}_{x}}.
\eeno
Then $\mathcal{P} \nabla_{x} \cdot (\mathcal{P}^{\perp} u^{\epsilon} \otimes \mathcal{P} u^{\epsilon} + \mathcal{P} u^{\epsilon} \otimes \mathcal{P}^{\perp} u^{\epsilon}) \in L^{\infty}(\mathbb{R}_{+}; H^{2}_{x})$.
Note that
\beno
\|\mathcal{P}^{\perp} u^{\epsilon} \otimes (\mathcal{P} u^{\epsilon} - u)\|_{C(\mathbb{R}_{+}; H^{1}_{x})} \lesssim \|\mathcal{P}^{\perp} u^{\epsilon}\|_{L^{\infty}(\mathbb{R}_{+}; H^{2}_{x})} \|\mathcal{P} u^{\epsilon} - u\|_{C(\mathbb{R}_{+}; H^{1}_{x})}
\eeno
From this together with \eqref{condition=on-convergence-of-Pu}, we have
\ben \label{part1-strongly}
\mathcal{P}^{\perp} u^{\epsilon} \otimes (\mathcal{P} u^{\epsilon} - u) \to 0 \text{ strongly in } C(\mathbb{R}_{+}; H^{1}_{x})
\een
From $u \in L^{\infty}(\mathbb{R}_{+}; H^{2}_{x})$ together with \eqref{condition=on-convergence-of-Pu}, we have
\ben \label{part2-weakly}
\mathcal{P}^{\perp} u^{\epsilon} \otimes  u \to 0 \text{ weakly-* in } L^{\infty}(\mathbb{R}_{+}; H^{1}_{x})
\een
By \eqref{part1-strongly}  and \eqref{part2-weakly}, we have
\beno
\mathcal{P}^{\perp} u^{\epsilon} \otimes \mathcal{P} u^{\epsilon} \to 0 \text{ weakly-* in } L^{\infty}(\mathbb{R}_{+}; H^{1}_{x}).
\eeno
Similarly, we can derive
\beno
\mathcal{P} u^{\epsilon} \otimes \mathcal{P}^{\perp} u^{\epsilon} \to 0 \text{ weakly-* in } L^{\infty}(\mathbb{R}_{+}; H^{1}_{x}).
\eeno
The previous two results  give \eqref{cross-product-term}.
We now derive
\ben \label{square-product-part}
\mathcal{P} \nabla_{x} \cdot ( \mathcal{P}^{\perp} u^{\epsilon} \otimes \mathcal{P}^{\perp} u^{\epsilon})
\to 0 \text{ in the distributional sense.}
\een
Applying $\mathcal{P}^{\perp}$ to \eqref{eq-u-epsilon-2}, we get
\beno
\partial_{t} \mathcal{P}^{\perp} u^{\epsilon} +  \f{T^{\f{1}{2}}}{\epsilon} \nabla_{x}  (\rho^{\epsilon} + \theta^{\epsilon}) + \f{T^{\f{1}{2}}}{\epsilon} \f{3}{m_{2}}
\mathcal{P}^{\perp} \nabla_{x} \cdot \langle B, f^{\epsilon}_{2} \rangle= 0.
\eeno
Adding \eqref{eq-rho-epsilon-2} and \eqref{eq-theta-epsilon-2}, we get
\beno
\partial_{t} (\rho^{\epsilon} + \theta^{\epsilon})+ \f{T^{\f{1}{2}}}{\epsilon} \f{2}{3} K_{A} \nabla_{x} \cdot \mathcal{P}^{\perp} u^{\epsilon} + \f{T^{\f{1}{2}}}{\epsilon}
\f{K_{A}}{C_{A}}\f{2K_{A}m_{0}-m_{2}}{m_{2}}\nabla_{x} \cdot \langle A, f^{\epsilon}_{2} \rangle = 0.
\eeno
For simplicity, let
\beno  \mathcal{F}_{\epsilon} \colonequals  T^{\f{1}{2}}
\f{K_{A}}{C_{A}}\f{2K_{A}m_{0}-m_{2}}{m_{2}}\nabla_{x} \cdot \langle  f^{\epsilon}_{2}, A \rangle,
\quad \mathcal{G}_{\epsilon} \colonequals  T^{\f{1}{2}} \f{3}{m_{2}}
\mathcal{P}^{\perp} \nabla_{x} \cdot \langle B, f^{\epsilon}_{2} \rangle.
\eeno
Then we get
\beno
\epsilon \partial_{t} [(\rho^{\epsilon} + \theta^{\epsilon})\mathcal{P}^{\perp} u^{\epsilon}] =  T^{\f{1}{2}} \f{2}{3} K_{A} (\nabla_{x} \cdot \mathcal{P}^{\perp} u^{\epsilon}) \mathcal{P}^{\perp} u^{\epsilon}
+ \mathcal{F}_{\epsilon} \mathcal{P}^{\perp} u^{\epsilon} + \f{1}{2} T^{\f{1}{2}} \nabla_{x}(\rho^{\epsilon} + \theta^{\epsilon})^{2} +
(\rho^{\epsilon} + \theta^{\epsilon}) \mathcal{G}_{\epsilon}.
\eeno
Using
\beno
\nabla_{x} \cdot ( \mathcal{P}^{\perp} u^{\epsilon} \otimes \mathcal{P}^{\perp} u^{\epsilon}) = \f{1}{2} \nabla_{x}(|\mathcal{P}^{\perp} u^{\epsilon}|^{2}) + (\nabla_{x} \cdot \mathcal{P}^{\perp} u^{\epsilon}) \mathcal{P}^{\perp} u^{\epsilon},
\eeno
and taking projection $\mathcal{P}$,
we get
\beno
T^{\f{1}{2}} \f{2}{3} K_{A} \mathcal{P} \nabla_{x} \cdot ( \mathcal{P}^{\perp} u^{\epsilon} \otimes \mathcal{P}^{\perp} u^{\epsilon}) =  \epsilon \partial_{t} \mathcal{P}[(\rho^{\epsilon} + \theta^{\epsilon})\mathcal{P}^{\perp} u^{\epsilon}] -
\mathcal{P} [\mathcal{F}_{\epsilon} \mathcal{P}^{\perp} u^{\epsilon} + (\rho^{\epsilon} + \theta^{\epsilon}) \mathcal{G}_{\epsilon}]
\eeno
By \eqref{L-2-t-HNx-L2-cb} and \eqref{L-infty-t-HNx-L2-cb}, 
\beno
\mathcal{F}_{\epsilon}, \mathcal{G}_{\epsilon} \to 0 \text{ strongly in } L^{2}(\mathbb{R}_{+}; H^{1}_{x}),
\quad (\rho^{\epsilon}, u^{\epsilon}, \theta^{\epsilon}) \in L^{\infty}(\mathbb{R}_{+}; H^{2}_{x}).
\eeno
Then we arrive at \eqref{square-product-part}.
From \eqref{cross-product-term}, \eqref{square-product-part}, and \eqref{error-term},
we obtain 
\ben \label{R-u-epsilon-2-0}
R_{\epsilon, u} \to 0 \text{ in the distributional sense.}
\een

By \eqref{partial-t-u-epsilon}, \eqref{Delta-p-u-epsilon-2-Delta-u},
\eqref{P-u-epsilon-tensor-2-u-tensor} and \eqref{R-u-epsilon-2-0}, it holds that
\ben \label{exact-meaning-of-weak-solution-u}
- \int_{\mathbb{R}^{3}} \mathcal{P}u_{0} \cdot \psi(0,x) \mathrm{d}x - \int_{0}^{T_{*}}\int_{\mathbb{R}^{3}}  u \cdot \partial_{t} \psi \mathrm{d}t \mathrm{d}x =  \int_{0}^{T_{*}}\int_{\mathbb{R}^{3}} (- T^{\f{1}{2}} \mathcal{P} (u \cdot \nabla_{x} u)
+ \mu_{\lambda, T} \Delta_{x} u) \cdot \psi \mathrm{d}t \mathrm{d}x.
\een
That is, $u$ satisfies ${\eqref{NSF-problem}}_{2}$ with the initial condition $u|_{t=0} = \mathcal{P}u_{0}$.

Plugging \eqref{equation-for-Lf} into \eqref{combination-rho-theta}, using \eqref{A-with-macro} and \eqref{A-with-Gamma11}, we get
\beno
\partial_{t} (\f{K_{\lambda}\theta^{\epsilon} - \rho^{\epsilon}}{K_{\lambda}+1}) - T \f{K_{A}}{C_{A}(K_{\lambda}+1)} (\kappa_{1,\lambda, T} \Delta_{x} \theta^{\epsilon} + \kappa_{2,\lambda, T} \Delta_{x} (\rho^{\epsilon}+\theta^{\epsilon}))
  \\+ T^{\f{1}{2}} \f{K_{A}}{C_{A}(K_{\lambda}+1)} \nabla_{x} \cdot(
C_{A} \theta^{\epsilon} u^{\epsilon} + C_{*} (\rho^{\epsilon} + \theta^{\epsilon}) u^{\epsilon})
 + T^{\f{1}{2}}
\f{K_{A}}{C_{A}(K_{\lambda}+1)}
 \nabla_{x} \cdot \langle  \hat{A}, R_{\epsilon} \rangle
= 0
\eeno
Recalling \eqref{defintion-of-KA}, defining
\ben
\kappa_{\lambda, T} \colonequals   \f{T \kappa_{1,\lambda, T}}{C_{A}},
\een
we have
\beno
\partial_{t} (\f{K_{\lambda}\theta^{\epsilon} - \rho^{\epsilon}}{K_{A}})
+ T^{\f{1}{2}} u^{\epsilon} \cdot \nabla_{x}\theta^{\epsilon}
= \kappa_{\lambda, T} \Delta_{x} \theta^{\epsilon} + R_{\epsilon,\theta},
\eeno
where
\beno
R_{\epsilon,\theta}\colonequals  \f{T \kappa_{2,\lambda, T}}{C_{A}} \Delta_{x} (\rho^{\epsilon}+\theta^{\epsilon}) -  T^{\f{1}{2}} \theta^{\epsilon} \nabla_{x} \cdot u^{\epsilon} - T^{\f{1}{2}} \f{C_{*} }{C_{A}}  \nabla_{x} \cdot((\rho^{\epsilon} + \theta^{\epsilon}) u^{\epsilon}) -
\f{T^{\f{1}{2}}}{C_{A}}
 \nabla_{x} \cdot \langle  \hat{A}, R_{\epsilon} \rangle.
\eeno
For any $T_{*}>0$, and any test function $\psi(t,x) \in C^{1}(0,T_{*}; C_{c}^{\infty}(\mathbb{R}^{3}))$ with $\psi(T_{*},x) =0$,
then similarly to \eqref{partial-t-u-epsilon}, we have
\ben \label{partial-t-theta-epsilon}
\int_{0}^{T_{*}}\int_{\mathbb{R}^{3}} \partial_{t} (\f{K_{\lambda}\theta^{\epsilon} - \rho^{\epsilon}}{K_{\lambda}+1}) \psi \mathrm{d}t \mathrm{d}x
\to - \int_{\mathbb{R}^{3}} (\f{K_{\lambda}\theta_{0} - \rho_{0}}{K_{\lambda}+1}) \psi(0,x) \mathrm{d}x - \int_{0}^{T_{*}}\int_{\mathbb{R}^{3}}  \theta \partial_{t} \psi \mathrm{d}t \mathrm{d}x.
\een
Similarly to \eqref{Delta-p-u-epsilon-2-Delta-u}, \eqref{P-u-epsilon-tensor-2-u-tensor} and \eqref{R-u-epsilon-2-0}, we have
\ben \label{Delta-theta-epsilon-2-Delta-theta}
\Delta_{x} \theta^{\epsilon} \to \Delta_{x} \theta \text{ in distributional sense},
\\
\label{u-cdot-nabla-theta}
u^{\epsilon} \cdot \nabla_{x}\theta^{\epsilon} \to u \cdot \nabla_{x} \theta \text{ strongly in } C(\mathbb{R}_{+}; L^{2}_{x}),
\\ \label{R-theta-epsilon-2-0}
R_{\epsilon, \theta} \to 0 \text{ in the distributional sense.}
\een
By \eqref{partial-t-theta-epsilon}, \eqref{Delta-theta-epsilon-2-Delta-theta}, \eqref{u-cdot-nabla-theta} and \eqref{R-theta-epsilon-2-0}, we have
\ben \label{exact-meaning-of-weak-solution-theta}
- \int_{\mathbb{R}^{3}} (\f{K_{\lambda}\theta_{0} - \rho_{0}}{K_{\lambda}+1}) \psi(0,x) \mathrm{d}x - \int_{0}^{T_{*}}\int_{\mathbb{R}^{3}}  \theta \partial_{t} \psi \mathrm{d}t \mathrm{d}x =  \int_{0}^{T_{*}}\int_{\mathbb{R}^{3}} (- T^{\f{1}{2}} u \cdot \nabla_{x}\theta
+ \kappa_{\lambda, T} \Delta_{x} \theta) \psi \mathrm{d}t \mathrm{d}x.
\een
That is, $\theta$ satisfies ${\eqref{NSF-problem}}_{3}$ with the initial condition $\theta|_{t=0} = \f{K_{\lambda}\theta_{0} - \rho_{0}}{K_{\lambda}+1}$.

We already obtain all the results in Theorem \ref{thm-hydrodynamical-limit}.
\end{proof}

\section{Appendix} \label{appendix}
In this appendix, we see why $\hat{A}, \hat{B}$ must be of the form in \eqref{inverse-function-form-lambda-T} in order to satisfy $\tilde{\mathcal{L}}^{\lambda,T} (\hat{A}) = A, \tilde{\mathcal{L}}^{\lambda,T} (\hat{B}) = B$. 

Let $O_{3}$ be the isometry on $\mathbb{R}^{3}$.
Given $R \in O_{3}$, we
define the rotation operator $T_{R}$ on functions $f: \mathbb{R}^{3} \to \mathbb{R}$ by
\beno
T_{R} f (v) = f(R v).
\eeno
Thanks to the intrinsic structure of $\tilde{\mathcal{L}}^{\lambda,T}$, it commutes with rotation operator.
\begin{lem} \label{commutation}[Commutative property] Let $R \in O_{3}$, then
\beno
T_{R} \tilde{\mathcal{L}}^{\lambda,T}  =  \tilde{\mathcal{L}}^{\lambda,T} T_{R}.
\eeno
\end{lem}
As a result of Lemma \ref{commutation}, for any $R \in O_{3}$,
\ben \label{TR-preserves-orthogonal}
f \in (\ker\tilde{\mathcal{L}}^{\lambda,T})^{\perp} \quad \Leftrightarrow \quad T_{R} f \in (\ker\tilde{\mathcal{L}}^{\lambda,T})^{\perp}.
\een

For any $g \in (\ker\tilde{\mathcal{L}}^{\lambda,T})^{\perp}$, the following problem has at most a unique solution
\ben \label{unique-solution}
\tilde{\mathcal{L}}^{\lambda,T} f = g, \quad f \in (\ker\tilde{\mathcal{L}}^{\lambda,T})^{\perp}.
\een

We recall two elementary results in linear algebra.

\begin{lem} \label{invariance-gives-vector-form} Let $f: \mathbb{R}^{3} \to \mathbb{R}^{3}$. If for any
	$R \in O_{3}$,
	\beno
	f \circ R = R \circ f.
	\eeno
	Then there exists a radial function $\alpha = \alpha (|v|): \mathbb{R}^{3} \to \mathbb{R}$ such that
	\beno
	f(v) = \alpha (|v|) v.
	\eeno
\end{lem}

\begin{lem} \label{invariance-gives-matrix-form} Let $f: \mathbb{R}^{3} \to \mathbb{R}^{3 \times 3}$. Suppose for any
	$R \in O_{3}$,
	\beno
	f \circ R = R f R^{-1}
	\eeno
	as functions $\mathbb{R}^{3} \to \mathbb{R}^{3 \times 3}$. Here the right-hand side is interpreted as matrix multiplication.
	Moreover, for any $v \in \mathbb{R}^{3}$,
	\ben \label{symmetric-traceless-pointwise}
	\text{ the matrix } f(v) \text{ is symmetric and traceless}.
	\een
	Then there exists a radial function $\beta = \beta (|v|) : \mathbb{R}^{3} \to \mathbb{R}$ such that
	\beno
	f(v) = \beta (|v|) (v  \otimes v - \f{|v|^{2}}{3}I_{3}).
	\eeno
\end{lem}

Now we are ready to establish the following result.
\begin{thm} \label{form-of-inverse} Recall \eqref{definition-of-vector-and-matrix-bos} for $A,B$.
	Let $\hat{A}, \hat{B} \in (\ker\tilde{\mathcal{L}}^{\lambda,T})^{\perp}$ be solutions to
	\ben \label{inverse-function-definition}
	\tilde{\mathcal{L}}^{\lambda,T}(\hat{A}) = A, \quad \tilde{\mathcal{L}}^{\lambda,T}(\hat{B}) = B.
	\een
	The functions $\hat{A}, \hat{B}$ of $A, B$
 must be of the following form
\ben \label{inverse-function-form}
\hat{A} = \alpha_{\lambda,T}(|v|) A, \quad \hat{B} = \beta_{\lambda,T}(|v|) B,
\een
for some radial functions $\alpha_{\lambda,T}(|v|), \beta_{\lambda,T}(|v|)$ depending on the operator $\tilde{\mathcal{L}}^{\lambda,T}$.
\end{thm}
\begin{proof} We claim for any $R \in O_{3}$,
\ben \label{commute-with-rotation}
\hat{A} \circ R = R \circ \hat{A},
\een
as functions $\mathbb{R}^{3} \to \mathbb{R}^{3}$. Then by Lemma \ref{invariance-gives-vector-form}, we conclude the existence of $\alpha_{\lambda,T}(|v|)$.

For any function $f : \mathbb{R}^{3} \to \mathbb{R}^{3}$, $f \circ R = T_{R} f$. Here $R \circ f$ is a just linear combination of $f$.
Note that
\beno
\tilde{\mathcal{L}}^{\lambda,T} (\hat{A} \circ R) = \tilde{\mathcal{L}}^{\lambda,T} T_{R} \hat{A} = T_{R} \tilde{\mathcal{L}}^{\lambda,T}  \hat{A} = T_{R} A = A \circ R,
\\
\tilde{\mathcal{L}}^{\lambda,T} (R \circ \hat{A}) = R \circ \tilde{\mathcal{L}}^{\lambda,T} \hat{A} = R \circ A.
\eeno
The first line uses Lemma \ref{commutation}, while the second uses the linearity of $\tilde{\mathcal{L}}^{\lambda,T}$.
It is easy to check
\beno
R \circ A = A \circ R, \quad \hat{A} \circ R \in (\ker\tilde{\mathcal{L}}^{\lambda,T})^{\perp}, \quad R \circ \hat{A}  \in (\ker\tilde{\mathcal{L}}^{\lambda,T})^{\perp}.
\eeno
The first one is obvious.
The second one is given by \eqref{TR-preserves-orthogonal}. For the third one, note that
 each element of $R \circ \hat{A}$ is just a linear combination of $\hat{A}_{i}$.
By uniqueness of problem \eqref{unique-solution}, we get \eqref{commute-with-rotation}.

We claim for any $R \in O_{3}, v \in \mathbb{R}^{3}$,
\ben \label{commute-with-rotation-matrix-type}
\hat{B} \circ R = R \hat{B} R^{-1}
\een
as functions $\mathbb{R}^{3} \to \mathbb{R}^{3 \times 3}$.
Moreover, for any $v \in \mathbb{R}^{3}$,
\ben \label{symmetric-traceless}
\hat{B}(v) \text{ is symmetric and traceless}.
\een
Then by Lemma \ref{invariance-gives-matrix-form}, we conclude the existence of $\beta_{\lambda,T}(|v|)$.

For any function $f : \mathbb{R}^{3} \to \mathbb{R}^{3 \times 3}$, $f \circ R = T_{R} f$.
Here $R \circ f$ is just a linear combination of $f$.
Note that
\beno
\tilde{\mathcal{L}}^{\lambda,T} (\hat{B} \circ R) = \tilde{\mathcal{L}}^{\lambda,T} T_{R} \hat{B} = T_{R} \tilde{\mathcal{L}}^{\lambda,T}  \hat{B} = T_{R} B = B \circ R,
\\
\tilde{\mathcal{L}}^{\lambda,T} (R \hat{B} R^{-1})= R (\tilde{\mathcal{L}}^{\lambda,T} \hat{B}) R^{-1}= R B R^{-1}.
\eeno
It is easy to check
\beno
B \circ R = R B R^{-1}, \quad \hat{B} \circ R \in (\ker\tilde{\mathcal{L}}^{\lambda,T})^{\perp}, \quad R \hat{B} R^{-1}  \in (\ker\tilde{\mathcal{L}}^{\lambda,T})^{\perp}.
\eeno
For the first one, note that
\beno
(B \circ R)(v) = B (Rv) = Rv  \otimes Rv - \f{|v|^{2}}{3}I_{3} = Rv  \otimes v R^{-1} - \f{|v|^{2}}{3} R I_{3} R^{-1} = R B(v) R^{-1}.
\eeno
The second one is given by \eqref{TR-preserves-orthogonal}.
For the third one, each element of $R \hat{B} R^{-1}$ is just linear combination of $\hat{B}_{ij}$.
By uniqueness of problem \eqref{unique-solution}, we get \eqref{commute-with-rotation-matrix-type}.

As $B(v)$  is symmetric and traceless, we have
\beno
\tilde{\mathcal{L}}^{\lambda,T} (\hat{B} - \hat{B}^{\mathrm{T}}) = B - B^{\mathrm{T}} = 0, \quad \hat{B} - \hat{B}^{\mathrm{T}} \in (\ker\tilde{\mathcal{L}}^{\lambda,T})^{\perp},
\\
\tilde{\mathcal{L}}^{\lambda,T} ( \mathrm{Tr}(\hat{B})) = \mathrm{Tr}(B) = 0, \quad \mathrm{Tr}(\hat{B}) \in (\ker\tilde{\mathcal{L}}^{\lambda,T})^{\perp}.
\eeno
Uniqueness of problem \eqref{unique-solution} gives
\beno
\hat{B} - \hat{B}^{\mathrm{T}}= 0, \quad \mathrm{Tr}(\hat{B}) = 0.
\eeno
That is, we get \eqref{symmetric-traceless}.
\end{proof}

 {\bf Acknowledgments.} This work was partially supported by National Key Research and Development Program of China under the grant 2021YFA1002100. 
 Ling-Bing He was supported by NSF of China  
 under the grant 12141102. Ning Jiang was supported by NSF of China under the grants 11971360 and 11731008, and also supported by the Strategic Priority Research Program of Chinese Academy of Sciences under the grant XDA25010404.
 Yu-Long Zhou was  partially supported by NSF of China under the grant 12001552,  Science and Technology Projects in Guangzhou under the grant 202201011144, and
 Youth Talent Support Program of Guangdong Provincial Association for Science and Technology under the grant SKXRC202311.

\bibliographystyle{siam}
\bibliography{bbe2nsf-7}

\end{document}